\documentclass[nonatbib,11pt]{elsarticle}

\makeatletter
\let\c@author\relax

\makeatother

\journal{{\color{black}EJOR}}

\usepackage{pifont}

\newcommand{\cmark}{\ding{51}} 
\newcommand{\xmark}{\ding{55}} 

\makeatletter
\def\ps@pprintTitle{%
  \let\@oddhead\@empty
  \let\@evenhead\@empty
  \let\@oddfoot\@empty
  \let\@evenfoot\@oddfoot
}
\makeatother
\usepackage[usenames,dvipsnames]{color}

\usepackage[utf8]{inputenc}
\usepackage[english]{babel}
\usepackage{amsthm}
\usepackage{thmtools}
\usepackage{thm-restate}
\usepackage{tabularx}
\usepackage{graphicx}
\usepackage{adjustbox}
\usepackage{hyperref}
\usepackage{blindtext}
\usepackage{amsmath,amssymb,amsfonts,amsthm}
\usepackage{cleveref}
\usepackage{fancyhdr}
\usepackage[table]{xcolor}
\usepackage{tikz}
\usetikzlibrary{shapes.geometric, positioning}
\usepackage{pgfplots}
\usepackage{siunitx}
\usepackage[linesnumbered,ruled,vlined]{algorithm2e}

\usepackage[short]{optidef}
\usepackage{float}
\graphicspath{{figures/}}
\usepackage{amssymb}
\usepackage{color}
\usepackage{mathtools}
\usepackage{subcaption}
\usepackage{etex}
\usepackage{lipsum} 

\usepackage{multirow}
\usepackage{longtable}
\usepackage{supertabular}

\usetikzlibrary{positioning}

\usepackage{hyperref}

\usepackage[style=apa, backend=biber, hyperref=true]{biblatex}
\usepackage{graphics} 
\usepackage{setspace} 
\usepackage{graphicx} 
\usepackage{lscape}
\usepackage{makeidx} 
\usepackage{pdflscape} 
\usepackage{amsmath,amssymb,amsfonts} 
\usepackage{mathptmx} 
\usepackage[final]{pdfpages} 
\usepackage{calc} 
\usepackage{fancyhdr} 
\usepackage{multicol} 
\usepackage[T1]{fontenc}
\usepackage[utf8]{inputenc} 
\usepackage[english]{babel}
\usepackage{sectsty} 
\usepackage{titlesec} 
\usepackage{threeparttable, supertabular} 
\usepackage{endnotes} 
\usepackage{rotating, rotfloat} 
\usepackage{fullpage} 
\usepackage[usenames,dvipsnames]{color} 
\usepackage{colortbl} 
\usepackage{eurosym}
\usepackage[bottom, flushmargin]{footmisc} 
\usepackage{booktabs} 
\usepackage{gensymb}
\usepackage{tabularx} 
\usepackage[format=plain, textfont=normal, justification=justified, singlelinecheck=false, font=normal]{caption}

\allowdisplaybreaks

\begin{document}

\begin{frontmatter}

\title{A Dynamic Strategic Plan for the Transition to a Clean Bus Fleet using Multi-Stage Stochastic Programming with a Case Study in Istanbul}

\author[1]{Neman Karimi}

\ead{neman.karimi@sabanciuniv.edu}

\address[1]{Faculty of Engineering and Natural Sciences, Sabanc{\i} University, Istanbul, Turkey}

\author[1]{Burak Kocuk\corref{cor1}}
\ead{burak.kocuk@sabanciuniv.edu}

\author[1,2]{Tu\u{g}\c{c}e Y\"{u}ksel}
\ead{tugce.yuksel@sabanciuniv.edu}

\address[2]{Smart Mobility and Logistics Lab, Sabanc{\i} University, Istanbul, Turkey}

\cortext[cor1]{Corresponding author}

\begin{abstract}
In recent years, the transition to clean bus fleets has accelerated. Although this transition might bring environmental and economic benefits, it requires a long-term strategic plan due to the large investment costs involved. This paper proposes a multi-stage stochastic program to optimize strategic plans for the clean bus fleet transition that explicitly considers the uncertainty scenarios in the cost and efficiency improvements of clean buses. Our optimization model minimizes the total expected cost subject to emission targets, budget restrictions and several other operational considerations. We propose a new forecasting approach that captures the correlation between these improvements to obtain realistic future pathways for Battery Electric Buses (BEBs) and Hydrogen Fuel Cell Buses (HFCBs), which are then given to the multi-stage stochastic program as scenarios. We also utilize a physics-based model for BEBs to accurately capture their energy consumption and recharging needs. As a case study, we focus on the complex public bus network of Istanbul, which aims to transition to a clean bus fleet by 2050. Utilizing real datasets, we solve a five-stage stochastic program spanning a 25-year planning horizon that involves 256 scenarios to obtain dynamic strategic plans that can be used by the policy makers. Our results suggest that BEBs are more advantageous than HFCBs, even in slow BEB but fast HFCB development scenarios. We also conduct several sensitivity analyses to understand the effects of the intermediate emission targets, budget limitations and energy prices.
\end{abstract}

\begin{keyword}
\texttt{OR in Energy, Bus fleet transition, Zero-emission vehicles, Multi-stage stochastic programming}
\end{keyword}

\end{frontmatter}

\section{Introduction}\label{Introduction}

Transport-related CO$_2$
emissions account for over 21\% of global emissions, with the majority coming from internal combustion engine vehicles used in road transport \parencite{IEAWorldEnergyOutlook2023}. Due to the significant impact of these emissions on climate change, ambitious targets have been set to transition the transport sector from fossil fuels to sustainable energy sources. For instance, the United States aims to eliminate all emissions in the transport sector by 2050, while the European Union pledges to reduce transport-related emissions by 90\% by the same year \parencite{USgoal, EUgoal}. Although the adoption of clean energy alternatives such as electric and fuel cell vehicles has increased in recent years, a lot more effort is still needed to keep the long-term net-zero emission goals within reach. This is particularly important for medium and heavy-duty vehicles such as trucks and buses, which need accelerated transitions \parencite{IEA2023roadmap}. Despite making up only 8\% of all vehicles (excluding two- and three-wheelers), trucks and buses account for over 35\% of road transport emissions \parencite{IEATrucksAndBuses}. To accelerate the shift of these vehicles to clean technology alternatives, 33 countries have joined the Global Memorandum of Understanding on Zero-Emission Medium- and Heavy-Duty Vehicles as of 2023, committing to 30\% zero-emission truck and bus sales by 2030 and 100\% by 2040 \parencite{IEAEVoutlook2024}.

Such policies and pledges necessitate strategic planning to minimize the total costs of transitioning heavy-duty vehicle fleets, such as city buses. Due to the currently high investment costs of zero-emission buses—including their purchase and the costs related to building the required infrastructure— carefully constructed strategic plans could save millions of dollars. These plans should focus on identifying the most cost-effective mix of bus technologies to purchase, the most favorable infrastructure technology to utilize (such as overnight charging vs. fast charging for electric buses), and the optimal timing for bus replacements. To this end, case studies have been conducted for several cities worldwide, including cities in 
China \parencite{zhang2022}, Singapore \parencite{zhou2023}, Germany \parencite{dirks2022},
France \parencite{pelletier2019}, the US \parencite{islam2019},  
and Austria \parencite{friess2024}, to optimize the bus fleet transition to clean energy buses. 

Most of the relevant studies focused on electric buses (EB) as the clean energy option. With today's technology, EBs are widely considered to be the most cost-effective zero-emission alternative. However, their integration into fleets is complicated by a number of barriers, the most significant being their limited driving range. 
Currently, battery electric buses (BEBs), which rely only on electricity stored in on-board batteries, are the most common type of EBs. There are different options available for BEBs in terms of battery capacity which can affect operational planning during their use. BEBs equipped with large and heavy batteries might be recharged only at night due to the long recharging time, or those that have smaller batteries can utilize fast charging stations and possibly recharge multiple times within the day to ensure demand satisfaction. 
Assuming that existing bus dispatch timetables are maintained, it is crucial to ensure that buses relying on overnight charging (ONC) have sufficient battery capacity to complete the scheduled trips.
Moreover, fast-charging (FC) buses require an additional consideration in recharge scheduling, which depends on the frequency of assigned trips, availability of charging stations and the time needed for recharging. To this end, different studies have attempted to optimize the decisions related to the investment and placement of charger technology, investments on electric buses with different battery capacities and recharge scheduling~\parencite{perumal2022}.

Hydrogen fuel cell buses (HFCBs) are also gaining attention in the literature. Although they may not yet be competitive with EBs, efforts to reduce their purchase and operational costs could make them promising options for new city buses in near future. These buses are similar to diesel buses in terms of their operation, as they usually require one daily refueling, which takes less than 10 minutes \parencite{anandarajah2013}.

Many efforts are being made to improve both BEBs and HFCBs in terms of reducing costs and increasing efficiency. For BEBs, the anticipated increase in energy density of batteries will improve their driving range, allowing them to operate on lines that were previously impractical or not cost-effective. Improvements in HFCBs are also expected to increase fuel cell efficiency, thereby decreasing fuel consumption and lowering operational costs. However, these improvements in cost and efficiency are often overlooked in strategic transition plans for bus fleets. Anticipated cost reductions are only rarely considered  as yearly reductions in the purchase cost of BEBs and in a deterministic manner. Efficiency-related improvements and their implications are even more scarce in the literature. In~\textcite{he2023}, reductions in energy consumption of newly purchased buses and increased recharging efficiency were included, but again as deterministic parameters. 

To the best of our knowledge, no study has included long-term technological advances in BEBs and HFCBs, which are uncertain by nature, affecting their costs and efficiency as stochastic scenarios. To fill this gap, we present an optimization model to strategically plan the bus fleet transition to clean energy technologies, treating cost- and efficiency-related technological advances in batteries and fuel cells as scenarios in a large-scale multi-stage stochastic programming model. After forecasting the advances and clustering them into scenarios in order to capture their inherent correlation, we test our model on the large bus fleet of Istanbul Municipality, which aims to transition to a clean bus fleet by 2050 as part of its Sustainable Urban Mobility Plan~\parencite{IstanbulSUMP2022},
with more than 6,500 buses operating on over {830} lines. We also perform sensitivity analysis to provide insights for bus fleet owners and managers.

\subsection{Literature Review}

The literature on strategic planning for transitioning to clean energy buses can be categorized into two main approaches: ‘one-step transition’ and ‘multi-step transition.’ In the one-step transition approach, strategic decisions, including investments in fleet and infrastructure, are made at a single point in time. In contrast, the multi-step transition approach involves making decisions at multiple points throughout the planning period, usually because of budget constraints or long-term emission reduction goals. This gradual transition approach is better suited for taking into account evolving technologies and changing parameters; therefore, we adopt this approach in the current paper. In what follows, we first mention some one-step transition studies before delving into the literature on multi-step transition.

Existing studies with a one-step transition approach are mostly focused on electric buses and their necessary infrastructure. Different charging technologies, including recharging at stations, lane charging, and battery swap technologies were evaluated against each other in some studies \parencite{bi2017, chen2018}. However, the majority of studies focused on charging stations and involved the deployment of fast-charging and overnight depot charging stations in their strategic plans. This includes optimal placement of fast-charging stations \parencite{xylia2017}, decisions about battery capacity \parencite{kunith2017}, and investigating the effect of factors like demand uncertainties \parencite{an2020} and energy consumption uncertainty \parencite{benoliel2021} on such decisions. Some researchers integrated more operational-level problems such as electric bus scheduling and charge scheduling, with strategic decisions including infrastructure planning, battery sizing and fleet sizing, see, e.g., \textcite{rogge2018, yildirim2021, he2022,  shehabeldeen2024}.

On the other hand, most studies with a multi-step approach developed optimization models similar to those for the parallel equipment replacement problem. In these models, purchasing and salvaging decisions are made in each period to minimize total costs while ensuring demand is met throughout the planning horizon. 
\textcite{keles2004} were among the first to apply such models to bus fleets, including different competing technologies. Later studies incorporated emission-related costs \parencite{feng2014} or constraints \parencite{emiliano2020}, along with decisions for bus-to-task assignments \parencite{stasko2010}.
\textcite{islam2019} included BEBs as clean technology options and accounted for their charger costs. An integer linear program developed by \textcite{pelletier2019} aimed to minimize total fleet management costs, incorporating various charging technologies for electric buses, midlife costs for BEB batteries, and bus-to-route assignments. \textcite{tang2021} addressed diesel-electric replacement ratios to manage potential fleet size increases due to range anxiety associated with BEBs.

Some recent multi-step transition studies integrated strategic decisions with various operational-level decisions. \textcite{dirks2022} addressed charging station deployment and battery sizing decisions, along with operational decisions such as assigning electric buses to routes and tracking battery state of charge. \textcite{zhang2022} incorporated seasonal variations in electric bus consumption into their model. \textcite{zhou2023} used an aggregated demand approach and considered external costs, including climate and health impacts, battery replacement costs, and evaluated different types of electric and hybrid buses.
\textcite{li2022}, focused on determining charger locations, planning bus route electrification, and assigning buses to charging stations in a Build-Operate-Transfer setting.
\textcite{he2023} accounted for future technological advances such as reduced battery and charger costs, and improved charging efficiency in a deterministic setting. This study included siting and sizing of fast-charging stations along with recharge schedules for new BEBs.

In Table \ref{tab:literature}, we compare various studies on multi-step bus fleet transition and highlighting their key aspects.
To the best of our knowledge, no study has incorporated the uncertainty of future ZEB improvements into strategic planning for bus fleet transitions. As   both BEBs and HFCBs are still developing,  incorporating these uncertainties into planning ensures that strategies remain adaptable and cost-effective as advancements occur. The current absence of HFCBs in many strategic plans is largely due to their higher costs and lower efficiency compared to BEBs, an aspect that  should be  reconsidered as  the   future technological improvements   could make HFCBs more competitive.
Additionally, no  study has considered the impact of route-specific energy consumption of BEBs, seasonal variations in their consumption rate, and diesel-electric replacement ratios
together in strategic planning. To fill these gaps in the literature, we present a multi-stage stochastic program that incorporates technological developments while also considering detailed operational considerations, providing a comprehensive and adaptive solution for bus fleet transition planning.
\setlength{\tabcolsep}{.1pt}
\begin{table}[h!]
    \centering
    \caption{Studies on multi-step bus fleet transition.
    }
    \label{tab:literature}
    \scriptsize
    \begin{tabular}{|c|c|c|c|c|c|c|c|c|c|c|c|c|}
        \hline
        & & \multicolumn{3} {c|} {Decisions} & \multicolumn{4} {c|} { EB Operational Feasibility} & \multicolumn{3} {c|} { Technological Change} \\
        \hline
        Reference & Bus & Fleet & {Charger} & Task & Calculating & Recharge  & Seasonal & Diesel-electric & Cost & Efficiency & Stochasticity 
        \\
        & Technologies & Invest. & Invest. & Assign. &  Electricity & Scheduling  & Effect & Replacement  & & & 
        \\
        & & & & & Consumption & &  & Ratio & & & \\
        \hline
        \textcite{stasko2010} & DB & \cmark & - & \cmark & - & - & - & - & \xmark & \xmark & \xmark \\
        \textcite{feng2014} &  DB/HEB & \cmark & - & \xmark & - & - & - & - & \cmark & \xmark & \xmark \\
        \textcite{islam2019} & DB/HEB/BEB & \cmark & \cmark & \xmark & \xmark & \xmark & \xmark & \xmark & \cmark & \xmark & \xmark \\
        \textcite{pelletier2019} & DB/CNG/BEB & \cmark & \cmark & \cmark & \xmark & \xmark & \xmark & \xmark & \cmark & \xmark & \xmark \\
        \textcite{emiliano2020} & DB & \cmark & - & \xmark & - & - & - & - & \xmark & \xmark & \xmark \\
        \textcite{tang2021} & DB/HEB/BEB  & \cmark & \cmark & \xmark & \xmark & \xmark & \xmark & \cmark & \cmark & \cmark & \xmark \\
        \textcite{li2022} & BEB  & \cmark & \cmark & \cmark & \xmark & \cmark & \xmark & - & \cmark & \xmark & \xmark \\
        \textcite{dirks2022} & DB/BEB & \cmark & \cmark & \cmark & \xmark & \cmark & \xmark & - & \cmark & \xmark & \xmark \\
        \textcite{zhang2022} & DB/BEB & \cmark &  \cmark & \cmark & \cmark & \cmark & \cmark & - & \cmark & \xmark & \xmark  \\
        \textcite{zhou2023} & DB/HEB/BEB & \cmark & \cmark & \xmark & \xmark & \xmark & \xmark & \xmark & \cmark & \xmark  & \xmark\\
        \textcite{he2023} & BEB & \cmark & \cmark & \cmark & \cmark & \cmark & \xmark & - & \cmark & \cmark & \xmark \\
        \hline
        This Paper &  DB/BEB/HFCB & \cmark & \cmark & \cmark & \cmark & \cmark & \cmark & \cmark & \cmark & \cmark & \cmark \\
        \hline
    \end{tabular}
    
       DB: Diesel, BEB: Battery Electric, HFCB: Hydrogen Fuel Cell, CNG: Compressed Natural Gas, HEB: Hybrid Electric  {\color{white} dummydummydumm}
\end{table}

\subsection{Our Approach and Contributions}

In this paper, we develop a multi-stage stochastic program, in which the future costs and efficiencies of competing clean technologies are represented as nodes in a scenario tree. Decisions regarding purchasing and salvaging the buses and assigning them to routes will be made on a yearly basis, and a dynamic transition plan will guide the fleet managers to make optimal decisions in each scenario throughout the planning horizon. We then test our model on a large bus fleet in Istanbul, including BEBs utilizing ONC or FC technology, as well as HFCBs. To account for operational feasibility of replacing diesel buses with electric options, we pre-calculate the energy consumed by each type and version of BEB on a given route. Using a heuristic method, we estimate the diesel-electric replacement ratio for each electric bus, on each route, and under each scenario. We provide an extensive case study involving real data collected from various sources, and conduct sensitivity analyses regarding the underlying stochastic process, the emission targets, budget limitations and energy prices.
Our main contributions to the literature are listed as below.
\begin{itemize}
  \setlength\itemsep{-0.5em}
    \item We introduce a multi-stage stochastic program for the fleet replacement problem, and investigate the effect of uncertain technological advances on the optimal long-term transition plan of vehicle fleets.
    \item 
    We forecast the technological improvements of both BEBs and HFCBs in terms of  cost reductions and efficiency improvements, and cluster them into a number of scenarios.
    \item We test our model on the large bus network of Istanbul for a 25-year planning horizon with five stages, involving details such as the energy consumption calculation of several versions of BEBs on each bus line,  energy consumption changes with seasonal variations, and ensuring operational feasibility by finding the diesel-electric replacement ratios and recharge scheduling for FC electric buses.
\end{itemize}
We note that although we focus on a bus fleet transition problem in Istanbul, the idea of modeling technology advances as scenarios and incorporating them into a multi-stage stochastic program is applicable to other energy transition problems, emphasizing both the novelty and the generalization potential of our work.

The rest of this paper is organized as follows: Section~\ref{sec:method} provides a detailed explanation of the problem and our methodology. Section~\ref{sec:data} describes the data requirements and how we obtain the parameters  for our case study. Section~\ref{sec:case} presents the case study and sensitivity analysis. Finally, Section~\ref{sec:conclusion} concludes the paper.

\section{Methodology}
\label{sec:method}

This section provides a detailed explanation of our methodological approach. We begin with a conceptual description of the problem  in Section~\ref{sec:method-description-framework}, followed by the mathematical formulation in Section~\ref{sec:method-formulation}. Next, in Section~\ref{sec:method-energy}, we explain the method for estimating the amount of energy consumed by the BEBs with different configurations, when assigned to each scheduled trip. Finally, in Section~\ref{sec:method-algorithm}, we describe our \textit{Bus-to-Route Simulator}, which is run to obtain key parameters of our optimization model   if a specific bus is assigned to a route.

\subsection{Problem Description}
\label{sec:method-description-framework}

Bus fleet operators face the challenge of determining the optimal long-term transition plan from DBs to zero-emission buses (ZEBs) such as BEBs and HFCBs. The goal is to minimize the total costs of managing the fleet over a planning horizon subject to emission  targets, budget limitations and other operational requirements. Let us denote \(\mathcal{T} = \{0, 1, \ldots, T\}\) as the time periods in the planning horizon, where decisions on  purchasing, salvaging, and assigning buses to routes are made at the beginning of each time period. Time period 0 is used to initiate the existing fleet, and no actual decisions are made in this time period. From time period 1 onwards, bus fleet operators must decide for each bus whether to keep it (if it has not reached its economic lifetime) or to replace it with one of the available bus technologies in the market.

Accounting for the technological advancements of ZEBs is essential for planning a cost-effective transition as these advancements can significantly improve performance and reduce costs. However, the uncertain nature of these advancements makes it challenging to determine the optimal timing for ZEB investments. To tackle this complex task, we introduce a multi-stage stochastic program with $S$ stages, where each stage consists of a collection of time periods. We follow the classical scenario tree based representation~\parencite{shapiro2021} where nodes within each stage represent possible states of all technologies over these periods in terms of cost and efficiency. Although our formulation below works for any given scenario tree $\mathcal{N}$, we specify how we construct it according to our case study in Section~\ref{sec:scenario-gen}.

Before formally providing our formulation, we list our  
assumptions below:
\begin{itemize}
  \setlength\itemsep{-0.5em}
    \item Each time period represents a year, as bus fleet operators typically make decisions based on the yearly budget. We further divide each year into subperiods (seasons) to account for potential differences in route demand and variations in BEB and HFCB consumption rates due to seasonal factors.

    \item The time periods within each stage share the same stochastic characteristics.

    \item Different versions may exist for each type of technology, reflecting variations such as the brand and model of the bus.
    For example, one version of a BEB might be a 12-meter model with a 280 \unit{\kWh} battery capacity, using fast-charging stations, with a 12-year lifespan, and priced at 440,000 USD. 

    \item Investment costs of BEBs include charger  cost upon purchase and battery replacements cost at their mid-life.  The manufacturer commits to providing battery replacements   to account for battery degradation.
    
    \item  All BEBs, regardless of charging type, start the day with a full battery.
    
    \item A ZEB cannot be salvaged before reaching half of its economic lifetime.

    \item Demand is known in advance and is deterministic. 
    
    \item Assuming the passenger capacity is the same for all buses of the same length, each bus can be replaced only by another bus of the same length.

     \item As the cost and energy density of lithium-ion batteries (Li-ion) improve, larger batteries will be installed on BEBs while maintaining the same overall weight, but with a smaller unit cost. This will allow a BEB to cover longer distances. 
    
    \item As the cost and efficiency  of fuel cell systems improve, purchase cost as well as O\&M cost of HFCBs will decrease due to the fact that improved efficiency of fuel cell systems will reduce the energy consumption per unit distance.

    \item Recharging takes a fixed amount of time, depending on the battery capacity of the BEB versions. We assume charging power also advances alongside BEB technological improvements so  that the recharging time remains consistent even after the improvement.
    
    \item Deadheading (trips with no passengers) between terminals is allowed, and  buses consume
    less energy in such trips.
    However, deadheading trips to and from garages are disregarded.

    \item Only tailpipe emissions are considered and the emissions related to electricity and hydrogen production are not included.
\end{itemize}

\vspace{-5mm}

\subsection{Mathematical Formulation}
\label{sec:method-formulation}

In this section, we present the mathematical formulation of our multi-stage stochastic program for transitioning to a clean bus fleet. The model aims to optimize decisions regarding the purchasing, salvaging, and assigning of buses to routes while considering technological advancements over a planning horizon.
We present the index sets,   decision variables and parameters used in our mathematical formulation   in Tables~\ref{tab:Sets}, ~\ref{tab:Variables} and~\ref{tab:Parameters}, respectively.

\begin{table}[H]
    \centering
    \caption{List of index sets.}
    \label{tab:Sets}\scalebox{.95}{
    \begin{tabular}{|c|l||c|l|}
    \hline
    \textbf{Set} & \textbf{Description} & \textbf{Set} & \textbf{Description} \\
    \hline
    $\mathcal{N}$ & Set of  nodes, $\mathcal{N} = \{0, 1, \ldots, N\}$ &
     $\mathcal{T}$ & Set of time periods, $\mathcal{T} = \{0, 1, \ldots, T\}$  \\
    $\mathcal{N}_f$ & Set of  nodes in the final stage &
    $\mathcal{T}_n$ & Set of time periods containing node $n$ \\
    $\mathcal{J}$ & Set of bus types  &
    $\mathcal{T}_{(t']}$ & Set of  time periods \( t \) s.t. $t \le t'$ and $t + \omega_{(j,t)} > t'$ \\
    $\mathcal{K}_j$ & Set of versions for bus type $j$ &
    $\mathcal{T}_{[t')}$ & Set of  time periods \( t \) s.t. $t < t'$ and $t + \omega_{(j,t)} \ge t'$ \\
    $\mathcal{R}$ & Set of routes &
    $\mathcal{T}_{(t')}$ & Set of  time periods \( t \) s.t. $t < t'$ and $t + \omega_{(j,t)} > t'$ \\
    $\mathcal{Q}$ & Set of subperiods of each time period &
    $\mathcal{T}^-_{[t')}$ & Set of  time periods \( t \) s.t. $t+\lceil\omega_{(j,k,t)}/2\rceil \leq t'$ and $t + \omega_{(j,k,t)} \ge t'$ \\
    \hline
    \end{tabular}}
\end{table}

\begin{table}[h]
    \centering
    \caption{List of decision variables.}
    \label{tab:Variables}\scalebox{.95}{
    \begin{tabular}{|c|p{15cm}|}
    \hline
    \textbf{Variable} & \textbf{Description} \\
    \hline
    $v^+_{(j,k,t),n}$ & Number of version $k \in \mathcal{K}_j$ of bus type $j \in \mathcal{J}$ purchased in time period $t \in \mathcal{T}_n$ and node $n \in \mathcal{N}$. \\
    $v_{(j,k,t),t',n}$ & Number of version $k \in \mathcal{K}_j$ of bus type $j \in \mathcal{J}$ purchased in time period $t \in \mathcal{T}$ and available in time period $t' \in \mathcal{T}_n$ and node $n \in \mathcal{N}$, subject to $t \leq t' < t + \omega_{(j,k,t)}$. \\
    $v'_{(j,k,t),t',r,q,n}$ & Number of version $k \in \mathcal{K}_j$ of bus type $j \in \mathcal{J}$ purchased in time period $t \in \mathcal{T}$ and assigned to route $r \in \mathcal{R}$ in subperiod $q \in \mathcal{Q}$ of time period $t' \in \mathcal{T}_n$ and node $n \in \mathcal{N}$, subject to $t \leq t' < t + \omega_{(j,k,t)}$. \\
    $v^-_{(j,k,t),t',n}$ & Number of version $k \in \mathcal{K}_j$ of bus type $j \in \mathcal{J}$ purchased in time period $t \in \mathcal{T}$ and salvaged in time period $t' \in \mathcal{T}_n$ and node $n \in \mathcal{N}$, subject to $t + \lceil\omega_{(j,k,t)}/2\rceil \leq t' \leq t + \omega_{(j,k,t)}$. \\
    \hline
    \end{tabular}}
\end{table}

\begin{table}[H]
    \centering
    \caption{List of parameters.}
    \label{tab:Parameters}\scalebox{.95}{
    \begin{tabular}{|c|p{15cm}|}
    \hline
    \textbf{Parameter} & \textbf{Description} \\
    \hline
    $\pi_n$ & Probability of node $n \in \mathcal{N}$ \\
    $\mu_{(n,t)}$ & The node corresponding to the ancestor of node $n \in \mathcal{N}$ at time $t \in \mathcal{T}$ \\
    $\beta_{nominal}$ & Nominal discount rate applied to time periods \\
    $\zeta$ & Inflation rate applied to time periods\\
    $\beta_{real}$ & Real discount rate, calculated as $\beta_{real} = \frac{1 + \beta_{nominal}}{1 + \zeta} - 1$\\
    $\beta$ & Discount factor, calculated as $\beta = \frac{1}{1 + \beta_{real}}$\\
    $\delta^{+}_{(j,k,t),n}$ & Investment cost of version $k \in \mathcal{K}_j$ of bus type $j \in \mathcal{J}$ purchased in time period $t \in \mathcal{T}_n$ at node $n \in \mathcal{N}$ \\
    $\delta^-_{(j,k,t),t',n}$ & Salvage value of version $k \in \mathcal{K}_j$ of bus type $j \in \mathcal{J}$ purchased in time period $t \in \mathcal{T}$ and salvaged in time period $t' \in \mathcal{T}_n$ at node $n \in \mathcal{N}$ \\
    $\delta_{(j,k,t),t',r,q,n}$ & Operational and maintenance costs of version $k \in \mathcal{K}_j$ of bus type $j \in \mathcal{J}$ purchased in time period $t \in \mathcal{T}$, operated in time period $t' \in \mathcal{T}_n$, assigned to route $r \in \mathcal{R}$ in subperiod $q \in \mathcal{Q}$, and at node $n \in \mathcal{N}$ \\
    $\omega_{(j,k,t)}$ & Economic life of version $k \in \mathcal{K}_j$ of bus type $j \in \mathcal{J}$ purchased in time period $t \in \mathcal{T}_n$ \\
    $\lambda_{(j,k,t), t', r, q, n}$ & Demand satisfaction ratio of version $k \in \mathcal{K}_j$ of bus type $j \in \mathcal{J}$ purchased in time period $t \in \mathcal{T}$, operated in time period $t' \in \mathcal{T}_n$, assigned to route $r \in \mathcal{R}$ in subperiod $q \in \mathcal{Q}$, and at node $n \in \mathcal{N}$ \\
    $\Delta_{t, r, q}$ & Demand in subperiod $q \in \mathcal{Q}$ of route $r \in \mathcal{R}$ at time period $t \in \mathcal{T}$ \\
    $\gamma_t$ & Available budget for time period $t \in \mathcal{T}$ \\
    $\epsilon_{(j,k,t),t'}$ & Emissions of version $k \in \mathcal{K}_j$ of bus type $j \in \mathcal{J}$ purchased in time period $t \in \mathcal{T}$ and operated in time period $t' \in \mathcal{T}$ \\
    $ \eta_{t} $ & Maximum allowable emissions for time period $t \in \mathcal{T}$ \\
    $\phi_{j,k}$ & Initial fleet size of version $k \in \mathcal{K}_j$ of bus type $j \in \mathcal{J}$ \\
    $ \psi $ & Maximum desired average age of the fleet at the end of the planning horizon \\
    \hline
    \end{tabular}}
\end{table}

We now present our multi-stage stochastic program as below:
\begin{subequations}\label{eq:stochasticModel}
\begin{align}
    \min \
    & \ \sum_{n=1}^{N} \pi_n \sum_{j \in \mathcal{J}} \sum_{k \in \mathcal{K}_j} \left( \sum_{t \in \mathcal{T}_n} \beta^{(t-1)} \delta^{+}_{(j,k,t),n} v^{+}_{(j,k,t),n} + \sum_{t' \in \mathcal{T}_n} \sum_{t \in \mathcal{T}_{(t']}}\sum_{r \in \mathcal{R}} \sum_{q \in \mathcal{Q}} \beta^{(t'-1)} \delta_{(j,k,t),t',r,q,n} v'_{(j,k,t),t',r,q,n} \right. \nonumber \\
    & \ \left. - \sum_{t' \in \mathcal{T}_n} \sum_{t \in \mathcal{T}_{[t')}} \beta^{(t'-1)} \delta^-_{(j,k,t),t',n} v^-_{(j,k,t),t',n} \right)  
    \label{eq:objFunc}
\end{align}
\vspace{-8mm}
\begin{align}
  \text{s.t.}  & \ \sum_{j \in \mathcal{J}} \sum_{k \in \mathcal{K}_j}\sum_{t \in \mathcal{T}_{(t']}} \lambda_{(j,k,t), t', r, q, n} \ v'_{(j,k,t), t', r, q, n} \geq \Delta_{t', r, q}, \quad t' \in \mathcal{T}_n, r \in \mathcal{R}, q \in \mathcal{Q}, n \in \mathcal{N} \label{eq:demandConstr} \\
    & \ \sum_{r \in \mathcal{R}} v'_{(j,k,t),t',r,q, n} \leq v_{(j,k,t),t', n}, \quad j \in \mathcal{J}, k \in \mathcal{K}_j, q \in \mathcal{Q}, t' \in \mathcal{T}_n, t \in \mathcal{T}_{(t']}, n \in \mathcal{N} \label{eq:assignConstr} \\
    & \ v_{(j,k,t), t, n} = v^+_{(j,k,t), n}, \quad j \in \mathcal{J}, k \in \mathcal{K}_j, t \in \mathcal{T}_n, n \in \mathcal{N} \label{eq:balance1Constr} \\
    & \ v_{(j,k,t),t',n} = v_{(j,k,t),t'-1,\mu_{(n,t'-1)}} -  v_{(j,k,t),t',n}^-, \quad n\neq0, j \in \mathcal{J}, k \in \mathcal{K}_j, t' \in \mathcal{T}_n, t \in \mathcal{T}_{(t']}, n \in \mathcal{N}
    \label{eq:balance2Constr} \\
    & \ v_{(j,k,t),t',n}^- = v_{(j,k,t),t'-1,\mu_{(n,t'-1)}}, \quad n\neq0, j \in \mathcal{J}, k \in \mathcal{K}_j, t' \in \mathcal{T}_n, t \in \mathcal{T}_{(t']}, t'=t+\omega_{(j,k,t)}, n \in \mathcal{N} \label{eq:lifetimeConstr} \\
    & \ \sum_{j \in \mathcal{J}} \sum_{k \in \mathcal{K}_j} \delta^+_{(j,k,t), n} \ v^+_{(j,k,t), n} \leq \gamma_{t}, \quad t \in \mathcal{T}_n, n \in \mathcal{N} \label{eq:budgetConstr}\\
    & \ \sum_{j \in \mathcal{J}} \sum_{k \in \mathcal{K}_j} \sum_{t \in \mathcal{T}_{(t']}} \sum_{r \in \mathcal{R}} \sum_{q \in \mathcal{Q}} \epsilon_{(j,k,t),t'} \ v'_{(j,k,t),t',r,q,n} \leq \eta_{t'}, \quad t' \in \mathcal{T}_n, n \in \mathcal{N} \label{eq:emissionConstr} \\
    & \ v^+_{(j,k,0),0} = \phi_{j,k}, 
     \ v^-_{(j,k,0),0,0} = 0, \quad j \in \mathcal{J}, k \in \mathcal{K}_j \label{eq:initalFleetConstr} \\
    & \ \sum_{j \in \mathcal{J}} \sum_{k \in \mathcal{K}_j} \sum_{t \in \mathcal{T}_{(T]}} (T-t+1) \ v_{(j,k,t),T,n} \leq \psi \ \sum_{j \in \mathcal{J}} \sum_{k \in \mathcal{K}_j} \sum_{t \in \mathcal{T}_{(T]}} v_{(j,k,t), T, n}, \quad n \in \mathcal{N}_{f} \label{eq:maxAverageAgeConstr}\\
    & \  v_{(j,k,t),n}^+ \in \mathbb{Z}_+, \quad j \in \mathcal{J}, k \in \mathcal{K}_j, t \in \mathcal{T}_n, n \in \mathcal{N} \label{eq:VarDomain1}
    \\
    & \ v_{(j,k,t),t',n}^- \in \mathbb{Z}_+, \quad j \in \mathcal{J}, k \in \mathcal{K}_j, t' \in \mathcal{T}_n, t \in \mathcal{T}_{[T)}, n \in \mathcal{N} \label{eq:VarDomain2}
    \\
    & \ v_{(j,k,t),t',n} \in \mathbb{Z}_+, \ v'_{(j,k,t),t',n} \in \mathbb{Z}_+, \quad j \in \mathcal{J}, k \in \mathcal{K}_j, t' \in \mathcal{T}_n, t \in \mathcal{T}_{(T]}, n \in \mathcal{N}. \label{eq:VarDomain3}
\end{align}
\end{subequations}
The terms in the objective function~\eqref{eq:objFunc} represent the investment costs of the buses, the operational and maintenance costs, and the salvage revenue, respectively.
Constraints~\eqref{eq:demandConstr} require that the number of buses assigned to a route in each time period and node meets the required demand. Constraints~\eqref{eq:assignConstr} specify that the number of assigned buses should not exceed the number of buses available. Constraints~\eqref{eq:balance1Constr} and~\eqref{eq:balance2Constr} maintain the balance of buses for each type and version in every time period. Constraints~\eqref{eq:lifetimeConstr} ensure that the buses will be salvaged once they reach their economic lifetime. Constraints~\eqref{eq:budgetConstr} restrict the total purchasing costs to stay within the allocated budget for each time period, and constraints~\eqref{eq:emissionConstr} ensure that emissions remain within the limits for each time period. Constraints~\eqref{eq:initalFleetConstr} define the initial fleet size for each type and version of the bus, and make sure that no buses are salvaged at stage 0. To mitigate end-of-horizon effects, constraints~\eqref{eq:maxAverageAgeConstr} limit the average age of the buses at the end of the planning horizon. Finally, constraints~\eqref{eq:VarDomain1}, \eqref{eq:VarDomain2}, and \eqref{eq:VarDomain3} ensure that decision variables take non-negative integer values.

\subsection{Energy Requirement Calculation}
\label{sec:method-energy}

In order to ensure operational feasibility of assigning buses to routes, and optimizing the assignment decisions, it is crucial to account for the amount of energy consumed by each version of BEBs on scheduled trips of different routes. 
We now present our physics-based approach  to calculate the energy requirement of a specific BEB assigned to a given service trip, which impacts the O\&M costs and the demand satisfaction ratio parameters.
Table~\ref{tab:energyParameters} provides the notations used in this section. 
\begin{table}[h]
    \centering
    \caption{List of parameters for traction power and battery power calculations.}
    \label{tab:energyParameters}
    \begin{tabular}{|c|p{6cm}||c|p{6.5cm}|}
    \hline
    \textbf{Parameter} & \textbf{Description} & \textbf{Parameter} & \textbf{Description} \\
    \hline
    $f_r$ & Rolling resistance coefficient &
    $\rho_{\text{air}}$ & Air density (\si{\kilo\gram/\cubic\metre}) \\
    $C_D$ & Drag coefficient &
    $g$ & Gravitational acceleration (\si{\metre/\second\squared}) \\
    $A_f$ & Frontal area of the bus (\si{\metre\squared}) &
    $\eta_t$ & Transmission efficiency \\
    $\eta_m$ & Motor and inverter efficiency &
    $\eta_{\text{rb}}$ & Regenerative braking efficiency \\
    $m$ & Mass of the bus (\si{\kilo\gram}) &
    $m_{\text{eq}}$ & Equivalent mass (\si{\kilo\gram}) \\
    $\alpha$ & Road grade (\si{\radian}) &
    $a(\tau)$ & Bus acceleration rate  at time $\tau$ (\si{\metre/\second\squared})      \\
    $v(\tau)$ & Speed of the bus at time $\tau$ (\si{\metre/\second}) &  
    $P_w(\tau)$ & Traction power at time $\tau$ (\si{\watt}) \\
    \hline
    \end{tabular}
\end{table}

To calculate the energy requirement for a bus on a specific service trip, we first divide the trip into a set  \(\mathcal{S}\) of segments, where each segment represents the path between two consecutive bus stops.
We assume that the bus will accelerate between the time interval $[0,\tau_1]$ with a constant rate of $a > 0$, maintain a constant speed between $[\tau_1,\tau_2]$, and decelerate with rate $-a$ between $[\tau_2,\tau_1+\tau_2]$. For each segment, we minimize the duration $\tau_1+\tau_2$ to cover its distance   subject to the constraints on maximum speed and maximum power, where the  \textit{traction power} \( P_w(\tau) \) at time $\tau$ is calculated as
\begin{equation*}
    P_w(\tau) = \left(m g \sin(\alpha) + f_r m g \cos(\alpha) + 0.5 \rho_{\text{air}} C_D A_f v(\tau)^2 + m_{\text{eq}} a(\tau)\right) v(\tau) . 
\end{equation*}
Here, traction power \( P_w(\tau) \) refers to the power required at the wheels  to overcome resistance forces acting on the bus when traveling at a speed of \( v(\tau) \), and  accelerate an equivalent mass of \( m_{\text{eq}} \), which accounts for the inertial resistance of the rotating masses in the vehicle, with an acceleration of \( a(\tau) \).
Required traction power is provided by the   battery at a higher rate due to losses at the electric motor and transmission elements before reaching the wheels. In addition, during braking or travelling downhill, BEBs can recuperate some of the power that otherwise would be lost as heat via their regenerative braking system. The power required from the battery during traction (or recuperated by the battery during regenerative braking) \( P_{\text{bat}}(\tau) \) is given~by 
\begin{equation*}
    P_{\text{bat}}(\tau) = 
    \begin{cases} 
    \frac{P_w(\tau)}{\eta_t \eta_m}, & \text{if } P_w(\tau) \ge 0 \\
    P_w(\tau) \cdot \eta_t \cdot \eta_m \cdot \eta_{\text{rb}}, & \text{if } P_w(\tau) < 0
    \end{cases}.
\end{equation*}
To find the energy requirement \(E_{\text{s}}\) of segment \(s\), we integrate \( P_{\text{bat}}(\tau) \) over the segment duration as
\(
E_{\text{s}} = \int_{0}^{\tau_1+\tau_2} P_{\text{bat}}(\tau) \, d\tau  
\). 
Finally, the total energy required during the trip is the sum of the energy requirements for all segments, computed as
\(
E_{\text{trip}} = \sum_{s\in\mathcal{S}} E_{\text{s}}
\)
. We also account for the variation in energy requirements due to seasonal changes by   multiplying the nominal consumption value by a constant that depends on the ambient temperature. 
{{We provide the details of power and energy calculations in the Supplementary Material.}}

\subsection{Bus-to-Route Simulator}
\label{sec:method-algorithm}

We develop a simple simulator, Algorithm~\ref{AssignBuses_Algorithm}, to estimate key parameters  if a version $k$ of bus type $j$ is assigned to a certain route $r$ in a subperiod $q$. Since the outcome of such an assignment changes over the planning horizon due to technological advances, we run this algorithm for every node $n$ in the scenario tree in an offline manner before solving  our optimization model.

Inputs of the Bus-to-Route Simulator are the the specifications of a bus (such as the type, version, charging scheme in case of BEBs, energy consumption, etc.) and the characteristics of the route (such as the trip schedule, road profile, etc.). The simulator aims to minimize the necessary number of buses in a heuristic manner for each assignment while considering operational feasibility of each bus type.

\begin{algorithm}
\caption{Bus-to-Route Simulator.}
\label{AssignBuses_Algorithm}
\KwIn{Bus parameters, trip schedule, trip information, terminal-terminal trip information}
\KwOut{Bus assignments, recharge details for fast-charging (FC) BEBs, deadheading details}

\For{each trip in the scheduled trips}{
    \eIf{no buses are yet assigned}{
        Assign the first bus to the first trip.
        Move on to the next trip.
    }{
        \For{each assigned bus}{
            \eIf{the bus's garage location and length match with that of the scheduled trip and the bus can reach the trip's start point on time after completing its previous trip, considering deadheading}{
                \If{bus is electric}{
                    \If{battery is insufficient for the trip}{
                        \If{bus type is ONC}{
                            Skip this bus due to insufficient battery.
                        }
                        \ElseIf{bus type is FC}{
                            \If{the bus can start the trip with the recharge time added}{
                                \If{bus has enough energy to get to the starting point of the trip}{
                                    Plan a recharge at the starting point of the trip.
                                }
                                \Else{
                                    Plan a recharge where the bus finishes its previous assigned trip.
                                }
                                Record recharge details (location, start and end times).
                                Assign the trip to the bus.
                                Update cumulative energy consumed.
                                Reset battery capacity.
                                Move on to the next trip.
                            }
                            \Else{
                                Skip this bus due to insufficient time for recharging.
                            }
                        }
                    }
                }
                Assign the trip to this bus.
                \If{bus is electric}{
                    Update cumulative energy consumed.
                    Update remaining battery capacity.
                }
                Move on to the next trip.
            }{Skip this bus.}
        }
    }
    \If{no current buses can be assigned}{
        Assign a new bus to this trip.
    }
}

\For{each bus}{
    Calculate the assigned distances and energy consumed (if bus is electric), including deadheading.
    \If{the total assigned distance is less than 10 km}{
        Disregard that bus.
    }
}
Record recharge summary for all fast-charging electric buses and all locations.
\end{algorithm}

The algorithm outputs the number of necessary buses and their daily task assignments, using which we obtain several parameters used in our optimization model: 
i) Demand satisfaction ratio $\lambda_{(j,k,t), t', r, q, n}$, defined as the number of buses needed of a specific type-version pair divided by the number of DBs needed.  
ii) Charger-to-bus ratio for fast-charging BEBs, used in the   calculation of the investment cost $\delta^+_{(j,k,t), n}$.
iii) Average daily distance covered, used in the calculation of the O\&M cost  $\delta_{(j,k,t), t', r, q, n}$.

\section{Data Collection and Processing}
\label{sec:data}

In this section, we present the data used in our case study in Section~\ref{sec:case}. We start by discussing the bus network data provided by the public bus fleet operator of Istanbul in Section~\ref{sec:iett}, including a description of the preprocessing steps. Section~\ref{sec:buses} details the initial bus-specific cost information included in our case study.
{In Section~\ref{sec:tech_imp}, we provide the data required for forecasting future technological advancements in BEBs and HFCBs, along with an explanation of our forecasting approach.} Finally, in Section~\ref{sec:scenario-gen}, we explain how the scenarios of our stochastic program are obtained using the resulting projections. {We note that the Supplementary Material contains more detailed datasets, analyses and results.}

\subsection{Public Bus Transit in Istanbul}
\label{sec:iett}

The buses operated by the Istanbul Electricity, Tram and Tunnel Establishments (IETT), the authority responsible for public bus transportation in Istanbul, carry nearly 5 million people daily and cover approximately 1.2 million kilometers \parencite{iett-faaliyet2023}. IETT provided us the following datasets:

\begin{enumerate}
    \item Trip Schedule: This dataset includes the details of service trips scheduled for both March 15, 2023 (Winter Schedule) and August 3, 2023 (Summer Schedule). The summer schedule is assumed to be used only for 92 days during the summer, while the winter schedule is utilized for the remainder of the year. The specific information provided in this dataset is the route code, trip ID, scheduled route start times, scheduled distance, vehicle depot (when available) and vehicle length group.

    \item Stop Sequence with Coordinates: This dataset provides the sequence of stops for each route and route type, along with their latitudes and longitudes.
    Furthermore, we use the Open-Elevation API \parencite{open-elevation} to determine  the elevation of each stop. 

    \item Vehicle Information:
     This dataset includes information on 3,351 buses operated by IETT, detailing the number of buses at each depot, along with their brand, model, and manufacturing year.
    It does not cover
     3,076 buses under a private   bus brand of IETT,
     for which the detailed data is not available.
\end{enumerate}

We process the above datasets to obtain the following pieces of information.
\begin{enumerate}
    \item {Estimation of Missing Distances:} As we do not have the full information about the specific routing of each vehicle, we estimate some distances  related to the deadheading services using the Haversine approximation when needed in Algorithm~\ref{AssignBuses_Algorithm} in addition to the distance between consecutive stops.

    \item {Estimation of Travel Time:}
    Through our preliminary study and  dataset received from IETT, we determine that the average speed of buses is approximately 25 \unit{\kilo\metre/\hour}. Consequently, we assume a maximum speed of 30 \unit{\kilo\metre/\hour}.

    \item {Calculating Bus Demand:}
    We run our simulator, Algorithm~\ref{AssignBuses_Algorithm}, to estimate the number of buses needed of each version   to meet the scheduled trips in each route, for both summer and winter schedules. 
    We present the summary of results in Table~\ref{tab:demandSummary} for DBs. As the total number of buses owned by IETT in 2023 is reported to be 6,652 \parencite{iett-faaliyet2023}
    the algorithm provides a reliable estimation of the actual demand.
    \begin{table}[h]
    \centering
    \caption{Total DB demand by 
    length for winter and summer schedules.}
    \label{tab:demandSummary}
    \begin{tabular}{|c|c|c|c|c|c|c|}
        \hline
        \textbf{Bus Length (\unit{\meter})} & \textbf{ 6.5-8 } & \textbf{ 8-9 } & \textbf{ 10-11 } & \textbf{ 11-14 } & \textbf{ 14-19 } & \textbf{ Total } \\ \hline
        \textbf{ Winter/Summer Schedule } & 280/269 & 26/23 & 12/15 & 4540/4349 & 1695/1418 & \textbf{ 6553/6071 } \\ 
        \hline
    \end{tabular}
\end{table}

    \item {Route Aggregation:}  As a preprocessing step, we aggregate routes based on a metric related to the demand satisfaction ratio (DSR). For each route, we calculate the minimum DSR for each BEB version across all bus lengths and seasons using today’s technology. These minimum DSRs are then averaged across all BEB versions to create a single metric for each route.
    Routes are then grouped into 12 clusters based on this metric as follows: $\{ \{1.00\}, [0.95, 1.00), [0.90, 0.95), \dots, [0.50, 0.55), [0.38, 0.50) \}$.
    As an example, cluster with the metric of 1.00 
    represents the routes where even the BEB with the smallest battery capacity can cover the scheduled trips in all seasons and all bus length groups used. {Since only 13 routes have a metric smaller than 0.50, we group them together into one cluster.} 

\end{enumerate}

\vspace{-5mm}

\subsection{Initial Bus-Specific Costs}
\label{sec:buses}

We consider four bus models across different length groups: an 8\unit{\meter} model for the 6.5-8\unit{\meter} group, a 10\unit{\meter} model for both the 8-9 and 10-11\unit{\meter} groups, a 12\unit{\meter} model for the 11-14\unit{\meter} group, and an 18\unit{\meter} model for the 14-19\unit{\meter} group.
Bus purchase costs are estimated from recent tenders in Turkey and Europe,  information gathered from local bus manufacturers and technical reports, and given in Table~\ref{tab:purchase_costs}.
The investment costs of BEBs include battery replacement and charger costs, in addition to the bus purchase costs. We assume that each BEB requires one battery replacement at the end of year six of its operation that is adjusted for the discount rate. The initial battery pack cost, denoted by~$\text{BC}$, is assumed to be \si{500 USD/\kWh}. For charging infrastructure, our study includes regular chargers with \si{50\kW} charging power, priced at 20,000 USD per unit for ONC buses, and fast chargers with \si{350\kW} charging power, priced at 45,000 USD per unit for fast charging buses.
Each ONC bus requires one dedicated charger while charger-to-bus ratio for FC buses with a specific version is estimated using Algorithm~\ref{AssignBuses_Algorithm} and the total charger cost is distributed across all buses. Recharging times for FC buses are based on a full recharge, assuming a 90\% charging efficiency.
Finally, we assume that  salvage values depreciate yearly by 15\% for all buses.
\begin{table}[H]
    \centering
    \caption{Purchase costs (USD) for DBs, BEBs, and HFCBs.}
    \label{tab:purchase_costs}\scalebox{.88}{
    \begin{tabular}{|c|c|c|c|c|c|c|c|c|c|c|c|c|}
    \hline
    \textbf{Model} & \textbf{DB} & \multicolumn{5}{c|}{\textbf{BEB Fast Charging}} & \multicolumn{5}{c|}{\textbf{BEB Overnight Charging}} & \textbf{HFCB} \\
    \cline{3-12}
    \textbf{Length} & & \textbf{140 \unit{\kWh}} & \textbf{210 \unit{\kWh}} & \textbf{280 \unit{\kWh}} & \textbf{350 \unit{\kWh}} & \textbf{420 \unit{\kWh}} & \textbf{280 \unit{\kWh}} & \textbf{350 \unit{\kWh}} & \textbf{420 \unit{\kWh}} & \textbf{490 \unit{\kWh}} & \textbf{560 \unit{\kWh}} & \\
    \hline
    8\unit{\meter} & 135,000 & 305,000 & 340,000 & - & - & - & 375,000 & 410,000 & - & - & - & 500,000 \\
    10\unit{\meter} & 170,000 & 340,000 & 375,000 & 410,000 & - & - & 410,000 & 445,000 & 480,000 & - & - & 600,000 \\
    12\unit{\meter} & 200,000 & 370,000 & 405,000 & 440,000 & 475,000 & - & 440,000 & 475,000 & 510,000 & 545,000 & - & 700,000 \\
    18\unit{\meter} & 300,000 & 470,000 & 505,000 & 540,000 & 575,000 & 610,000 & 540,000 & 575,000 & 610,000 & 645,000 & 680,000 & 1,000,000 \\
    \hline
    \end{tabular}}
\end{table}
\noindent

The O\&M costs include energy, maintenance, and driver costs. We use Algorithm~\ref{AssignBuses_Algorithm} to determine the average daily distance the buses of a specific  version cover when assigned a specific route. This information along with the average energy requirement for BEBs is used to calculate the total energy costs and maintenance costs, where the unit costs are given in Table~\ref{tab:maintenance_energy_summary}. 
We assume that the energy consumption estimation is done in spring/fall, and increase BEB and HFCB consumption rates by 15\% for winter and 5\% for summer, consistent with Istanbul's climate. 
Regarding the driver cost,  the information gathered from IETT Activity Reports and Financial Statements suggest that
the approximate salary is nearly double the minimum wage. Additionally, the number of drivers is  about 20\% more than the number of  buses. Therefore we assume the driver cost per bus to be 2.4 times the minimum wage, approximately 40 USD/day.

\begin{table}[H]
    \small
    \centering
    \caption{Maintenance costs, energy consumption and energy costs for all bus types.}
    \label{tab:maintenance_energy_summary}
    \begin{tabular}{|c|c|c|c|}
    \hline
    \textbf{Bus Type} & \textbf{Energy Consumption} & \textbf{Maintenance Cost (USD/\unit{\kilo\meter})} & \textbf{Energy Cost (USD per unit)} \\
    \hline
    DB & 0.435 \si{\liter/\km} \parencite{ma2021examining} & 0.58 \parencite{holland2021environmental} & 1.29 USD/\unit{\liter} \parencite{petrolofisi2024}\\
    BEB & Varies with route and version  & 0.34 \parencite{holland2021environmental} & 0.16 USD/\unit{\kWh} \parencite{encazip2024} \\
    HFCB & 0.09 \si{\kg/\km} \parencite{ajanovic2021prospects} & 0.29 \parencite{collins2022fuel} & 10.00 USD/\unit{\kg}  \\
    \hline
    \end{tabular}
\end{table}

\subsection{Technological Change Forecasts}
\label{sec:tech_imp}

In this section, we examine the cost and efficiency trends of BEBs and HFCBs. 
We quantify these trends for each year with respect to five years prior, consistent with our case study (Section~\ref{sec:case}). Specifically, for each year, the efficiency improvement rate is calculated as the natural logarithm of the year’s value divided by the value from five years earlier, while cost improvement is calculated as the negative logarithm of the same ratio.
We then plot the scatter plot of these improvement rates, which highlights distinct patterns and trends across five-year periods. Clustering techniques are applied to group the data points based on their similarities, with each cluster representing a possible scenario of technological advancement. The probabilities associated with each cluster indicate the proportion of data points that belong to that cluster.

\subsubsection{BEB Technological Change}\label{sec:tech_imp-beb}

\begin{figure}[H]
\centering

\begin{subfigure}{0.32\textwidth}
\centering
    \begin{tikzpicture}
    \pgfplotsset{
        width=\textwidth,
        height=6cm,
        xmin=1990, xmax=2023,
        y axis style/.style={
            yticklabel style=#1,
            ylabel style=#1,
            y axis line style=#1,
            ytick style=#1
       }
    }
          \pgfplotsset{every axis y label/.append style={rotate=0,yshift=-.45cm}}
    \begin{axis}[
           xticklabel style={/pgf/number format/1000 sep=},
      axis y line*=left,
      y axis style=black!75!black,
      ymin=100, ymax=10000,
              ymode=log,
      xlabel=Year,
      ylabel=Battery Cost (2023 USD/\unit{\kWh})
    ]
    \addplot[mark=*,black] 
      coordinates{
(1991,8805.2)
(1992,7166.22)
(1993,5756.38)
(1994,6451.6)
(1995,6068.66)
(1996,5151.87)
(1997,4398.27)
(1998,3541.91)
(1999,2741.66)
(2000,2576.53)
(2001,1832.35)
(2002,1327.7)
(2003,1014.27)
(2004,909.97)
(2005,788.02)
(2006,675.05)
(2007,644.22)
(2008,673.49)
(2009,580.73)
(2010,533.49)
(2011,499.68)
(2012,507.12)
(2013,472.73)
(2014,409.52)
(2015,280.86)
(2016,235.02)
(2017,170.31)
(2018,142.6)
(2019,119.58)
(2020,111.57)
(2021,107.44)
(2022,119.23)
(2023,100.63)
    };
    \end{axis}

      \pgfplotsset{every axis y label/.append style={rotate=180,yshift=5.02cm}}
    \begin{axis}[
      axis y line*=right,
      axis x line=none,
      ymin=90, ymax=1000,
              ymode=log,
      ylabel=Energy Density(\si{\watt\hour/\kilogram}),
      y axis style=black!75!black
    ]
    \addplot[mark=o,black] 
      coordinates{
(1991,98.2)
(1992,102.75)
(1993,107.86)
(1994,112.98)
(1995,118.66)
(1996,126.62)
(1997,130.59)
(1998,144.8)
(1999,156.37)
(2000,160.13)
(2001,179.87)
(2002,188.18)
(2003,193.37)
(2004,209.99)
(2005,202.72)
(2006,207.92)
(2007,222.46)
(2008,218.3)
(2009,222.46)
(2010,230.77)
(2011,238.04)
(2012,243.23)
(2013,243.23)
(2014,251.54)
(2015,248.43)
(2016,259.85)
(2017,276.47)
(2018,293.09)
(2019,313.87)
(2020,334.64)
(2021,389.7)
(2022,444.75)
(2023,499.8)
    };
    \end{axis}
    
    \end{tikzpicture}
\caption{Cell cost (left, solid circle) vs. energy density (right, hollow circle)   in log-scale.
}
\label{fig:scatter_plot_cost_vs_density}
\end{subfigure}
\hfill
\begin{subfigure}{0.32\textwidth}
\centering
\begin{tikzpicture}
          \pgfplotsset{every axis y label/.append style={rotate=0,yshift=-.3cm}}
    \begin{axis}[
        width=\textwidth,
        height=6cm,
        xlabel={Energy Density Impr. Rate},
        ylabel={Cost Impr.  Rate},
        grid=major,
        xmin=0, xmax=0.6,
        ymin=0, ymax=1.5,
        xtick={0,0.1,...,0.6},
        ytick={0,0.3,...,1.5},
        minor y tick num=2,
        scatter/classes={
            1={mark=triangle,blue},
            2={mark=diamond,green},
            3={mark=square,orange},
            representative1={mark=triangle*,draw=blue,fill=blue,scale=1.5},
            representative2={mark=diamond*,draw=green,fill=green,scale=1.5},
            representative3={mark=square*,draw=orange,fill=orange,scale=1.5}
        }
    ]
    \addplot[scatter,only marks,scatter src=explicit symbolic] 
    table[meta=label] {
        x y label
        0.37 1.20 1
        0.29 1.25 1
        0.29 1.10 1
        0.24 1.18 1
        0.14 1.00 1
        0.13 1.09 1
        0.19 1.20 1
        0.22 1.23 1
        0.33 0.86 1
        0.30 0.86 1
        0.35 1.03 1
        0.30 0.92 1
        0.41 0.78 1
    };
    \addplot[scatter,only marks,scatter src=explicit symbolic] 
    table[meta=label] {
        x y label
        0.09 0.75 3
        0.17 0.72 3
        0.07 0.64 3
        0.25 0.54 3
        0.24 0.49 3
        0.29 0.49 3
        0.12 0.41 3
        0.06 0.45 3
        0.13 0.39 3
        0.14 0.30 3
        0.09 0.24 3
        0.11 0.35 3
        0.12 0.35 3
        0.48 0.36 3
        0.53 0.35 3
    };

    \addplot[scatter,only marks,scatter src=explicit symbolic] 
    table[meta=label] {
        x y label 
        0.272711 1.054408 representative1
        0.192712 0.455076 representative3
    };
    
    \node[below] at (axis cs:0.2718,1.065) {\scriptsize \textbf{F}};
    \node[below] at (axis cs:0.192,0.439) {\scriptsize \textbf{S}};
    
    \end{axis}
\end{tikzpicture}
\caption{Five-year energy density vs. cost improvement rates for BEBs with two clusters.}
\label{fig:scatter_plot_2_clusters}
\end{subfigure}
\hfill
\begin{subfigure}{0.32\textwidth}
\centering
\begin{tikzpicture}
    \label{fig:3by2BEB}
          \pgfplotsset{every axis y label/.append style={rotate=0,yshift=-.3cm}}
    \begin{axis}[
        width=\textwidth,
        height=6cm,
        xlabel={Energy Density Impr.  Rate},
        ylabel={Cost Impr.  Rate},
        grid=major,
        xmin=0, xmax=0.6,
        ymin=0, ymax=1.5,
        xtick={0,0.1,...,0.6},
        ytick={0,0.3,...,1.5},
        minor y tick num=2,
        scatter/classes={
            1={mark=triangle,blue},
            2={mark=diamond,green}, 
            3={mark=square,orange},
            representative1={mark=triangle*,draw=blue,fill=blue,scale=1.5},
            representative2={mark=diamond*,draw=green,fill=green,scale=1.5},
            representative3={mark=square*,draw=orange,fill=orange,scale=1.5}
        }
    ]
    \addplot[scatter,only marks,scatter src=explicit symbolic] 
    table[meta=label] {
        x y label
        0.37 1.20 1
        0.29 1.25 1
        0.29 1.10 1
        0.24 1.18 1
        0.14 1.00 1
        0.13 1.09 1
        0.19 1.20 1
        0.22 1.23 1
        0.35 1.03 1
    };
    \addplot[scatter,only marks,scatter src=explicit symbolic] 
    table[meta=label] {
        x y label
        0.33 0.86 2
        0.30 0.86 2
        0.30 0.92 2
        0.41 0.78 2
        0.09 0.75 2
        0.17 0.72 2
        0.07 0.64 2
    };
    \addplot[scatter,only marks,scatter src=explicit symbolic] 
    table[meta=label] {
        x y label
        
        0.25 0.54 3
        0.24 0.49 3
        0.29 0.49 3
        0.12 0.41 3
        0.06 0.45 3
        0.13 0.39 3
        0.14 0.30 3
        0.09 0.24 3
        0.11 0.35 3
        0.12 0.35 3
        0.48 0.36 3
        0.53 0.35 3
    };
    \addplot[scatter,only marks,scatter src=explicit symbolic] 
    table[meta=label] {
        x y label
        0.246349 1.143194 representative1
        0.23666 0.791069 representative2
        0.213477 0.392256 representative3
    };

    \node[below] at (axis cs:0.2463,1.16) {\scriptsize \textbf{F}};
    \node[below] at (axis cs:0.2385,0.79) {\scriptsize \textbf{M}};
    \node[below] at (axis cs:0.213,0.392) {\scriptsize \textbf{S}};
    
    \end{axis}
\end{tikzpicture}
\caption{Five-year energy density vs. cost improvement rates for BEBs with three clusters.}
\label{fig:scatter_plot_3_clusters}
\end{subfigure}

\caption{BEB technology improvement charts. Cluster centers are marked with solid marks in (b) and (c), and can be distinguished by the abbreviations of \textbf{S}low, \textbf{F}ast and \textbf{M}edium.
}
\label{fig:scatter_plot_comparison}
\end{figure}
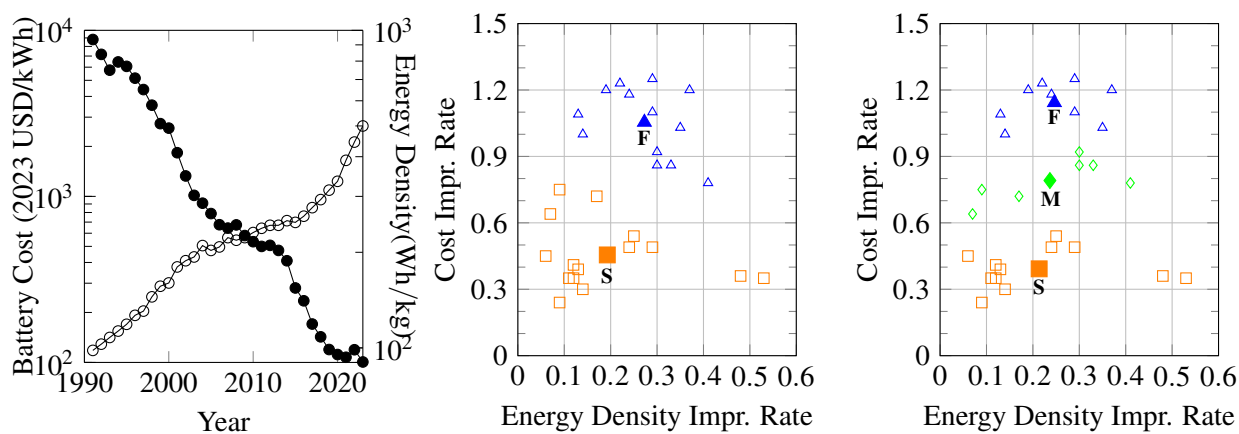

The cost of  Li-ion  batteries per kilowatt-hour has significantly decreased over the years  while their energy density has steadily increased. Figure~\ref{fig:scatter_plot_cost_vs_density} presents the data on energy density and battery cost from 1991 onwards, extracted from the charts in reference~\textcite{walter2023} using the WebPlotDigitizer tool~\parencite{autometris}.
Given the cost and energy density values, we compute the five-year improvement ratios as described above. Then, these ratios are clustered into two and three different groups by minimizing the sum of squared Euclidean distances between each data point and its assigned cluster centers can be seen in Figures~\ref{fig:scatter_plot_2_clusters} and \ref{fig:scatter_plot_3_clusters}. 
As an example, let us consider the two-cluster setting, in which case the probabilities for the fast improvement cluster (\textbf{F}) and the slow improvement cluster (\textbf{S}) are 0.46 and 0.54, respectively.
In fast improvement cluster,  the energy density and cost improvement rates are 0.27 and 1.05, respectively. Taking the exponential of these rates, we get an energy density multiplier of 1.31 
and a cost multiplier of~0.35.

\vspace{-8mm}

\subsubsection{HFCB Technological Change}\label{sec:tech_imp-hfcb}

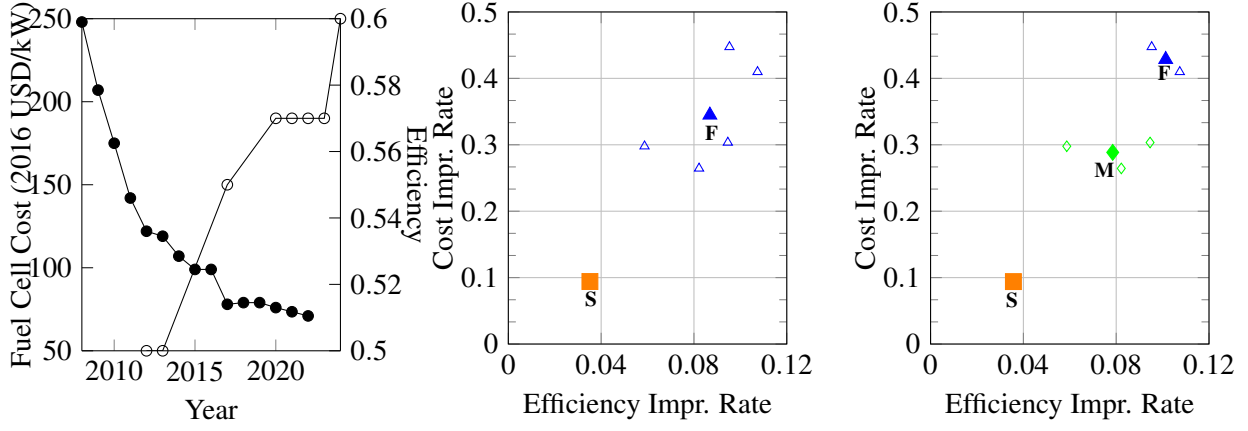
\begin{figure}[H]
\centering

\begin{subfigure}{0.32\textwidth}
\centering
    \begin{tikzpicture}
    \pgfplotsset{
        width=0.95\textwidth,
        height=6cm,
        xmin=2008, xmax=2024,
        y axis style/.style={
            yticklabel style=#1,
            ylabel style=#1,
            y axis line style=#1,
            ytick style=#1
       }
    }
          \pgfplotsset{every axis y label/.append style={rotate=0,yshift=-.45cm}}
    \begin{axis}[
           xticklabel style={/pgf/number format/1000 sep=},
      axis y line*=left,
      y axis style=black!75!black,
      ymin=50, ymax=250,
      xlabel=Year,
      ylabel=Fuel Cell Cost (2016 USD/\unit{\kW})
    ]
    \addplot[mark=*,black] 
      coordinates{
(2008,248)
(2009,207)
(2010,175)
(2011,142)
(2012,122)
(2013,119)
(2014,107)
(2015,99)
(2016,99)
(2017,78)
(2018,79)
(2019,79)
(2020,76)
(2021,73.5)
(2022,71)
    };
    \end{axis}

      \pgfplotsset{every axis y label/.append style={rotate=180,yshift=5.2cm}}
    \begin{axis}[
      axis y line*=right,
      axis x line=none,
      ymin=0.50, ymax=0.60,
      ylabel=Efficiency,
      y axis style=black!75!black,
        yticklabel style={
        /pgf/number format/fixed,
        /pgf/number format/precision=2
        },
        scaled y ticks=false,
    ]
    \addplot[mark=o,black] 
      coordinates{
(2012,0.50)
(2013,0.50)
(2017,0.55)
(2020,0.57)
(2021,0.57)
(2022,0.57)
(2023,0.57)
(2024,0.60)
    };
    \end{axis}
    
    \end{tikzpicture}
\caption{Fuel cell system cost (left, solid circle) vs. efficiency (right, hollow circle).
}
\label{fig:scatter_plot_cost_vs_density-fc}
\end{subfigure}
\hfill
\begin{subfigure}{0.32\textwidth}
\centering
\begin{tikzpicture}
          \pgfplotsset{every axis y label/.append style={rotate=0,yshift=-.4cm}}
    \begin{axis}[
        width=\textwidth,
        height=6cm,
        xlabel={Efficiency Impr. Rate},
        ylabel={Cost Impr.  Rate},
        grid=major,
        xmin=0, xmax=0.12,
        ymin=0, ymax=0.5,
        xtick={0,0.04,...,0.12},
        ytick={0,0.1,...,0.5},
        xticklabel style={
        /pgf/number format/fixed,
        /pgf/number format/precision=2
        },
        scaled x ticks=false,
        minor y tick num=2,
        scatter/classes={
            1={mark=triangle,blue},
            2={mark=diamond,green},
            3={mark=square,orange},
            representative1={mark=triangle*,draw=blue,fill=blue,scale=1.5},
            representative2={mark=diamond*,draw=green,fill=green,scale=1.5},
            representative3={mark=square*,draw=orange,fill=orange,scale=1.5}
        }
    ]
    \addplot[scatter,only marks,scatter src=explicit symbolic] 
    table[meta=label] {
        x y label
    0.0953	0.4473 1
    0.1074	0.4097 1
    0.0946	0.3034 1
    0.0587	0.2978 1
    0.0822	0.2644 1 
    };
    \addplot[scatter,only marks,scatter src=explicit symbolic] 
    table[meta=label] {
        x y label
    0.0357	0.0940 3
    };
    \addplot[scatter,only marks,scatter src=explicit symbolic] 
    table[meta=label] {
        x y label 
        0.086857 0.344751 representative1
        0.035160 0.094123 representative3
    };
    
    \node[below] at (axis cs:0.0876,0.345) {\scriptsize \textbf{F}};
    \node[below] at (axis cs:0.0357,0.094) {\scriptsize \textbf{S}};
    
    \end{axis}
\end{tikzpicture}
\caption{Five-year efficiency vs. cost improvement rates for HFCBs with two clusters.}
\label{fig:scatter_plot_2_clusters-fc}
\end{subfigure}
\hfill
\begin{subfigure}{0.32\textwidth}
\centering
\begin{tikzpicture}
          \pgfplotsset{every axis y label/.append style={rotate=0,yshift=-.4cm}}
    \begin{axis}[
        width=\textwidth,
        height=6cm,
        xlabel={Efficiency Impr.  Rate},
        ylabel={Cost Impr.  Rate},
        grid=major,
        xmin=0, xmax=0.12,
        ymin=0, ymax=0.5,
        xtick={0,0.04,...,0.12},
        ytick={0,0.1,...,0.5},
        xticklabel style={
        /pgf/number format/fixed,
        /pgf/number format/precision=2
        },
        scaled x ticks=false,
        minor y tick num=2,
        scatter/classes={
            1={mark=triangle,blue},
            2={mark=diamond,green}, 
            3={mark=square,orange},
            representative1={mark=triangle*,draw=blue,fill=blue,scale=1.5},
            representative2={mark=diamond*,draw=green,fill=green,scale=1.5},
            representative3={mark=square*,draw=orange,fill=orange,scale=1.5}
        }
    ]
    \addplot[scatter,only marks,scatter src=explicit symbolic] 
    table[meta=label] {
        x y label
    0.0953	0.4473 1
    0.1074	0.4097 1
    };
    \addplot[scatter,only marks,scatter src=explicit symbolic] 
    table[meta=label] {
        x y label
    0.0946	0.3034 2
    0.0587	0.2978 2
    0.0822	0.2644 2 
    };
    \addplot[scatter,only marks,scatter src=explicit symbolic] 
    table[meta=label] {
        x y label
    0.0357	0.0940 3
    };
    \addplot[scatter,only marks,scatter src=explicit symbolic] 
    table[meta=label] {
        x y label
        0.101335 0.4285 representative1
        0.078503 0.288547 representative2
        0.0357 0.094 representative3
    };

    \node[below] at (axis cs:0.1008,0.435) {\scriptsize \textbf{F}};
    \node[below] at (axis cs:0.075,0.288) {\scriptsize \textbf{M}};
    \node[below] at (axis cs:0.0351,0.093) {\scriptsize \textbf{S}};
    
    \end{axis}
\end{tikzpicture}
\caption{Five-year efficiency vs. cost improvement rates for HFCBs with three clusters.}
\label{fig:scatter_plot_3_clusters-fc}
\end{subfigure}

\caption{HFCB technology improvement charts.
}
\label{fig:scatter_plot_comparison-fc}
\end{figure}

Fuel cell systems used in HFCBs have also improved in terms of cost and efficiency. However, since these systems are only being used for heavy duty purposes recently, the data is scarce. Instead, we  use the cost figures for light duty fuel cell systems~\parencite{FCcost} and the efficiency of a specific producer (Ballard) as a proxy~as shown in Figure~\ref{fig:scatter_plot_cost_vs_density-fc}. We then follow a similar approach to Li-ion batteries: We  compute five-year cost~and efficiency improvements (for the common years in two datasets), and obtain the clusters  in Figures~\ref{fig:scatter_plot_2_clusters-fc} and~\ref{fig:scatter_plot_3_clusters-fc}. 

\subsection{Scenario Tree Generation}\label{sec:scenario-gen}

We now formalize how we obtain the scenarios of our multi-stage stochastic program  using the technological change forecasts. Our approach is built on a two-step procedure: In the first step, we  obtain \textit{technology trees} for each bus type illustrating  future projections while in the second step, we combine these individual projections to build the scenario tree.

\begin{algorithm}[H]
\caption{Technology Tree Construction Algorithm}
\label{alg:techTree}
\KwIn{Improvement distribution of technology $j$, the number of stages $S$.}
\KwOut{Technology tree.}

Set $\theta_{j1}^c=\theta_{j1}^e=\theta_{j1}^p=1$, $\text{ID}=1$ and $\mathcal{L}=\{1\}$.

\For{$s=1,\dots,S-1$}{

    Set $\mathcal{L}'=\emptyset$.

    \For{$\ell \in \mathcal{L}$}{

        \For{$b=1,\dots,B_j$}{
        
           Set $\text{ID}=\text{ID}+1$, $h=\text{ID}$.

           Set $\theta_{jh}^c = \theta_{j\ell}^c \Theta_{jb}^c$, 
            $\theta_{jh}^e = \theta_{j\ell}^e \Theta_{jb}^e$,       $\theta_{jh}^p = \theta_{j\ell}^p \Theta_{jb}^p$.

           Set $\mathcal{L}'=\mathcal{L}' \cup \{ \text{ID} \}$.
        
        }
        
    }
    $\mathcal{L}=\mathcal{L}'$.

}

\end{algorithm}

Let us start with the first step. For a bus technology  $j \in \mathcal{J}$, let us denote its cost  and efficiency improvement multipliers for each stage by random variables $\Theta_j^c$ and $\Theta_j^e$, respectively, which are assumed to be independent from each other. We will assume that each pair of random variables take values from a joint discrete probability distribution with a sample space consisting of elements $(\Theta_{jb}^c, \Theta_{jb}^e)$ w.p. $\Theta_{jb}^p$, $b=1,\dots,B_j$. 
In particular, for a technology $j$ with a given support size $B_j$, we use the clustering approach from Section~\ref{sec:tech_imp} with $B_j$ clusters, and use efficiency and cost multipliers of each cluster $b$ as $(\Theta_{jb}^c, \Theta_{jb}^e)$ values and the fraction of points in each cluster $b$ as $\Theta_{jb}^p$. 
We run Algorithm~\ref{alg:techTree} to construct a perfect $B_j$-ary tree denoted by $\mathcal{N}(j)$, which we will call as the technology tree. 
Figure~\ref{fig:tech_trees} illustrates exemplary technology trees for DBs, BEBs and HFCBs for $S=2$ stages with respect to the data reported in Section~\ref{sec:tech_imp}.

{
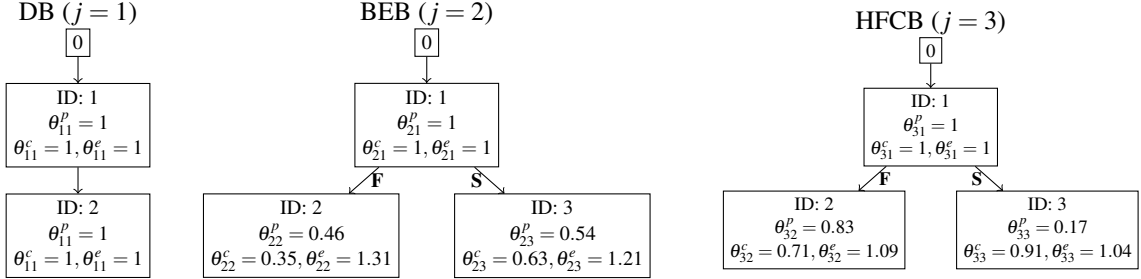
\begin{figure}[H]
\centering
\begin{minipage}{0.15\textwidth}
\centering
\makebox[0pt]{\text{\small DB ($j=1$)}}\par
\begin{tikzpicture}[scale=0.7, transform shape, ->]
\node (0) [rectangle, draw, align=center] {0};
\node (1) [rectangle, draw, align=center, minimum size=1cm][below=0.5cm and 0.0cm of 0] {ID: 1 \\  $\theta_{11}^p=1$ \\ 
$\theta_{11}^c=1, \theta_{11}^e=1$ };
\node (2) [rectangle, draw, align=center, minimum size=1cm][below=0.5cm and 0.0cm of 1] {ID: 2 \\  $\theta_{11}^p=1$ \\ 
$\theta_{11}^c=1, \theta_{11}^e=1$ };
\draw (0) -- (1);
\draw (1) -- (2);
\end{tikzpicture}
\end{minipage}
\begin{minipage}{0.4\textwidth}
\centering
\makebox[0pt]{\text{\small BEB ($j=2$)}}\par
\begin{tikzpicture}[scale=0.7, transform shape, node distance=3.5cm, ->]
\node (0) [rectangle, draw, align=center] {0};
\node (1) [rectangle, draw, align=center, minimum size=1cm][below=0.5cm and 0.0cm of 0] {ID: 1 \\  $\theta_{21}^p=1$ \\ 
$\theta_{21}^c=1, \theta_{21}^e=1$ };
\node (2) [rectangle, draw, align=center, minimum size=1cm][below left=0.5cm and -0.85cm of 1] {ID: 2 \\  $\theta_{22}^p=0.46$ \\ $\theta_{22}^c=0.35, \theta_{22}^e=1.31$ };
\node (3) [rectangle, draw, align=center, minimum size=1cm][below right=0.5cm and -0.85cm of 1] {ID: 3 \\  $\theta_{23}^p=0.54$ \\ $\theta_{23}^c=0.63, \theta_{23}^e=1.21$ };  
\draw (0) -- (1);    
\draw (1) -- (2) node [midway, right] {\textbf{F}} ;
\draw (1) -- (3) node [midway, left] {\textbf{S}} ;
\end{tikzpicture}
\end{minipage}
\begin{minipage}{0.4\textwidth}
\centering
\makebox[0pt]{\text{\small HFCB  ($j=3$)}}\par
\begin{tikzpicture}[scale=0.65, transform shape, node distance=3.5cm, ->]
\node (0) [rectangle, draw, align=center] {0};
\node (1) [rectangle, draw, align=center, minimum size=1cm][below=0.5cm and 0.0cm of 0] {ID: 1 \\  $\theta_{31}^p=1$ \\ 
$\theta_{31}^c=1, \theta_{31}^e=1$ };
\node (2) [rectangle, draw, align=center,  minimum size=1cm][below left=0.5cm and -0.85cm of 1] {ID: 2 \\  $\theta_{32}^p=0.83$ \\ $\theta_{32}^c=0.71, \theta_{32}^e=1.09$ };
\node (3) [rectangle, draw, align=center, minimum size=1cm][below right=0.5cm and -0.85cm of 1] {ID: 3 \\  $\theta_{33}^p=0.17$ \\ $\theta_{33}^c=0.91, \theta_{33}^e=1.04$ };
\draw (0) -- (1);   
\draw (1) -- (2) node [midway, right] {\textbf{F}} ;
\draw (1) -- (3) node [midway, left] {\textbf{S}} ;
\end{tikzpicture}
\end{minipage}
\caption{Technology trees of DB ($B_1=1$), BEB ($B_2=2$) and HFCB ($B_3=2$) with $S=2$ stages.}
\label{fig:tech_trees}
\end{figure}
}

\begin{figure}[H]
\centering
\begin{tikzpicture}[scale=0.65, transform shape, node distance=2cm and 2cm, ->]
\node (0) [rectangle, draw, align=center] {0};
\node (1) [rectangle, draw, align=center, below of=0] {ID: 1 \\ $\pi_1=1$ \\ 
$\vartheta_{11}^c=1, \vartheta_{11}^e=1$ \\
$\vartheta_{21}^c=1, \vartheta_{21}^e=1$ \\
$\vartheta_{31}^c=1, \vartheta_{31}^e=1$};
\node (2) [rectangle, draw, align=center, below left of=1, xshift=-8cm, yshift=-2cm]  {ID: 2 \\ $\pi_1=0.3818$ \\ 
$\vartheta_{12}^c=1, \vartheta_{12}^e=1$ \\
$\vartheta_{22}^c=0.35, \vartheta_{22}^e=1.31$ \\
$\vartheta_{32}^c=0.71, \vartheta_{32}^e=1.09$};
\node (3) [rectangle, draw, align=center, below left of=1, xshift=-1.25cm, yshift=-2cm]{ID: 3 \\ $\pi_1=0.0782$ \\ 
$\vartheta_{13}^c=1, \vartheta_{13}^e=1$ \\
$\vartheta_{22}^c=0.35, \vartheta_{23}^e=1.31$ \\
$\vartheta_{33}^c=0.91, \vartheta_{33}^e=1.04$};
\node (4) [rectangle, draw, align=center, below right of=1, xshift=1.25cm, yshift=-2cm] {ID: 4 \\ $\pi_1=0.4482$ \\ 
$\vartheta_{14}^c=1, \vartheta_{14}^e=1$ \\
$\vartheta_{24}^c=0.63, \vartheta_{24}^e=1.21$ \\
$\vartheta_{34}^c=0.71, \vartheta_{34}^e=1.09$};
\node (5) [rectangle, draw, align=center, below right of=1, xshift=8cm, yshift=-2cm]{ID: 5 \\ $\pi_1=0.0918$ \\ 
$\vartheta_{15}^c=1, \vartheta_{15}^e=1$ \\
$\vartheta_{25}^c=0.63, \vartheta_{25}^e=1.21$ \\
$\vartheta_{35}^c=0.91, \vartheta_{35}^e=1.04$};
\draw (0) -- (1) ;
\draw (1) -- (2) node [midway, above] {\textbf{F}\textbf{F}} ;
\draw (1) -- (3) node [near start, left] {\textbf{F}\textbf{S}} ;
\draw (1) -- (4) node [near start, right] {\textbf{S}\textbf{F}} ;
\draw (1) -- (5) node [midway, above] {\textbf{S}\textbf{S}} ;
\end{tikzpicture}
\caption{Scenario tree with two stages.}
\label{fig:scenario_tree}
\end{figure}
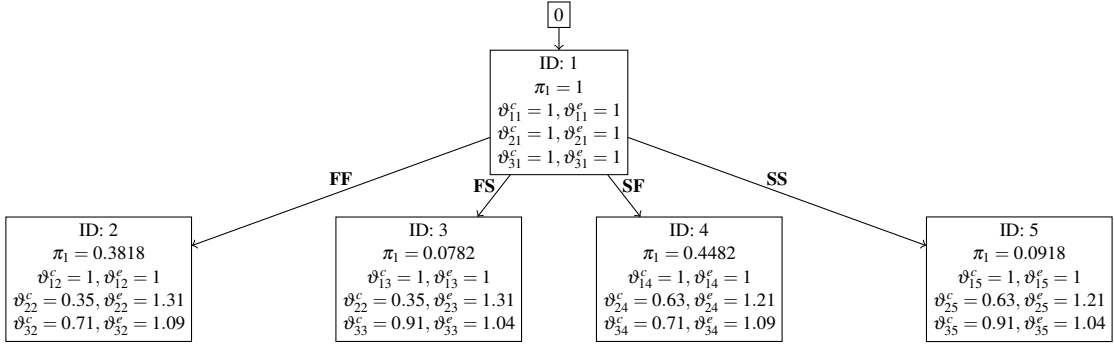

In the second step, once the technology trees are at hand for each $j\in\mathcal{J}$, we obtain the scenario tree~$\mathcal{N}$ by taking the Cartesian product of each node of the technology tree $\mathcal{N}(j)$ at the same level. We label the nodes in the scenario tree in  lexicographic order with respect to the node IDs in the technology tree to establish the bijection $n \in \mathcal{N} \leftrightarrow (n_1, n_2, \dots, n_{|\mathcal{J}|}) \in \mathcal{N}(1) \times \mathcal{N}(2) \times \cdots \mathcal{N}(|\mathcal{J}|)   $. 
This construction can be best explained with an example: Let us consider the scenario tree in Figure~\ref{fig:scenario_tree}, which is obtained from the technology trees in Figure~\ref{fig:tech_trees}. 
The labels of nodes 2, 3, 4, 5 in the scenario tree corresponds to the node ID triplets (2,2,2), (2,2,3), (2,3,2), (2,3,3) in the technology trees.

\noindent
Finally, for each $n\in\mathcal{N}$ and $j\in\mathcal{J}$, we set 
\(
    \pi_n = \prod_{j\in\mathcal{J}} \theta_{j,n_j}^p \text{ and }
    ( \vartheta_{jn}^c, \vartheta_{jn}^e ) =  ( \theta_{j,n_j}^c, \theta_{j,n_j}^e ) .
\)
We note that in a general situation with  $|\mathcal{J}|$ many technologies with respective branches of $B_j$ over $S$ stages, the scenario tree will contain $(\prod_{j=1}^J B_j)^{S-1}$ leaf nodes (or scenarios), which grows exponentially.

We are now ready to explain how the uncertain parameters 
are affected in each scenario. We use the abbreviation $\text{IC}$ for  investment cost, which does not change for  DBs over the planning horizon. 

For BEBs, technological advances affect the battery capacity (as we assume that the total bus weight remains unchanged)  and investment cost.  For  a specific version with initial battery capacity of $C_1$, the battery capacity in node $n$  becomes 
$C_n = C_1 \times \vartheta^e_{BEB, n}$. 
Then, we run Algorithm~\ref{AssignBuses_Algorithm} to determine the new DSR with battery capacity $ C_n$. 
Regarding the   IC of a BEB of length $L$ in node $n$ with battery capacity $C_n$, we use 
\begin{equation}\label{eq:IC-BEB-costEvolution}
    \text{IC}_{\text{BEB},n} (L, C_n) =  \text{IC}_{\text{DB},1} (L) + \text{CC} +
    (1+\beta^{\omega/2}) [ \text{PC} +  \text{BC} \times  C_n ] \vartheta^c_{\text{BEB}, n}.    
\end{equation}
Here, 
 $\text{PC}$ refers to the cost of other components of the electrified powertrain except the battery
  and $\text{CC}$ stands for the charger cost.

For HFCBs, technological advances affect the investment    and O\&M costs. In particular, the IC of a HFCB of length $L$ in node $n$ is computed as
\(
 \text{IC}_{\text{HFCB},n} (L) = \text{IC}_{\text{DB},1} (L) + 
  [ \text{IC}_{\text{HFCB},1} (L) - \text{IC}_{\text{DB},1} (L) ] \vartheta^c_{\text{HFCB}, n}
\) 
and its consumption rate per distance (CR), which affects the unit energy cost, becomes
\(
    \text{CR}_{\text{HFCB},n}  = \text{CR}_{\text{HFCB},1}  /   \vartheta^e_{\text{HFCB}, n}.
\)

\section{Case Study}
\label{sec:case}

This section presents our case study that focuses on the clean  fleet transition in Istanbul public bus network. After we provide our computational setup in Section~\ref{sec:case-case}, we present the detailed results about our base case in Section~\ref{sec:case-base}. Then, we provide a thorough sensitivity analyses regarding some key parameters in our case study in Section~\ref{sec:case-sensitivity}. {Finally, we compare our stochastic programming model with other approaches in Section~\ref{sec:case-comparison}.}

\subsection{Computational Setup}
\label{sec:case-case}

We conduct our computational experiments in the Python programming language using
  a 64-bit workstation with two Intel(R)
Xeon(R) Gold 6248R CPU (3.00GHz) processors (256 GB  RAM). 
Due to the large-scale nature of the multi-stage stochastic program~\eqref{eq:stochasticModel}, we relax the integrality restrictions and solve it as a linear program (LP)  utilizing  {Gurobi 11}.
{\color{black}
Since we cluster the routes into categories with a reasonably large demand, continuous variables provide an accurate approximation that can be rounded to obtain practical solutions.} In fact, in all cases considered,  the difference between the optimal value of the LP relaxation, denoted by $z_{\text{LP}}$ and the objective function value of the rounded solution, denoted by $z_{\text{round}}$, is at most 0.23\%. 

\subsection{{{Base Case}}}
\label{sec:case-base}

In the {\tt{Base Case}}, we obtain a dynamic strategic plan of Istanbul's transition to a clean fleet, spanning a planning horizon of 25 years (2025-2049) divided into five stages. We utilize the input data described in Section~\ref{sec:data}, and set the investment budget in million USD as $\gamma_t=\min\{50t, 250\}$ for $t=1,\dots,25$. 
We also include yearly emission targets, reducing the yearly emission gradually to zero up to year 2049 at a linear rate starting from year 2035. {We set the inflation rate $\zeta=0.04$ and the nominal discount rate as $\beta_{nominal}=0.05$.}
{The initial fleet consists of 2,072 7-year-old, 2,817 11-year-old, 318 14-year-old, 635 16-year-old, and 711 18-year-old DBs. The maximum lifetime for buses with an initial age of 11 years or less is assumed to be 16 years, while buses that are 14 years or older are assumed to have a maximum lifetime of 20 years. The remaining lifetime of the initial buses is adjusted accordingly based on these assumptions. Since the  fleet is quite old, we assume zero salvage value the existing buses.} For newly purchased buses, the economic lifetime is set at 12 years for BEBs and HFCBs, and 15 years for diesel buses. The maximum average fleet age at the end of the planning horizon is capped at 9 years.

\begin{figure}[H]
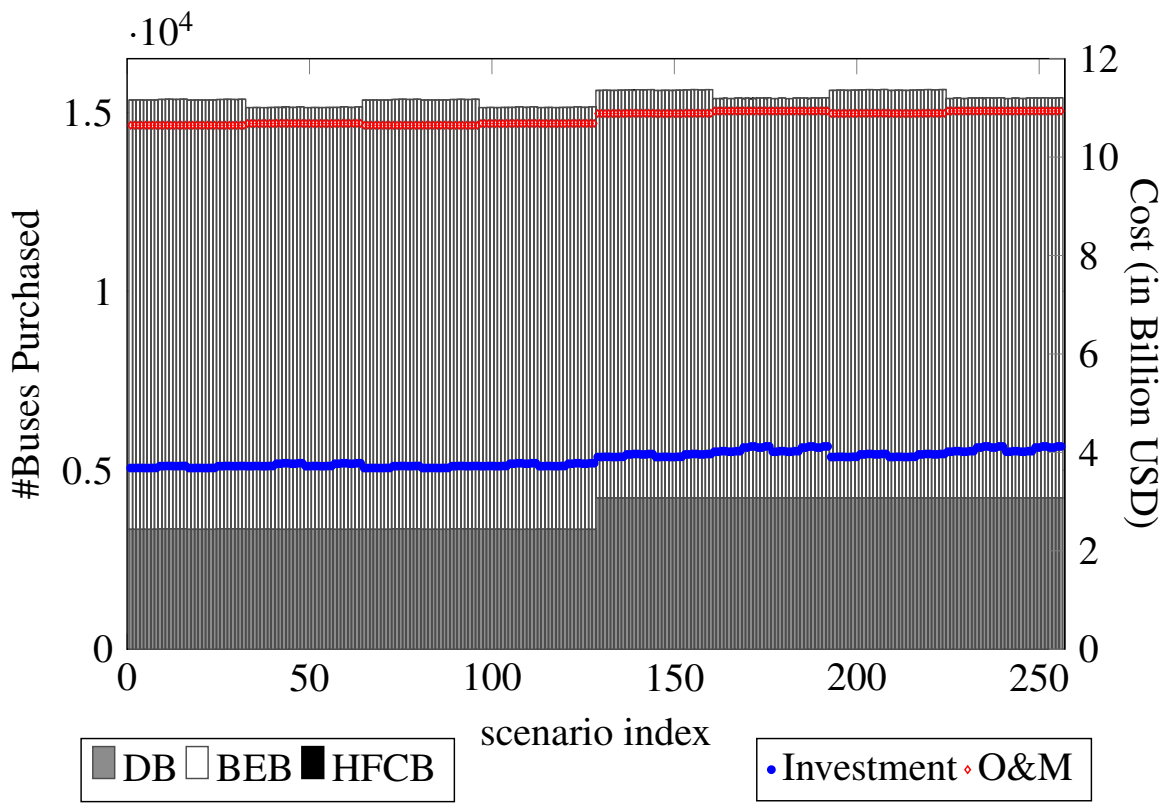

\centering

\caption{{\tt{Base Case}} results. Buses purchased are reported with bars  and the cost figures are reported with markers.}
\label{fig:base-case}
\end{figure}

We consider three technology alternatives: DBs, BEBs and HFCBs. We assume that the cost of DBs do not change, but the cost and efficiency of BEBs and HFCBs improve as explained in Section~\ref{sec:tech_imp} with respect to two clusters each. Therefore, we consider $(1\times2\times2)^{5-1}=256$ scenarios in our multi-stage stochastic program.
The resulting model contains more than 23 million variables and 3 million constraints.
The CPU time is 5,435 seconds, and the expected objective values are  $z_{\text{LP}}=14.126$B USD and $z_{\text{round}}=14.154$B USD, 
indicating only a $0.2\%$ difference.
For each scenario, we report the number of buses purchased of each technology type along with the total purchase and O\&M costs over the planning horizon in Figure~\ref{fig:base-case}. %
Results  indicate that the O\&M costs are more than twice of the investments costs, which is expected for public buses heavily utilized in Istanbul. In addition, BEBs dominate across all scenarios, while HFCBs are purchased in small numbers and only in a few scenarios.

\begin{figure}[H]
\centering

\begin{subfigure}{0.475\textwidth}
\centering
\begin{tikzpicture}[scale=1.1]
      \pgfplotsset{every axis y label/.append style={rotate=0,yshift=-0.4cm}}
\begin{axis}[
    title={Fast-Fast},
    x=6pt,
    xmin=0,
    xmax=26,
    ymin = 0,
    ymax = 8,     ytick pos = left,
    bar width=3pt, 
    ybar stacked,
    ylabel={\#Active Buses (in Thousand)},    
    xlabel={years},
]

 \addplot+ [ybar, fill=gray!90,draw=black!70] coordinates {
(1,6.52)
(2,6.522)
(3,6.523)
(4,6.523)
(5,6.523)
(6,5.734)
(7,4.959)
(8,4.242)
(9,3.622)
(10,2.94)
(11,2.005)
(12,1.77)
(13,1.77)
(14,0.868)
(15,0.348)
(16,0.325)
(17,0.227)
(18,0.208)
(19,0)
(20,0)
(21,0)
(22,0)
(23,0.003)
(24,0.003)
(25,0)
};

 \addplot+ [ybar, fill=gray!0,draw=black!70] coordinates {
(1,0.041)
(2,0.041)
(3,0.041)
(4,0.041)
(5,0.041)
(6,0.855)
(7,1.723)
(8,2.502)
(9,3.162)
(10,3.808)
(11,4.829)
(12,5.126)
(13,5.086)
(14,6.067)
(15,6.596)
(16,6.612)
(17,6.7)
(18,6.705)
(19,6.914)
(20,6.907)
(21,6.872)
(22,6.87)
(23,6.865)
(24,6.865)
(25,6.864)
};

 \addplot+ [ybar, black] coordinates {
(1,0)
(2,0)
(3,0)
(4,0)
(5,0)
(6,0)
(7,0)
(8,0)
(9,0)
(10,0)
(11,0)
(12,0)
(13,0)
(14,0)
(15,0)
(16,0)
(17,0)
(18,0)
(19,0)
(20,0)
(21,0)
(22,0)
(23,0)
(24,0)
(25,0)
};

\end{axis}
      \pgfplotsset{every axis y label/.append style={rotate=180,yshift=7.2cm}}
    \begin{axis}[
    x=6pt,
    xmin=0,
    xmax=26,
    ymin=0,
    ymax=0.7,     ytick={0, 0.1, 0.2, 0.3, 0.4, 0.5, 0.6,0.7},
    axis y line*=right,
    axis x line=none,
    ylabel=Cost (in Billion USD),
    ]

 \addplot[mark=*,blue, mark size=1pt, only marks] coordinates {
(1,0.05013221)
(2,0.09934476)
(3,0.14735267)
(4,0.19463147)
(5,0.24061138)
(6,0.23843625)
(7,0.23685161)
(8,0.23425979)
(9,0.23196815)
(10,0.22989374)
(11,0.22751475)
(12,0.06540969)
(13,0.00033442)
(14,0.22107594)
(15,0.21930168)
(16,0.21703647)
(17,0.21525760)
(18,0.21286601)
(19,0.21111669)
(20,0.06718602)
(21,0.12057391)
(23,0.00072914)
};

 \addplot[mark=diamond,red, mark size=1pt, only marks] coordinates {
(1,0.63623224)
(2,0.63053081)
(3,0.62452740)
(4,0.61871347)
(5,0.61284226)
(6,0.57508596)
(7,0.53115561)
(8,0.49591285)
(9,0.46447753)
(10,0.44173407)
(11,0.39775474)
(12,0.38346127)
(13,0.38132764)
(14,0.34840207)
(15,0.33608189)
(16,0.33188010)
(17,0.32767692)
(18,0.32420980)
(19,0.31972652)
(20,0.31713769)
(21,0.31556754)
(22,0.31251885)
(23,0.30956529)
(24,0.30663147)
(25,0.30387637)
};

    \end{axis}
\end{tikzpicture}
\label{fig:BaseFixed-Fast-Fast-scenario}
\end{subfigure}
\hfill
\begin{subfigure}{0.475\textwidth}
\centering
\begin{tikzpicture}[scale=1.1]
      \pgfplotsset{every axis y label/.append style={rotate=0,yshift=-.4cm}}
\begin{axis}[
    title={Fast-Slow},
    x=6pt,
    xmin=0,
    xmax=26,
    ymin=0,
    ymax = 8,     ytick pos = left,
    bar width=3pt, 
    ybar stacked,
    ylabel={\#active buses (in Thousand)},    
    xlabel={years},
]
 \addplot+ [ybar, fill=gray!90,draw=black!70] coordinates {
(1,6.52)
(2,6.522)
(3,6.523)
(4,6.523)
(5,6.523)
(6,5.734)
(7,4.959)
(8,4.242)
(9,3.622)
(10,2.94)
(11,2.005)
(12,1.77)
(13,1.77)
(14,0.868)
(15,0.348)
(16,0.325)
(17,0.227)
(18,0.208)
(19,0)
(20,0)
(21,0)
(22,0)
(23,0.003)
(24,0.003)
(25,0)
};

 \addplot+ [ybar, fill=gray!0,draw=black!70] coordinates {
(1,0.041)
(2,0.041)
(3,0.041)
(4,0.041)
(5,0.041)
(6,0.855)
(7,1.723)
(8,2.502)
(9,3.162)
(10,3.808)
(11,4.829)
(12,5.126)
(13,5.086)
(14,6.067)
(15,6.596)
(16,6.613)
(17,6.7)
(18,6.705)
(19,6.914)
(20,6.908)
(21,6.873)
(22,6.871)
(23,6.866)
(24,6.866)
(25,6.865)
};

 \addplot+ [ybar, black] coordinates {
(1,0)
(2,0)
(3,0)
(4,0)
(5,0)
(6,0)
(7,0)
(8,0)
(9,0)
(10,0)
(11,0)
(12,0)
(13,0)
(14,0)
(15,0)
(16,0)
(17,0)
(18,0)
(19,0)
(20,0)
(21,0)
(22,0)
(23,0)
(24,0)
(25,0)
};

\end{axis}
      \pgfplotsset{every axis y label/.append style={rotate=180,yshift=7.2cm}}
    \begin{axis}[
    x=6pt,
    xmin=0,
    xmax=26,
    ymin=0,
    ymax=0.7,     ytick={0, 0.1, 0.2, 0.3, 0.4, 0.5, 0.6,0.7},
    axis y line*=right,
    axis x line=none,
    ylabel=Cost (in Billion USD),
    ]

 \addplot[mark=*,blue, mark size=1pt, only marks] coordinates {
(1,0.05013221)
(2,0.09934476)
(3,0.14735267)
(4,0.19463147)
(5,0.24061138)
(6,0.23843625)
(7,0.23685161)
(8,0.23425979)
(9,0.23196815)
(10,0.22989374)
(11,0.22751475)
(12,0.06540969)
(13,0.00033442)
(14,0.22107594)
(15,0.21930168)
(16,0.21701023)
(17,0.21525760)
(18,0.21286601)
(19,0.21108485)
(20,0.06718602)
(21,0.12057391)
(23,0.00072914)
};

 \addplot[mark=diamond,red, mark size=1pt, only marks] coordinates {
(1,0.63623224)
(2,0.63053081)
(3,0.62452740)
(4,0.61871347)
(5,0.61284226)
(6,0.57512217)
(7,0.53105144)
(8,0.49585063)
(9,0.46440219)
(10,0.44171206)
(11,0.39776281)
(12,0.38347293)
(13,0.38136423)
(14,0.34837526)
(15,0.33608628)
(16,0.33193382)
(17,0.32774414)
(18,0.32428538)
(19,0.31984756)
(20,0.31727165)
(21,0.31557909)
(22,0.31255269)
(23,0.30959138)
(24,0.30667158)
(25,0.30395893)
};

    \end{axis}
\end{tikzpicture}
\label{fig:BaseFixed-Fast-Slow-scenario}
\end{subfigure}

\begin{subfigure}{0.475\textwidth}
\centering
\begin{tikzpicture}[scale=1.1]
      \pgfplotsset{every axis y label/.append style={rotate=0,yshift=-.4cm}}
\begin{axis}[
    title={Slow-Fast},
    x=6pt,
    xmin=0,
    xmax=26,
    ymin = 0,
    ymax= 8,
    ytick pos = left,
    bar width=3pt,
    ybar stacked,
    ylabel={\#active buses (in Thousand)},    
    xlabel={years},
]

 \addplot+ [ybar, fill=gray!90,draw=black!70] coordinates {
(1,6.52)
(2,6.522)
(3,6.523)
(4,6.523)
(5,6.523)
(6,6.023)
(7,5.408)
(8,4.855)
(9,4.331)
(10,3.995)
(11,3.304)
(12,2.58)
(13,2.58)
(14,1.86)
(15,1.201)
(16,0.41)
(17,0.323)
(18,0.322)
(19,0.025)
(20,0.025)
(21,0.008)
(22,0.006)
(23,0.006)
(24,0.007)
(25,0)
};

 \addplot+ [ybar, fill=gray!0,draw=black!70] coordinates {
(1,0.041)
(2,0.041)
(3,0.041)
(4,0.041)
(5,0.041)
(6,0.558)
(7,1.269)
(8,1.854)
(9,2.408)
(10,2.683)
(11,3.483)
(12,4.284)
(13,4.243)
(14,5.027)
(15,5.726)
(16,6.555)
(17,6.648)
(18,6.639)
(19,6.928)
(20,6.907)
(21,6.905)
(22,6.876)
(23,6.857)
(24,6.854)
(25,6.854)

};

 \addplot+ [ybar, black] coordinates {
(1,0)
(2,0)
(3,0)
(4,0)
(5,0)
(6,0)
(7,0)
(8,0)
(9,0)
(10,0)
(11,0)
(12,0)
(13,0)
(14,0)
(15,0)
(16,0)
(17,0)
(18,0)
(19,0)
(20,0)
(21,0)
(22,0.013)
(23,0.025)
(24,0.026)
(25,0.026)
};

\end{axis}
      \pgfplotsset{every axis y label/.append style={rotate=180,yshift=7.2cm}}
    \begin{axis}[
    x=6pt,
    xmin=0,
    xmax=26,
    ymin=0,
    ymax=0.7,     ytick={0, 0.1, 0.2, 0.3, 0.4, 0.5, 0.6,0.7},
    axis y line*=right,
    axis x line=none,
    ylabel=Cost (in Billion USD),
    ]

 \addplot[mark=*,blue, mark size=1pt, only marks] coordinates {
(1,0.05013221)
(2,0.09934476)
(3,0.14735267)
(4,0.19463147)
(5,0.24061138)
(6,0.23844826)
(7,0.23649614)
(8,0.23426266)
(9,0.23193217)
(10,0.23002190)
(11,0.22741165)
(12,0.22527387)
(14,0.22122040)
(15,0.21896552)
(16,0.21672812)
(17,0.21497116)
(18,0.21262321)
(19,0.21109562)
(20,0.09531080)
(21,0.20683116)
(22,0.16204243)
(23,0.00310219)
(24,0.00033243)
};

 \addplot[mark=diamond,red, mark size=1pt, only marks] coordinates {
(1,0.63623224)
(2,0.63053081)
(3,0.62452740)
(4,0.61871347)
(5,0.61284226)
(6,0.58665790)
(7,0.54966920)
(8,0.52042061)
(9,0.49343462)
(10,0.48354660)
(11,0.44554431)
(12,0.41274597)
(13,0.41016942)
(14,0.38003688)
(15,0.35654378)
(16,0.33250699)
(17,0.32794466)
(18,0.32456434)
(19,0.31983590)
(20,0.31777168)
(21,0.31580908)
(22,0.31384821)
(23,0.31096754)
(24,0.30801911)
(25,0.30547031)
};

    \end{axis}
\end{tikzpicture}
\label{fig:BaseFixed-Slow-Fast-scenario}
\end{subfigure}
\hfill
\begin{subfigure}{0.475\textwidth}
\centering
\begin{tikzpicture}[scale=1.1]
      \pgfplotsset{every axis y label/.append style={rotate=0,yshift=-.4cm}}
\begin{axis}[
    title={Slow-Slow},
    x=6pt,
    xmin=0,
    xmax=26,
    ymin = 0,
    ymax=  8,
    ytick pos = left,
    bar width=3pt, 
    ybar stacked,
    ylabel={\#active buses (in Thousand)},    
    xlabel={years},
]

 \addplot+ [ybar, fill=gray!90,draw=black!70] coordinates {
(1,6.52)
(2,6.522)
(3,6.523)
(4,6.523)
(5,6.523)
(6,6.023)
(7,5.408)
(8,4.855)
(9,4.331)
(10,3.995)
(11,3.304)
(12,2.58)
(13,2.58)
(14,1.86)
(15,1.201)
(16,0.41)
(17,0.323)
(18,0.322)
(19,0.025)
(20,0.025)
(21,0.008)
(22,0.006)
(23,0.006)
(24,0.007)
(25,0)
};

 \addplot+ [ybar, fill=gray!0,draw=black!70] coordinates {
(1,0.041)
(2,0.041)
(3,0.041)
(4,0.041)
(5,0.041)
(6,0.558)
(7,1.27)
(8,1.855)
(9,2.409)
(10,2.684)
(11,3.484)
(12,4.285)
(13,4.244)
(14,5.028)
(15,5.727)
(16,6.555)
(17,6.647)
(18,6.637)
(19,6.925)
(20,6.904)
(21,6.903)
(22,6.894)
(23,6.893)
(24,6.891)
(25,6.891)
};

 \addplot+ [ybar, black] coordinates {
(1,0)
(2,0)
(3,0)
(4,0)
(5,0)
(6,0)
(7,0)
(8,0)
(9,0)
(10,0)
(11,0)
(12,0)
(13,0)
(14,0)
(15,0)
(16,0)
(17,0)
(18,0)
(19,0)
(20,0)
(21,0)
(22,0)
(23,0)
(24,0)
(25,0)
};

\end{axis}
      \pgfplotsset{every axis y label/.append style={rotate=180,yshift=7.2cm}}
    \begin{axis}[
    x=6pt,
    xmin=0,
    xmax=26,
    ymin=0,
    ymax=0.7,     ytick={0, 0.1, 0.2, 0.3, 0.4, 0.5, 0.6,0.7},
      axis y line*=right,
      axis x line=none,
      ylabel=Cost (in Billion USD),
    ]

 \addplot[mark=*,blue, mark size=1pt, only marks] coordinates {
(1,0.05013221)
(2,0.09934476)
(3,0.14735267)
(4,0.19463147)
(5,0.24061138)
(6,0.23844826)
(7,0.23682949)
(8,0.23426266)
(9,0.23193217)
(10,0.23002190)
(11,0.22741165)
(12,0.22524382)
(14,0.22122040)
(15,0.21896552)
(16,0.21672812)
(17,0.21497116)
(18,0.21262321)
(19,0.21088017)
(20,0.09531080)
(21,0.20683116)
(22,0.16635223)
(24,0.00010833)
};

 \addplot[mark=diamond,red, mark size=1pt, only marks] coordinates {
(1,0.63623224)
(2,0.63053081)
(3,0.62452740)
(4,0.61871347)
(5,0.61284226)
(6,0.58683661)
(7,0.54982753)
(8,0.52051890)
(9,0.49356285)
(10,0.48356972)
(11,0.44546758)
(12,0.41275966)
(13,0.41024130)
(14,0.38003312)
(15,0.35649918)
(16,0.33242891)
(17,0.32782990)
(18,0.32454291)
(19,0.31977006)
(20,0.31778539)
(21,0.31569561)
(22,0.31373487)
(23,0.31066207)
(24,0.30770102)
(25,0.30515265)
};

    \end{axis}
\end{tikzpicture}
\label{fig:base-case-Slow-Slow-scenario}
\end{subfigure}

\caption{{\tt{Base Case}}-four specific scenarios.}

\label{fig:base-case-4scenarios}

\end{figure}
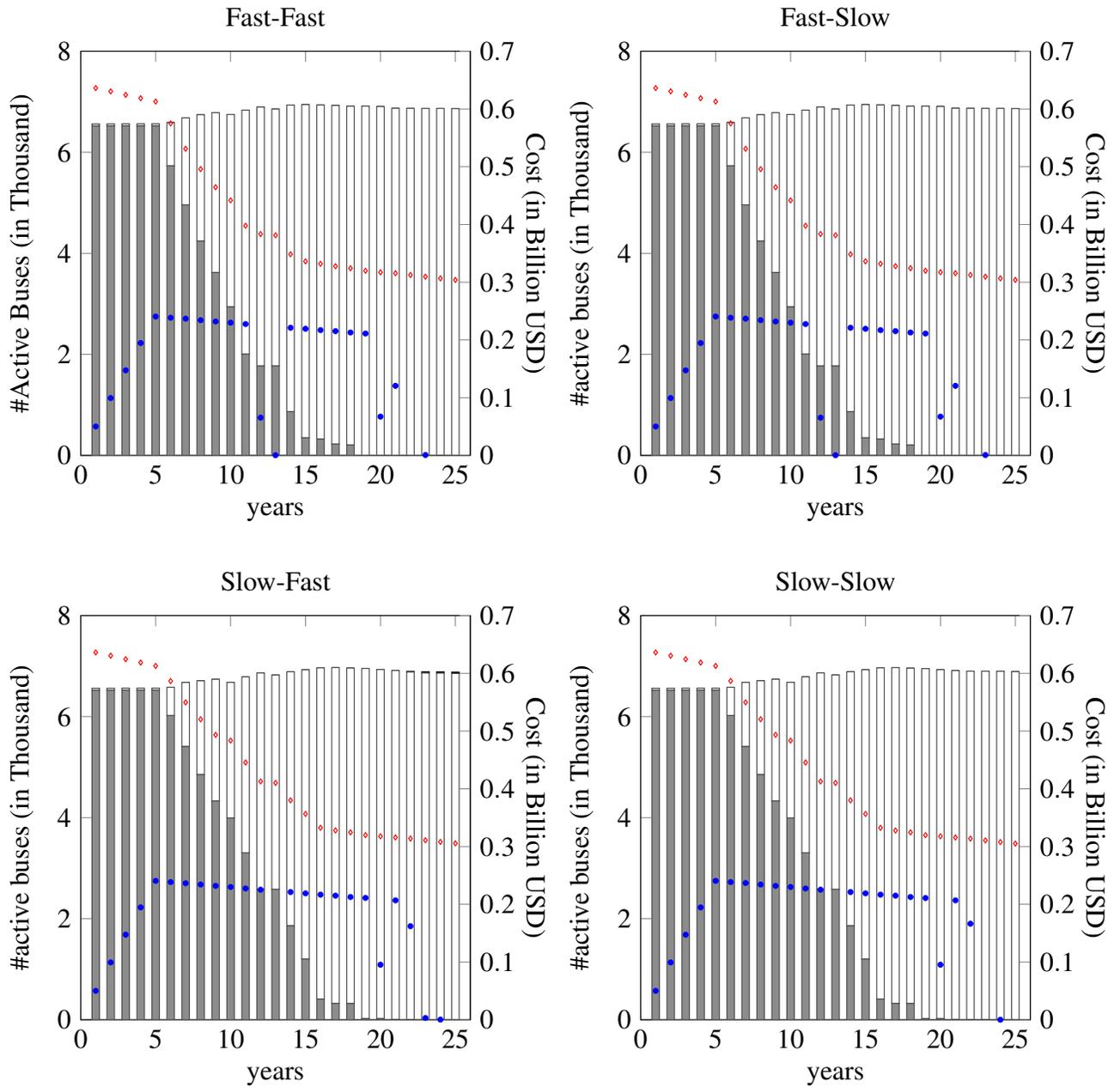

We also provide detailed results for four particular scenarios in Figure~\ref{fig:base-case-4scenarios}: i) Fast-Fast, ii) Fast-Slow, iii) Slow-Fast, iv) Slow-Slow. Here, Fast-Slow refers to the scenario where BEBs improves fast  while HFCBs improve slow in each stage (the other three scenarios are defined accordingly). These four scenarios provide extreme cases and help us illustrate the effect of stochasticity in the dynamic strategic plans.
As the initial fleet is quite old,  most of the  buses  need to be retired within the first six years and the budget limitation is most influential in this time frame for all scenarios. 
We recall that  heavily-utilized city buses typically have higher O\&M cost in their economic lifetime compared to the initial investment cost. The purchase of BEBs, even in the first year (albeit in small quantities due to budget limitations), shows their advantage over DBs in saving O\&M costs. This is especially true for routes where electrification can happen without needing many more buses, that is, routes with  DSR values close to 1.
{We also observe that the fleet size peaks in year 15 for the Fast-Fast and Fast-Slow scenarios, and in year 17 for the Slow-Fast and Slow-Slow scenarios. In the following years, the fleet size slightly decreases due to the higher battery capacities of BEBs purchased in the final stage.} 
The transition to ZEBs is almost complete by year 18 even under the Slow-Slow scenario. Therefore, the zero-emission target can be easily achieved and the adoption to BEBs is also economically justified  due to their low O\&M costs compared to DBs and HFCBs.
Note that HFCBs are only purchased in the last few time periods of the planning horizon, and only in scenarios  where the HFCB technology improves fast and BEB technology improves slowly. Even in the most favorable scenario for HFCBs,  less than 0.4\% of the fleet will be HFCBs in the end.

In Figure~\ref{fig:fleetEvolution},
We present the bus fleet composition for the first year of Stage 2 through Stage 5, corresponding to years 6, 11, 16, and 21, across two extreme scenarios: Fast-Fast and Slow-Fast. (The results of Fast-Slow and Slow-Slow are omitted, as they are almost identical to those of Fast-Fast and Slow-Fast, respectively.)
This figure shows the number of buses of each technology assigned to each cluster during the winter of those years.
Recall that Cluster 1 has a DSR value of 1, meaning that it can be electrified with the same number of DBs, while larger cluster labels indicate smaller DSR values, which correspond to routes that may require additional BEBs compared to DBs.
The results show that in all  scenarios, the transition to BEBs does not begin for clusters 5 through 12 within the first 6 years of the planning period due to their lower DSR values. Additionally, the last four clusters are still using diesel buses during the first 11 years. Cluster 4 is electrified more in the first 6 years than the previous three clusters, likely due to greater savings in O\&M costs and the higher age of the initial fleet corresponding to these clusters, among other factors.
{ We also observe that the few HFCBs purchased in the Slow-Fast scenario are bought in years 22-24 and are almost entirely assigned to cluster 12 to reduce the fleet size in that cluster.}

\subsection{Sensitivity Analysis}
\label{sec:case-sensitivity}

We now present the results of our extensive sensitivity analyses, in which we change some key deterministic parameters of the multi-stage stochastic program that are exogenously determined such as the  budget restrictions, emission targets and hydrogen prices.
The choice of these parameters are also motivated by our findings from the {\tt{Base Case}}: 
i) Budget restrictions are influencing the initial investment decisions, and forcing the model to choose DBs over BEBs due to the lower investment cost of the former despite the lower lifetime cost of the latter. 
ii) Emission targets are mostly redundant, which enables the model to make the best economical decisions that happen to address the environmental concerns as well.
iii) HFCBs are quite expensive in terms of their O\&M costs even without the consideration of potential infrastructure costs, which results in BEBs being the prominent choice for the transition.

\newcounter{groupcount}
\pgfplotsset{
    draw group line/.style n args={5}{
        after end axis/.append code={
            \setcounter{groupcount}{0}
            \pgfplotstableforeachcolumnelement{#1}\of\datatable\as\cell{%
                \def\temp{#2}
                \ifx\temp\cell
                    \ifnum\thegroupcount=0
                        \stepcounter{groupcount}
                        \pgfplotstablegetelem{\pgfplotstablerow}{[index]0}\of\datatable
                        \coordinate [yshift=#4] (startgroup) at (axis cs:\pgfplotsretval,0);
                    \else
                        \pgfplotstablegetelem{\pgfplotstablerow}{[index]0}\of\datatable
                        \coordinate [yshift=#4] (endgroup) at (axis cs:\pgfplotsretval,0);
                    \fi
                \else
                    \ifnum\thegroupcount=1
                        \setcounter{groupcount}{0}
                        \draw [
                            shorten >=-#5,
                            shorten <=-#5
                        ] (startgroup) -- node [anchor=north] {#3} (endgroup);
                    \fi
                \fi
            }
            \ifnum\thegroupcount=1
                        \setcounter{groupcount}{0}
                        \draw [
                            shorten >=-#5,
                            shorten <=-#5
                        ] (startgroup) -- node [anchor=north] {#3} (endgroup);
            \fi
        }
    }
}

\begin{figure}[H]
    \pgfplotstableread{
    1	457	41	0	1
    2	470	27	0	1
    3	200	297	0	1
    4	0	498	0	1
    5	350	16	0	2
    6	308	58	0	2
    7	89	280	0	2
    8	0	370	0	2
    9	366	182	0	3
    10	230	322	0	3
    11	17	543	0	3
    12	0	552	0	3
    13	501	616	0	4
    14	86	1036	0	4
    15	0	1127	0	4
    16	0	1119	0	4
    17	994	0	0	5
    18	63	991	0	5
    19	0	1048	0	5
    20	0	1040	0	5
    21	1493	0	0	6
    22	0	1600	0	6
    23	0	1575	0	6
    24	0	1555	0	6
    25	718	0	0	7
    26	82	694	0	7
    27	0	780	0	7
    28	0	777	0	7
    29	497	0	0	8
    30	403	104	0	8
    31	0	551	0	8
    32	0	537	0	8
    33	196	0	0	9
    34	196	1	0	9
    35	0	228	0	9
    36	0	219	0	9
    37	106	0	0	10
    38	106	0	0	10
    39	0	127	0	10
    40	0	125	0	10
    41	41	0	0	11
    42	41	0	0	11
    43	0	55	0	11
    44	0	54	0	11
    45	22	0	0	12
    46	21	1	0	12
    47	18	8	0	12
    48	0	29	0	12
    }\datatable

    \begin{subfigure}{0.32\textwidth}
        \centering
        \begin{tikzpicture}[scale=0.85]
            \pgfplotsset{every axis y label/.append style={rotate=0,yshift=0.0cm}}
            \begin{axis}[
                title={Fast-Fast},
                ylabel=label,
                xtick=data,
                width=18cm, height=8cm,
                xticklabels={\scriptsize 6 , \scriptsize 11, \scriptsize 16, \scriptsize 21, \scriptsize 6 , \scriptsize 11, \scriptsize 16, \scriptsize 21, \scriptsize 6 , \scriptsize 11, \scriptsize 16, \scriptsize 21, \scriptsize 6 , \scriptsize 11, \scriptsize 16, \scriptsize 21, \scriptsize 6 , \scriptsize 11, \scriptsize 16, \scriptsize 21, \scriptsize 6 , \scriptsize 11, \scriptsize 16, \scriptsize 21, \scriptsize 6 , \scriptsize 11, \scriptsize 16, \scriptsize 21, \scriptsize 6 , \scriptsize 11, \scriptsize 16, \scriptsize 21, \scriptsize 6 , \scriptsize 11, \scriptsize 16, \scriptsize 21, \scriptsize 6 , \scriptsize 11, \scriptsize 16, \scriptsize 21, \scriptsize 6 , \scriptsize 11, \scriptsize 16, \scriptsize 21, \scriptsize 6 , \scriptsize 11, \scriptsize 16, \scriptsize 21},
                enlarge y limits=false,
                enlarge x limits=0.01,
                ymin=0,
                ymax= 1700,
                ybar stacked,
                bar width=4pt,
                ylabel={\#Buses Assigned},
                draw group line={[index]4}{1}{\scriptsize \tt{Cluster1}}{-2.8ex}{1pt},
                draw group line={[index]4}{2}{\scriptsize \tt{Cluster2}}{-2.8ex}{2pt},
                draw group line={[index]4}{3}{\scriptsize \tt{Cluster3}}{-2.8ex}{2pt},
                draw group line={[index]4}{4}{\scriptsize \tt{Cluster4}}{-2.8ex}{2pt},
                draw group line={[index]4}{5}{\scriptsize \tt{Cluster5}}{-2.8ex}{2pt},
                draw group line={[index]4}{6}{\scriptsize \tt{Cluster6}}{-2.8ex}{2pt},
                draw group line={[index]4}{7}{\scriptsize \tt{Cluster7}}{-2.8ex}{2pt},
                draw group line={[index]4}{8}{\scriptsize \tt{Cluster8}}{-2.8ex}{2pt},
                draw group line={[index]4}{9}{\scriptsize \tt{Cluster9}}{-2.8ex}{2pt},
                draw group line={[index]4}{10}{\scriptsize \tt{Cluster10}}{-2.8ex}{2pt},
                draw group line={[index]4}{11}{\scriptsize \tt{Cluster11}}{-2.8ex}{2pt},
                draw group line={[index]4}{12}{\scriptsize \tt{Cluster12}}{-2.8ex}{2pt}
            ]
            \addplot [ybar, fill=gray!90, draw=black!70] table[x index=0,y index=1] \datatable;
            \addplot [ybar, fill=gray!30, draw=black!70] table[x index=0,y index=2] \datatable;
            \addplot [ybar, fill=black, draw=black] table[x index=0,y index=3] \datatable;
            \legend{DB, BEB, HFCB}
            \end{axis}
        \end{tikzpicture}
    \end{subfigure}
    \hfill
    \pgfplotstableread{
    1	455	41	0	1
    2	471	26	0	1
    3	200	297	0	1
    4	6	491	0	1
    5	365	0	0	2
    6	310	56	0	2
    7	88	282	0	2
    8	0	370	0	2
    9	540	0	0	3
    10	472	69	0	3
    11	49	507	0	3
    12	0	557	0	3
    13	599	517	0	4
    14	148	980	0	4
    15	0	1124	0	4
    16	0	1118	0	4
    17	995	0	0	5
    18	971	26	0	5
    19	0	1053	0	5
    20	0	1046	0	5
    21	1493	0	0	6
    22	41	1568	0	6
    23	0	1596	0	6
    24	0	1560	0	6
    25	718	0	0	7
    26	82	704	0	7
    27	0	790	0	7
    28	0	781	0	7
    29	496	0	0	8
    30	446	53	0	8
    31	8	546	0	8
    32	0	542	0	8
    33	197	0	0	9
    34	195	1	0	9
    35	0	231	0	9
    36	0	230	0	9
    37	106	0	0	10
    38	106	0	0	10
    39	0	128	0	10
    40	0	128	0	10
    41	41	0	0	11
    42	41	0	0	11
    43	43	0	0	11
    44	0	55	0	11
    45	22	0	0	12
    46	21	1	0	12
    47	21	1	0	12
    48	0	36	0	12
    }\datatableSF

    \begin{subfigure}{0.32\textwidth}
        \centering
        \begin{tikzpicture}[scale=0.85]
            \pgfplotsset{every axis y label/.append style={rotate=0,yshift=0.0cm}}
            \begin{axis}[
                            title={Slow-Fast},
                ylabel=label,
                xtick=data,
                width=18cm, height=8cm,
                xticklabels={\scriptsize 6 , \scriptsize 11, \scriptsize 16, \scriptsize 21, \scriptsize 6 , \scriptsize 11, \scriptsize 16, \scriptsize 21, \scriptsize 6 , \scriptsize 11, \scriptsize 16, \scriptsize 21, \scriptsize 6 , \scriptsize 11, \scriptsize 16, \scriptsize 21, \scriptsize 6 , \scriptsize 11, \scriptsize 16, \scriptsize 21, \scriptsize 6 , \scriptsize 11, \scriptsize 16, \scriptsize 21, \scriptsize 6 , \scriptsize 11, \scriptsize 16, \scriptsize 21, \scriptsize 6 , \scriptsize 11, \scriptsize 16, \scriptsize 21, \scriptsize 6 , \scriptsize 11, \scriptsize 16, \scriptsize 21, \scriptsize 6 , \scriptsize 11, \scriptsize 16, \scriptsize 21, \scriptsize 6 , \scriptsize 11, \scriptsize 16, \scriptsize 21, \scriptsize 6 , \scriptsize 11, \scriptsize 16, \scriptsize 21},
                enlarge y limits=false,
                enlarge x limits=0.01,
                ymin=0,
                ymax= 1700,
                ybar stacked,
                bar width=4pt,
                ylabel={\#Buses Assigned},
                draw group line={[index]4}{1}{\scriptsize \tt{Cluster1}}{-2.8ex}{1pt},
                draw group line={[index]4}{2}{\scriptsize \tt{Cluster2}}{-2.8ex}{2pt},
                draw group line={[index]4}{3}{\scriptsize \tt{Cluster3}}{-2.8ex}{2pt},
                draw group line={[index]4}{4}{\scriptsize \tt{Cluster4}}{-2.8ex}{2pt},
                draw group line={[index]4}{5}{\scriptsize \tt{Cluster5}}{-2.8ex}{2pt},
                draw group line={[index]4}{6}{\scriptsize \tt{Cluster6}}{-2.8ex}{2pt},
                draw group line={[index]4}{7}{\scriptsize \tt{Cluster7}}{-2.8ex}{2pt},
                draw group line={[index]4}{8}{\scriptsize \tt{Cluster8}}{-2.8ex}{2pt},
                draw group line={[index]4}{9}{\scriptsize \tt{Cluster9}}{-2.8ex}{2pt},
                draw group line={[index]4}{10}{\scriptsize \tt{Cluster10}}{-2.8ex}{2pt},
                draw group line={[index]4}{11}{\scriptsize \tt{Cluster11}}{-2.8ex}{2pt},
                draw group line={[index]4}{12}{\scriptsize \tt{Cluster12}}{-2.8ex}{2pt}
            ]

            \addplot [ybar, fill=gray!90, draw=black!70] table[x index=0,y index=1] \datatableSF;
            \addplot [ybar, fill=gray!30, draw=black!70] table[x index=0,y index=2] \datatableSF;
            \addplot [ybar, fill=black, draw=black] table[x index=0,y index=3] \datatableSF;
            \end{axis}
        \end{tikzpicture}
    \end{subfigure}

    \caption{Fleet composition in years 6, 11, 16, 21 under different scenarios with respect to route clusters.}
    \label{fig:fleetEvolution}
\end{figure}
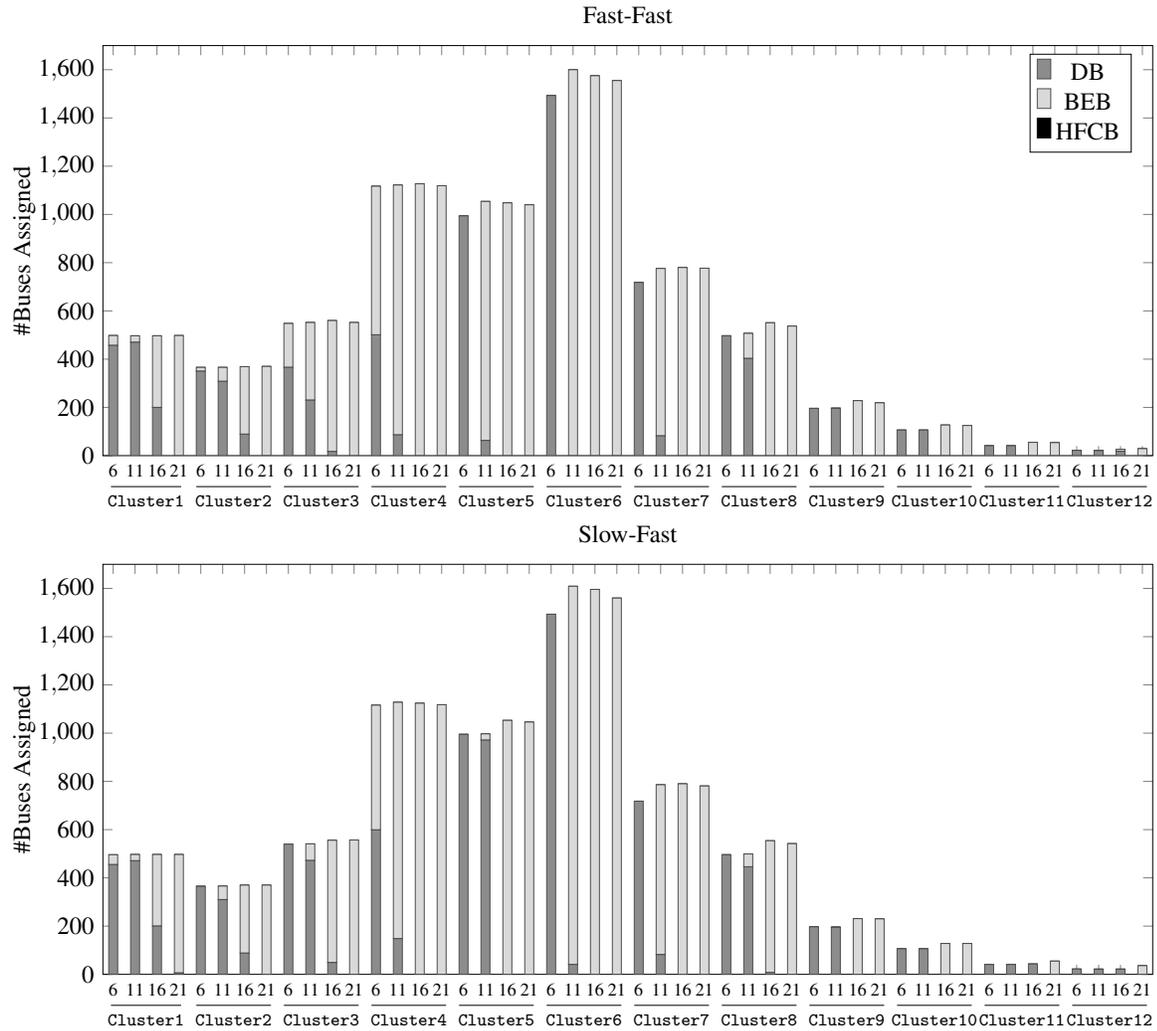

These observations motivate us to explore the answers of the following questions: 
i) If IETT can spare more budget annually, how much cost  benefit can be expected overall? 
ii) If even more strict emission targets are enforced, how much cost  increase can be expected overall? 
ii) If hydrogen prices drop drastically, can HFCBs become the prominent ZEB choice and how much cost  decrease can be expected overall? 

To answer these questions, we re-run our multi-stage stochastic model three times for the following cases by keeping everything else the same: 
\begin{enumerate}[i)]
    \item 
 {\tt{Relaxed Budget}}: Budget is chosen 300 Million USD annually.
    \item 
 {\tt{Strict Emissions}}: Zero emission target is enforced in year 2040 with intermediate targets starting from year 2030.
    \item 
 {\tt{Low Hydrogen Price}}: Hydrogen price is set at 2 USD/\unit{\kg}.
\end{enumerate}
In Table~\ref{tab:CaseResults}, we provide a summary  of the sensitivity analysis conducted while the detailed results are available in the Supplementary Material. We observe that ``\%Gap", defined as $100\times \frac{z_{\text{round}}-z_{\text{LP}}}{z_{\text{round}}}$, is consistently small, which suggests that our rounding scheme is quite successful in the new cases as well. In addition, we also report ``\%Change", which is the percentage change of a new case with respect to the {\tt{Base Case}}  in terms of the objective function value of the rounded solution. As expected, in the {\tt{Relaxed Budget Case}}, the change is negative, which indicates that an increased budget helps to obtain a more economical plan with 3.80\% lower expected cost. Interestingly, the change in the {\tt{Strict Emission Case}}  is very small, which suggests that an even earlier transition to ZEBs can be achieved with increasing the overall cost marginally by 0.08\%. Finally, the change in the expected total cost for the {\tt{Low Hydrogen Price Case}}  is  also small (-0.17\%), which shows the limited effect the hydrogen price has on the overall outcome.
\begin{table}[!ht]
\caption{Summary results for different cases.}
    \centering
    \setlength{\tabcolsep}{0.1pt}
    \small
    \begin{tabular}{|c|c|c|c|c|c|} 
    \hline
 \textbf{Case}   & \textbf{Time (sec)} & $z_{\text{LP}}$ & $z_{\text{round}}$  & \textbf{\%Gap} & \textbf{\%Change} \\
    \hline
    \tt{Base} &  \text{} 5,435 \text{} &  14,126,110,527  &  14,154,103,467  & 0.20 
    &  -  \\
    \hline
    \tt{Relaxed Budget} &  5,888  &  13,586,698,671  &  13,615,884,770 & 0.21 & -3.80 
    \\
    \hline
    \tt{Strict Emission} &  4,877  &  14,133,211,491  &  14,165,320,511  &  0.23 & \ 0.08  
    \\
    \hline
    \tt{Low Hydrogen Price} &  6,133  &  14,099,073,373  &  14,129,697,437  & 0.22 & -0.17 
    \\
    \hline
    \end{tabular}
    
    \label{tab:CaseResults}
\end{table}

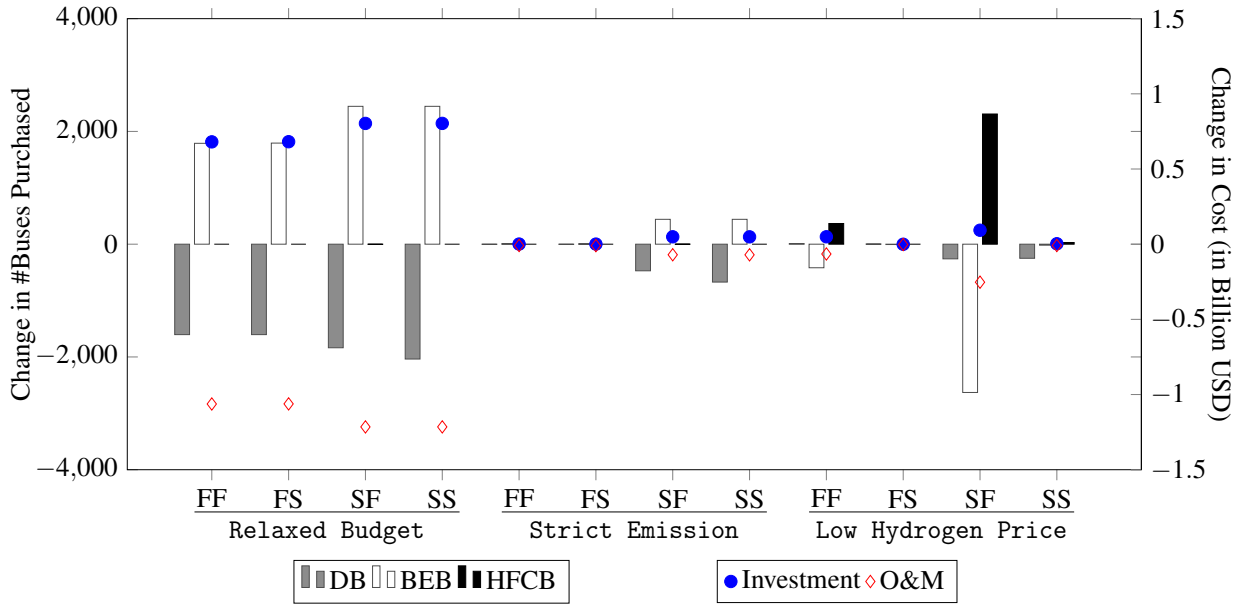
\begin{figure}[H]
\centering
\begin{tikzpicture}[scale=0.935]
      \pgfplotsset{every axis y label/.append style={rotate=0,yshift=0.25cm}}

\begin{axis}[
        ybar,
        bar width=6pt,
        width=16cm, height=8cm,
    ylabel=label,
    xtick=data,
    xticklabels={FF, FS, SF, SS, FF, FS, SF, SS, FF, FS, SF, SS},
    enlarge y limits=false,
    enlarge x limits=0.1,
    ymin=-4000,
    ymax= 4000,
    ylabel={Change in \#Buses Purchased},
    ybar,
    legend style={at={(0.3,-0.20)}, anchor=north,legend columns=-1},
    draw group line={[index]6}{1}{\tt{Relaxed Budget}}{-22ex}{7pt},
    draw group line={[index]6}{2}{\tt{Strict Emission}}{-22ex}{7pt},
    draw group line={[index]6}{3}{\tt{Low Hydrogen Price}}{-22ex}{7pt}
]
\addplot[fill=gray!90,draw=black!70] table[ybar, x index=0,y index=1] \datatable;
\addplot[fill=gray!0,draw=black!70] table[ybar, x index=0,y index=2] \datatable;
\addplot[fill=black] table[ybar, x index=0,y index=3] \datatable;
\addplot[fill=white,draw=white] table[ybar, x index=0,y index=6] \datatable;  
\legend{DB, BEB, HFCB}
\end{axis}

      \pgfplotsset{every axis y label/.append style={rotate=180,yshift=17cm}}

    \begin{axis}[
    ylabel={Change in Cost (in Billion USD)},
    bar width=6pt,
        width=16cm, height=8cm,
    xticklabels={FF, FS, SF, SS, FF, FS, SF, SS, FF, FS, SF, SS},
    axis x line=none,
    enlarge y limits=false,
    enlarge x limits=0.1,
    ymin=-1.5,ymax=1.5,
      axis y line*=right,
    legend style={at={(0.7,-0.20)}, anchor=north,legend columns=-1},
]

\addplot[mark=*,blue, mark size=2.5pt, only marks] table[x index=0,y index=4] \datatable;

\addplot[mark=diamond,red, mark size=2.5pt, only marks] table[x index=0,y index=5] \datatable;
\legend{Investment, O\&M}

\end{axis}

\end{tikzpicture}
\caption{Changes in bus purchases (shown with bars) and cost values (shown with markers) compared to the {\tt{Base Case}} for four specific scenarios in different cases. FF, FS, SF,  SS are respectively the abbreviations of Fast-Fast, Fast-Slow, Slow-Fast, Slow-Slow scenarios as described earlier.}
\label{fig:sensitivity-fourScenarios}
\end{figure}

We also report the changes in the fleet decomposition and cost figures for the new cases with respect to the {\tt{Base Case}} under four specific scenarios in Figure~\ref{fig:sensitivity-fourScenarios}.

\pgfplotstableread{
1	-1604	1789	0	0.680533317	-1.062876264	1
2	-1604	1793	0	0.681542463	-1.062220771	1
3	-1837	2445	-1	0.803139901	-1.215351023	1
4	-2037	2445	0	0.803326356	-1.215209635	1
5	-2	8	0	0.000810033	-0.008539448	2
6	-2	6	0	0.000523587	-0.009555154	2
7	-472	441	1	0.049047984	-0.070316511	2
8	-672	441	0	0.048784767	-0.070793603	2
9	11	-418	367	0.049253014	-0.065149238	3
10	3	-6	0	-0.000582901	-0.002340276	3
11	-260	-2630	2309	0.092785794	-0.252809884	3
12	-250	-18	30	0.002843665	-0.008701421	3
}\datatable

In the {\tt{Relaxed Budget Case}},  we observe an increase in the number of BEBs purchased across all scenarios as BEBs are economically more advantageous than DBs in their economic lifetime despite their higher investment cost. This is clear from the moderate increase in the investment cost and the substantial decrease in the O\&M cost, which is even more pronounced in the Slow-Fast and Slow-Slow scenarios.

In the {\tt{Strict Emission Case}}, slightly more BEBs are being purchased, with the difference only being non-negligible in scenarios when BEBs improve slowly. This is due to the need to purchase BEBs between 2030 and 2040 to meet emission constraints, even in scenarios where improvements are slow. 
The small changes with this respect to the {\tt{Base Case}}  illustrate that BEBs are not only environmentally benefical but also cost effective.

In the {\tt{Low Hydrogen Price Case}}, results indicate that the number of HFCBs purchased increases in most scenarios. In particular, for the four specific scenarios, HFCB purchases increase in three of them, but remain at zero in the Fast-Slow scenario. We also observe that HFCB purchases begin in year 6 in the Slow-Fast scenario, where HFCBs are expected to make up 24.4\% of the fleet by the final time period. Nevertheless, as the expected total cost reduction is very small overall, this case suggests that HFCBs will not be dominantly more preferable over BEBs, even if we ignore potential infrastructure costs and potential difficulties regarding hydrogen supply chain.

\subsection{Model Comparison}
\label{sec:case-comparison}

\subsubsection{Comparison with a Deterministic Model}

To analyze the effect of stochasticity, we solve a deterministic version of our problem and compare its solution with that of the multi-stage stochastic program. In the deterministic model, the technology improvements described in Section~\ref{sec:tech_imp} are represented as a single scenario, reflecting the average of all historical advancements. Specifically, the anticipated cost multipliers for BEBs and HFCBs are 0.48 and 0.74, respectively while the efficiency multipliers are 1.26 and 1.08, respectively. The model solves in 6.4 seconds, and the objective function values are obtained as $z_{\text{LP}} = 13.233$B USD and $z_{\text{round}} = 13.256$B USD. 
We observe that the deterministic problem's $z_{\text{round}}$ value is 0.91B USD (6.43\%) lower than that of the stochastic program. This indicates that the solution of deterministic problem cannot be feasible for the stochastic program and we quantify these infeasibilities as described below.

First, we  apply the solution of the deterministic problem to each scenario of the stochastic program. To do this, we adjust the battery capacities of BEBs in each scenario based on the capacities from the deterministic solution. However, since energy density improvements vary across scenarios, we need to adjust the bus weights accordingly to reflect the expected improvements in each stochastic scenario. Then, using equation~\ref{eq:IC-BEB-costEvolution}, we calculate the BEB investment costs, where $C_n$ (the battery capacity in node $n$) comes from the deterministic solution, and $\vartheta^c_{\text{BEB}, n}$ (the efficiency improvement) relates to each specific scenario. We observe that  the deterministic solution violates the budget constraints in 160 out of 256 scenarios of the stochastic program in at least one time period by more than 1\%. Considering the probabilities of these scenarios, this means that there is a 67.41\% probability that the deterministic solution will not be feasible due to insufficient budget. Note that these budget violations are even underestimated because, in scenarios where energy density improves more slowly than in the deterministic case, heavier batteries are required, necessitating more buses to meet demand.

Table~\ref{tab:budgetViolations} details the budget violations of the deterministic solution in the stochastic program.
Since no HFCBs are purchased in the deterministic problem, it suffices to consider a reduced scenario tree with $2^{5-1}=16$ scenarios, in which only  BEB improvements are considered. In particular, scenarios in Table~\ref{tab:budgetViolations} represent BEB improvements from stage 2 to stage 5 (e.g., FSSF refers to fast improvement in stage 2, slow in stages 3 and 4, and fast again in stage 5). The table shows the percentage of budget violations under different scenarios, and for each year within the stage. For example, the numbers 6.11, 6.76, 6.17, 4.87, 5.15 correspond to the percentage budget violations in years 6, 7, 8, 9, 10 respectively, and under the eight specified BEB improvement scenarios.
The observed infeasibilities highlight the importance of stochastic modelling and dynamic planning in strategic decision making. By accounting for uncertainty in technology improvements, the stochastic model avoids the infeasibilities that arise when applying the deterministic solution, ensuring a more robust and feasible plan.

\begin{table}[htbp]
\centering
    \caption{Percent violations of budget constraints in different stages under different scenarios. A non-empty cell gives the yearly violations in that stage while the symbol `-' indicates that there is no violation in the corresponding year. An empty cell indicates that there is no violation  in that stage.}
    \label{tab:budgetViolations}
\scriptsize
\begin{tabular}%
{|c|c|c|c|c|c|c|c|c|c|c|c|c|c|c|c|c|}
\hline
Stage 2  & \multicolumn{8}{c|}{} & \multicolumn{8}{c|}{6.11,6.76,6.17,4.87,5.15}  \\ \hline
Stage 3 & \multicolumn{4}{c|}{} & \multicolumn{4}{c|}{} & \multicolumn{4}{c|}{} &   \multicolumn{4}{c|}{11.96,11.98,-,11.43,9.39} \\ \hline
Stage 4 
& \multicolumn{2}{c|}{} & \multicolumn{2}{c|}{} & \multicolumn{2}{c|}{} &\multicolumn{2}{c|}{1.25,1.42,1.44,1.40,-}  & \multicolumn{2}{c|}{} & \multicolumn{2}{c|}{1.25,1.42,1.44,1.40,-} & \multicolumn{2}{c|}{1.25,1.42,1.44,1.40,-}  & \multicolumn{2}{c|}{11.32,12.23,12.36,12.18,-} \\ \hline
Stage 5 & & & & & & & & {2.60,-,-,-,-} & & & & {2.60,-,-,-,-} & & {2.60,-,-,-,-} & {2.60,-,-,-,-} & {9.32,-,-,-,-} \\ \hline
& FFFF & FFFS & FFSF & FFSS & FSFF & FSFS & FSSF & FSSS & SFFF & SFFS & SFSF & SFSS & SSFF & SSFS & SSSF & SSSS \\ \hline
\end{tabular}
\end{table}

\subsubsection{Comparison with  Other Stochastic Models}

To better understand the impact of scenarios on the optimal solution, we extend the scenario tree and solve the stochastic model using two extended versions. In the first version called the {\tt{3-by-2 Case}}, we include 3 branches for BEB improvements at each stage (corresponding to Figure~\ref{fig:scatter_plot_3_clusters}) and 2 branches for HFCBs. The second version called the {\tt{2-by-3 Case}} maintains 2 branches for BEBs and extends the number of  HFCB branches to 3 (corresponding to Figure~\ref{fig:scatter_plot_3_clusters-fc}). In both versions, the number of scenarios increases to $(1\times2\times3)^{5-1}=1296$, making the problem size significantly larger complex. To manage this, we simplify certain assignment decisions. Specifically, for route-bus length pairs with higher demand in the winter schedule, assignment variables for other seasons are excluded. Instead, we assume that buses assigned to a cluster remain there across all seasons, with their usage adjusted based on seasonal demand, to account for lower O\&M costs during other seasons. The detailed results are reported in the Supplementary Material.

We observe that in both extended cases the overall transition plan is similar to the {\tt{Base Case}} as BEBs are still the dominant choice. However, in the Fast-Fast and Fast-Slow scenarios of the {\tt{3-by-2 Case}}, the transition is nearly complete in year 14, as the Fast improvement branch of BEBs is faster in the {\tt{3-by-2 Case}} than the {\tt{Base Case}}. We also observe that HFCB purchases are almost identical in the {\tt{3-by-2 Case}} to the {\tt{Base Case}}, meaning only a few buses being purchased in the final stage and under the Slow-Fast scenario, out of the 4 extreme scenarios. This is because the Slow improvement branch of BEBs in the {\tt{3-by-2 Case}} is not much different than that of the Base Case. 

The {\tt{2-by-3 Case}} results are nearly identical to the {\tt{Base Case}}. We observe that in the four extreme scenarios, only in the Slow-Fast scenario HFCBs are purchased, and the number of purchases is identical to the Slow-Fast scenario of the {\tt{Base Case}}.

\section{Conclusion}
\label{sec:conclusion}

In this study, we propose a multi-stage stochastic program for the bus fleet replacement problem where we consider the uncertainties in the technology advancements and cost improvements of different bus technologies. We present a forecasting approach where we examine historical cost and efficiency trends of different technology options and compute the improvement ratios over the years. These ratios are then clustered into different groups which are provided to the multi-stage stochastic program as scenarios. We use our model to plan the transition of municipal bus network in Istanbul to clean alternatives (BEBs and HFCBs) over a planning horizon of 25 years with a zero-emission target in 2050. In our plan, we consider the changes in the energy consumption of BEBs on different routes and at different seasons, and we ensure operational feasibility by finding the diesel-electric replacement ratios and recharge scheduling for FC electric buses. We perform sensitivity analyses to test the effect of certain exogenously determined parameters and inputs on our results. We also compare our stochastic programming approach with a deterministic model, and show its advantage in providing feasible dynamic strategic transition plans.

Our results from the multi-stage stochastic program indicate that BEBs are viable alternatives to replace DBs to achieve zero emission bus fleet goals. BEBs can already satisfy the demand in most of the routes in Istanbul without the need to plan new bus lines or timetables even at today's technology level. Their lower O\&M cost compared to DBs make them advantageous in the long term across all scenarios we implemented. Although DBs are needed to be purchased initially due to the old age of the fleet and higher investment costs of BEBs, as costs improve BEBs can replace all DBs and transition to an all-electric fleet can be completed in less than 20 years in most scenarios. The advantages of BEBs also manifest themselves in the sensitivity analyses we performed. When the available budget is higher, more BEBs are purchased accross all scenarios, proving that once initial investment barrier is overcome, BEBs are better alternatives to DBs over their lifetime. We also observe that  under tighter emission constraints the number of BEBs increase only negligibly under most scenarios, showing BEBs are not only environmentally beneficial but
also cost effective. We acknowledge that other adoption challenges such as  potential need for grid investments, recharging scheduling and personnel training can pose certain barriers and slow down the transition in practice.

HFCBs do not prove to be preferable over BEBs across our scenarios- even under the most favorable scenario, less than 0.4\%  of the fleet consists of HFCBs. Higher investment and operational costs both play a role in this result. When the energy costs are decreased by one fifth, the fleet can consist of up to 24\% of HFCBs in the most favorable scenario but HFCBs still cannot be dominant due to their higher initial costs. 

There are some limitations of our work that should be considered while interpreting our results. In this study, we assumed that 
the municipality is not responsible for building hydrogen transmission and~storage infrastructure, or potential grid investments necessary for increased electricity demand.
In a setting where the municipality has to invest in infrastructure, transition to cleaner alternatives might be slower. In addition, our results are based on scenarios obtained using limited historical HFCB cost and technological development trends compared to BEBs. According to these trends, even slow  improvement in BEBs is faster than in HFCBs. However, our scenarios assume constant improvements based on past data and HFCBs are still at their early adoption stages being produced at limited numbers. It is possible that BEBs will reach a saturation in improvement sooner than HFCBs, which might change the results in favor of HFCBs in the~future. 

While this study focuses on strategical planning of the transition to zero emission buses, it can be enhanced by adding more tactical and operational aspects. For example, bus-to-route simulator might be improved by including recharge and maintenance scheduling optimization as well as battery degradation aspects. In our study, we assume each bus is changed with a bus of the same length, trip schedules do not change with transition and chargers are available at the last stop of  every route when needed. Each of these assumptions pose an optimization problem by themselves: optimal bus size selection and assignment problem, optimizing trip schedules and charge location optimization. Our future work will focus on these problems and  investigate how to integrate them into our model.

Finally, we will also investigate how to solve the large-scale multi-stage stochastic program more efficiently. This is especially important if one needs to include more technology options and construct larger technology trees that might provide a more realistic case study, albeit with a significantly larger scenario tree. In this case, obtaining structural results to systematically eliminate variables and constraints, and developing decomposition methods may prove to be necessary.

\subsection*{Acknowledgements}
This work was supported by the  Scientific and Technological Research Council of Turkey [grant number 222M243]. The authors thank IETT for providing detailed datasets related to the Istanbul public bus network, and Ahmet Emir {\c{S}}ener for his help related to Section~\ref{sec:tech_imp-beb}.

{
\printbibliography
}

\appendix

\section{Energy Requirement Calculation}

\subsection{Determining the Segment Duration}

To find $\tau_1$ and $\tau_2$, we can solve an optimization problem that minimizes the time spent within a segment of length $d$ and a fixed acceleration $a$ satisfying the condition that $V(\tau) \le V_{\max}$ and  
$$P_w(\tau) = \left(m g \sin(\alpha) + f_r m g \cos(\alpha) + 0.5 \rho_{\text{air}} C_D A_f v(\tau)^2 + m_{\text{eq}} a(\tau)\right) v(\tau) \le\min\{P_{\max}, \kappa V(\tau)\}$$ for $\tau \in [0, \tau_1+\tau_2]$.
To simplify, let us rewrite $P_w(\tau)$ as
\[
    P_w(\tau) = \left(A + B v(\tau)^2 + C a(\tau)\right) v(\tau).
\]
Since the air resistance is much smaller than the other components (that is, $B v(\tau)^2 \ll |A +C a(\tau)|$ for the bus specifications and speed profiles we are interested in), we can formulate this optimization problem as follows when $\alpha > 0$ (for $\alpha<0$, the power constraint is redundant):
\[
\min \{ \tau_1 + \tau_2 : a \tau_1 \le V_{\max}, \ a \tau_1 \tau_2 = d, \ \tau_2 \ge \tau_1, \ P_w(\tau_1) \le \min\{P_{\max} , \kappa a \tau_1) \} \}.
\]
We can eliminate $\tau_2$ by substituting $\frac{d}{a\tau_1}$:
\[
\min \left\{ \tau_1 + \frac{d}{a\tau_1} : a \tau_1 \le V_{\max}, \ \frac{d}{a\tau_1} \ge \tau_1, \ P_w(\tau_1) \le P_{\max}, \tau \le \sqrt{\frac{\kappa - A - Ca}{B a^2}} \right\}.
\]
Note that $\tau_1 + \frac{d}{a\tau_1}$ is a convex function with its minimizer at $  \tau_1 = \sqrt{\frac{d}{a}} $. Since the cubic polynomial $P_w(\tau_1)$ is increasing under our assumptions, there exists a single root for $P_w(\tau_1)=P_{\max}$, which we denote by $\tau_R$ (note that we can search for that root in the interval $ [0, \frac{V_{\max}}{a}]$).  Finally, we conclude that the optimal solution is
\[
\tau_1^* =   \min\left\{ \sqrt{\frac{d}{a}}, \frac{V_{\max}}{a},  \sqrt{\frac{\kappa - A - Ca}{B a^2}}, \tau_R \right\}  
\text{ and } \tau_2^* = \frac{d}{a \tau_1^*}.
\]

\subsection{Calculating the Energy Consumption for BEBs}

\begin{itemize}
    \item  {In the acceleration phase}, that is $\tau\in[0,\tau_1]$, we have $a(\tau)=a$, $v(\tau)=a\tau$, and we are interested in the following integral:
    \[
\int_{0}^{\tau_1} P_{\text{w}}(\tau) d\tau =
\int_{0}^{\tau_1} (A+Ba^2\delta^2+Ca)a\tau d \tau.
\]
\begin{itemize}
    \item If   $A+Ca \ge 0$, we have
\[
\begin{split}
\int_{0}^{\tau_1} P_{\text{bat}}(\tau) d\tau & =
\frac{1}{\eta_t \eta_m} \int_{0}^{\tau_1} P_{w}(\tau) d\tau =
\frac{1}{\eta_t \eta_m} \int_{0}^{\tau_1}  \left(A + B (a\tau)^2 + C a\right) (a\tau) d\tau  \\ &=
\frac{1}{\eta_t \eta_m}  \left [\frac12(A+Ca)a\tau_1^2 + \frac14 Ba^3\tau_1^4\right ].    
\end{split}
\]
    \item If $A+Ca <
    0$ and $\bar \tau := 
    \sqrt{\frac{-(A+Ca)}{Ba^2}}  \ge \tau_1$, we have 
    \[
\int_{0}^{\tau_1} P_{\text{bat}}(\tau) d\tau =
{\eta_t \eta_m \eta_{rb}} \int_{0}^{\tau_1} P_{w}(\tau) d\tau =
{\eta_t \eta_m \eta_{rb}}  \left [\frac12(A+Ca)a\tau_1^2 + \frac14 Ba^3\tau_1^4\right ].
\]
    \item     
    If $A+Ca <
    0$ and $\bar \tau := 
    \sqrt{\frac{-(A+Ca)}{Ba^2}}  < \tau_1$, we have 
    \[
    \begin{split}
        \int_{0}^{\tau_1} P_{\text{bat}}(\tau) d\tau 
    &=
     {\eta_t \eta_m \eta_{rb}} \int_{0}^{\bar \tau} (A+Ba^2\delta^2+Ca)a\delta d \delta
    +  \frac{1}{\eta_t \eta_m}   \int_{\bar \tau}^{\tau_1} (A+Ba^2\delta^2+Ca)a\delta d \delta  
     \\
    & = {\eta_t \eta_m \eta_{rb}} \left [ \frac12(A+Ca)a \bar \tau^2 + \frac14 Ba^3 \bar \tau^4 \right ] 
    +   \frac{1}{\eta_t \eta_m}  \left [ \frac12(A+Ca)a  (\tau_1^2-\bar \tau^2) + \frac14 Ba^3 (\tau_1^4-\bar \tau^4) \right ] . 
    \end{split}
    \]
\end{itemize}

\item {In the constant speed phase}, that is $\tau\in[\tau_1,\tau_2]$, we have $a(\tau)=0$ and $v(\tau)=a\tau_1$, and we are interested in the following integral:
\[
 \int_{\tau_1}^{\tau_2} P_{w}(\tau) d\tau =
\int_{\tau_1}^{\tau_2}  \left(A + B (a\tau_1)^2 \right) (a\tau_1) d\tau  
= \left(A + B (a\tau_1)^2 \right) (a\tau_1) (\tau_2-\tau_1) =: \mathcal{C}.
\]

\begin{itemize}
    \item If $\mathcal{C} \ge 0$, we have
    \[
    \int_{\tau_1}^{\tau_2} P_{\text{bat}}(\tau) d\tau = \frac{1}{\eta_t \eta_m} \int_{\tau_1}^{\tau_2} P_{w}(\tau) d\tau = \frac{1}{\eta_t \eta_m}  \mathcal{C}.
    \]
    \item If $\mathcal{C}  < 0 $, we have
    \[
    \int_{\tau_1}^{\tau_2} P_{\text{bat}}(\tau) d\tau = {\eta_t \eta_m \eta_{rb}} \int_{\tau_1}^{\tau_2} P_{w}(\tau) d\tau = {\eta_t \eta_m \eta_{rb}}  \mathcal{C}.
    \]
\end{itemize}

\item{In the deceleration phase}, that is $\tau\in[\tau_2,\tau_1+\tau_2]$, we have $a(\tau)=-a$ and $v(\tau)=a(\tau_1+\tau_2-t)$. Let us first apply change of variables $\delta = \tau_1+\tau_2-\tau$ so that we are interested in the following integral:
\[
\int_{\tau_2}^{\tau_1+\tau_2} P_{\text{w}}(\tau) d\tau =
\int_{0}^{\tau_1} P_{\text{w}}(\tau_1+\tau_2-\delta) d \delta = 
\int_{0}^{\tau_1} (A+Ba^2\delta^2-Ca)a\delta d \delta.
\]

\begin{itemize}
    \item If
     $A - Ca  \ge 0 $, we have
\[
\int_{\tau_2}^{\tau_1+\tau_2} P_{\text{bat}}(\tau) d\tau 
=
\frac{1}{\eta_t \eta_m}  \int_{0}^{\tau_1} (A+Ba^2\delta^2-Ca)a\delta d \delta
  = \frac{1}{\eta_t \eta_m}  \left [ \frac12 (A-Ca)a \tau_1^2 + \frac14 Ba^3 \tau_1^4 \right ].
\]

    \item If $A - Ca  <  0$ and $\bar \delta := \sqrt{\frac{Ca-A}{Ba^2}} \ge \tau_1$, we have
\[
\int_{\tau_2}^{\tau_1+\tau_2} P_{\text{bat}}(\tau) d\tau 
=
 {\eta_t \eta_m \eta_{rb}}  \int_{0}^{\tau_1} (A+Ba^2\delta^2-Ca)a\delta d \delta
  =  {\eta_t \eta_m \eta_{rb}} \left [ \frac12 (A-Ca)a \tau_1^2 + \frac14 Ba^3 \tau_1^4 \right ].
\]

    \item 
    If $A - Ca  <  0$ and $\bar \delta := \sqrt{\frac{Ca-A}{Ba^2}} < \tau_1$, we have
\[
\begin{split}
    \int_{\tau_2}^{\tau_1+\tau_2} P_{\text{bat}}(\tau) d\tau 
&=
 {\eta_t \eta_m \eta_{rb}} \int_{0}^{\bar \delta} (A+Ba^2\delta^2-Ca)a\delta d \delta
+  \frac{1}{\eta_t \eta_m}   \int_{\bar \delta}^{\tau_1} (A+Ba^2\delta^2-Ca)a\delta d \delta  
 \\
& = {\eta_t \eta_m \eta_{rb}} \left [ \frac12(A-Ca)a \bar \delta^2 + \frac14 Ba^3 \bar \delta^4 \right ] 
+    \frac{1}{\eta_t \eta_m}   \left [ \frac12(A-Ca)a  (\tau_1^2-\bar\delta^2) + \frac14 Ba^3 (\tau_1^4-\bar \delta^4) \right ] . 
\end{split}
\]
\end{itemize}

\end{itemize}

For the energy requirement calculations, we use the parameters listed in Table~\ref{tab:traction_power_parameters} for BEBs. Vehicle mass $m$ includes the body mass $m_{\text body}$ and the energy capacity multiplied with the unit battery mass $m_{\text bat}$. 
\begin{table}[H]
\centering
\caption{List of parameters for traction power and battery power calculations (acceleration and vehicle body mass are given for 8, 10, 12, 18 \unit{\meter} bus groups, respectively).}
\label{tab:traction_power_parameters}
\begin{tabular}{|c|c||c|c|}
\hline
\textbf{Parameter} & \textbf{Value} & \textbf{Parameter} & \textbf{Value} \\
\hline
$f_r$ & 0.01 & $\rho_{\text{air}}$ & 1.225 \si{\kilogram/\cubic\meter} \\
$C_D$ & 0.7 & $g$ & 9.81 \si{\meter/\square\second} \\
$A_f$ & $0.85 \times 3.25 \times 2.55 = 7.04 \si{\square\meter} $ & $\eta_t$ & 0.9 \\
 $\eta_m$ & 0.9    & $\eta_{\text{rb}}$ & 0.25 \\
$m_{\text body}$ & 9.5, 15, 16, 25 tonnes &  $m_{\text bat}$ & 5 \si{\kg/\kilo\watt\hour} \\
$m_{eq}$ & $1.1 \times m$ & $a $ & 2.1, 1.8, 1.7, 1.5 \si{\meter/\square\second}   \\
\hline
\end{tabular}
\end{table}

\newpage
\section{Data}

\subsection{IETT Data Details}

Here is a snapshot of the ``Trip Schedule" dataset:
  \begin{figure}[H]
    \centering
    \includegraphics[width=1\linewidth]{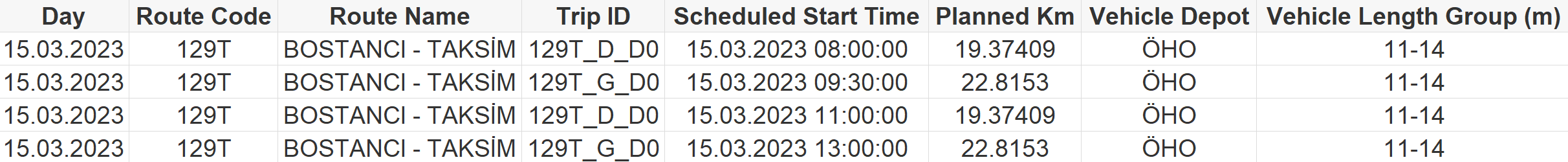}
    \caption{Snapshot of the ``Trip Schedule" dataset.}
    \label{fig:Trip_Snapshot}
    \end{figure}

    Summary statistics of the ``Trip Schedule" dataset are provided in Table~\ref{tab:BusSchedule}. 
\begin{table}[H]
\centering
\caption{Trip schedule summaries for March 15 and August 3, 2023. On March 15, service trips were planned for 2,863 route types with a total of 13,467 different bus stops visited, while on August 3, services were planned for 2,704 route types with 13,501 different bus stops visited.}

\label{tab:BusSchedule}
\small
\setlength{\tabcolsep}{0.1pt}
\begin{tabular}{|c|c|c|c|c|c|c|} 
\hline
\multirow{2}{*}{\textbf{Bus Length Group (\unit{\meter})}} & \multicolumn{3}{c|}{\textbf{Summer Schedule (August 3, 2023)}} & \multicolumn{3}{c|}{\textbf{Winter Schedule (March 15, 2023)}} \\
\cline{2-7}
& \textbf{Total Planned} & \textbf{Total Planned} & \textbf{Total Assigned} & \textbf{Total Planned} & \textbf{Total Planned} & \textbf{Total Assigned} \\
& \textbf{Services} & \textbf{Distance (\unit{\kilo\meter})} & \textbf{Routes} & \textbf{Services} & \textbf{Distance (\unit{\kilo\meter})} & \textbf{Routes} \\
\hline
6.5-8 & 2,898 & 63,525 & 38 & 2,948 & 61,200 & 37 \\
\hline
8-9 & 277 & 4,320 & 14 & 275 & 4,308 & 14 \\
\hline
10-11 & 76 & 3,086 & 8 & 64 & 2,635 & 6 \\
\hline
11-14 & 35,761 & 823,868 & 738 & 36,578 & 853,222 & 736 \\
\hline
11-14 natural gas & 1,464 & 26,448 & 129 & 1,554 & 30,263 & 145 \\
\hline
14-19 & 7,727 & 274,523 & 79 & 9,130 & 317,813 & 98 \\
\hline
\textbf{Total} & \textbf{48,203} & \textbf{1,192,684} & \textbf{820} & \textbf{50,549} & \textbf{1,269,441} & \textbf{819} \\
\hline
\end{tabular}
\end{table}

\noindent
    Here is a snapshot of the ``Stop Sequence with Coordinates" dataset:
    \begin{figure}[H]
        \centering
        \includegraphics[width=1\linewidth]{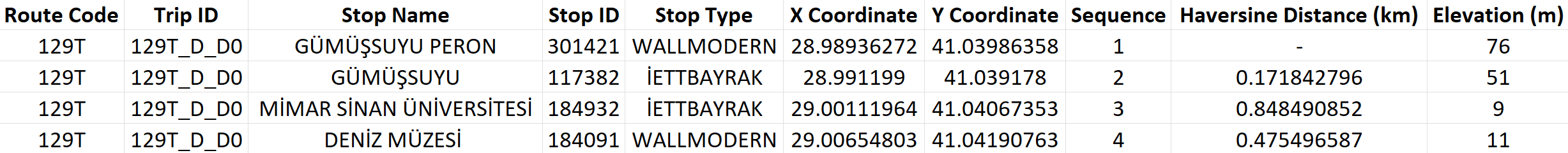}
        \caption{Snapshot of the ``Stop Sequence with Coordinates'' dataset.}
        \label{fig:Stop_Snapshot}
    \end{figure}

\newpage
\subsection{Technological Change}

\subsubsection{BEB Technological Change}

\begin{table}[h!]
\centering
\small
\caption{Technological change for BEBs over the years (energy density is given in \si{\watt\hour/\kilogram} and battery cell cost is given wrt 2023 USD/\unit{\kWh}).}
\label{tab:technoChange}
\begin{tabular}{|r|r|r||r|r|r||r|r|r|}
\hline
      Year & Energy Density &       Cost &       Year & Energy Density &       Cost &       Year & Energy Density &       Cost \\
\hline
      1991 &      98.20 &    8805.20 &       2002 &     188.18 &    1327.70 &       2013 &     243.23 &     472.73 \\
\hline
      1992 &     102.75 &    7166.22 &       2003 &     193.37 &    1014.27 &       2014 &     251.54 &     409.52 \\
\hline
      1993 &     107.86 &    5756.38 &       2004 &     209.99 &     909.97 &       2015 &     248.43 &     280.86 \\
\hline
      1994 &     112.98 &    6451.60 &       2005 &     202.72 &     788.02 &       2016 &     259.85 &     235.02 \\
\hline
      1995 &     118.66 &    6068.66 &       2006 &     207.92 &     675.05 &       2017 &     276.47 &     170.31 \\
\hline
      1996 &     126.62 &    5151.87 &       2007 &     222.46 &     644.22 &       2018 &     293.09 &     142.60 \\
\hline
      1997 &     130.59 &    4398.27 &       2008 &     218.30 &     673.49 &       2019 &     313.87 &     119.58 \\
\hline
      1998 &     144.80 &    3541.91 &       2009 &     222.46 &     580.73 &       2020 &     334.64 &     111.57 \\
\hline
      1999 &     156.37 &    2741.66 &       2010 &     230.77 &     533.49 &       2021 &     389.70 &     107.44 \\
\hline
      2000 &     160.13 &    2576.53 &       2011 &     238.04 &     499.68 &       2022 &     444.75 &     119.23 \\
\hline
      2001 &     179.87 &    1832.35 &       2012 &     243.23 &     507.12 &       2023 &     499.80 &     100.63 \\
\hline
\end{tabular}  
\end{table}

\begin{table}[h!]
\centering
\caption{Clustering results for BEBs.}
\label{tab:cluster_summary}
\begin{tabular}{|c|c|c|c|c|c|}
\hline
{Clusters} & {Cluster} & {Energy Density} & {Cost} & {Energy Density} & {Cost} \\
& {Probability} & {Improvement Rate} & {Improvement Rate} & {Change} & {Change} \\
\hline
\multirow{2}{*}{2	Clusters}	
&	\textbf{F}	(0.46)	&	0.27	&	1.05	&	1.31	&	0.35	\\
		&	\textbf{S}	(0.54)	&	0.19	&	0.45	&	1.21	&	0.63	\\
\hline													
\multirow{3}{*}{3	Clusters}	
&	\textbf{F}	(0.32	)&	0.25	&	1.14	&	1.28	&	0.32	\\
		&	\textbf{M}	(0.25)&	0.24	&	0.79	&	1.27	&	0.45	\\
		&	\textbf{S}	(0.429)	&	0.21	&	0.39	&	1.24	&	0.68	\\
\hline
\end{tabular}
\end{table}

\subsubsection{HFCB  Technological Change}

\begin{table}[h!]
\centering
\caption{Technological change for HFCBs over the years (fuel cell system cost is given wrt 2016 USD/\unit{\kW}).}
\label{tab:technoChange-fc}
\begin{tabular}{|r|r|r||r|r|r||r|r|r|}
\hline
     Year  & Efficiency &       Cost &      Year  & Efficiency &       Cost &      Year  & Efficiency &       Cost \\
\hline
      2010 &            &        175 &       2015 &     {\it } &         99 &       2020 &       0.57 &         76 \\
\hline
      2011 &            &        142 &       2016 &     {\it } &         99 &       2021 &       0.57 &     {\it } \\
\hline
      2012 &        0.50 &        122 &       2017 &       0.55 &         78 &       2022 &     {\it } &         71 \\
\hline
      2013 &        0.50 &        119 &       2018 &     {\it } &         79 &       2023 &       0.57 &            \\
\hline
      2014 &     {\it } &        107 &       2019 &     {\it } &         79 &       2024 &       0.60 &            \\
\hline
\end{tabular}  
\end{table}

\begin{table}[h!]
\centering
\caption{Clustering results for HFCBs.}
\label{tab:cluster_summary-fc}
\begin{tabular}{|c|c|c|c|c|c|}
\hline
{Clusters} & {Cluster} & {Efficiency} & {Cost} & {Efficiency} & {Cost} \\
& {Probability} & {Improvement Rate} & {Improvement Rate} & {Change} & {Change} \\
\hline
\multirow{2}{*}{2	Clusters}	&	\textbf{F}	(0.83) &	0.09	&	0.35	&	1.09	&	0.71	\\
		&	\textbf{S}	(0.17)	&	0.04	&	0.09	&	1.04	&	0.91	\\
\hline													
\multirow{3}{*}{3	Clusters}	&	\textbf{F}	(0.33)&	0.10	&	0.43	&	1.11	&	0.65	\\
		&	\textbf{M}	(0.50)	&	0.08	&	0.29	&	1.08	&	0.75	\\
		&	\textbf{S}	(0.17)	&	0.04	&	0.09	&	1.04	&	0.91	\\
\hline
\end{tabular}
\end{table}

\newpage
\subsection{Scenario Table}
\label{sec:scenarioTable}
Table~\ref{tab:scenarioTable} shows the technological improvements across different stages of each scenario, along with the scenario's probability (in percent). The columns labeled S2, S3, S4, and S5 refer to stages 2, 3, 4, and 5, respectively. In each stage, the first letter represents the BEB improvement, and the second letter refers to the HFCB improvement. For example, ``FS'' means that BEBs improve fast, while HFCBs improve slow.

\begin{table}[H]
    \centering
    \setlength{\tabcolsep}{1pt}
    \scriptsize
\begin{tabular}{|c|c|c|c|c|c||c|c|c|c|c|c||c|c|c|c|c|c||c|c|c|c|c|c|}
    \hline
        ID & S2 & S3 & S4 & S5 & Prob. (\%) & ID & S2 & S3 & S4 & S5 & Prob. (\%) & ID & S2 & S3 & S4 & S5 & Prob. (\%) & ID & S2 & S3 & S4 & S5 & Prob. (\%) \\ \hline
        1 & FF & FF & FF & FF & 2.12493 & 65 & FS & FF & FF & FF & 0.43523 & 129 & SF & FF & FF & FF & 2.49448 & 193 & SS & FF & FF & FF & 0.51092 \\ \hline
        2 & FF & FF & FF & FS & 0.43523 & 66 & FS & FF & FF & FS & 0.08914 & 130 & SF & FF & FF & FS & 0.51092 & 194 & SS & FF & FF & FS & 0.10465 \\ \hline
        3 & FF & FF & FF & SF & 2.49448 & 67 & FS & FF & FF & SF & 0.51092 & 131 & SF & FF & FF & SF & 2.92830 & 195 & SS & FF & FF & SF & 0.59977 \\ \hline
        4 & FF & FF & FF & SS & 0.51092 & 68 & FS & FF & FF & SS & 0.10465 & 132 & SF & FF & FF & SS & 0.59977 & 196 & SS & FF & FF & SS & 0.12284 \\ \hline
        5 & FF & FF & FS & FF & 0.43523 & 69 & FS & FF & FS & FF & 0.08914 & 133 & SF & FF & FS & FF & 0.51092 & 197 & SS & FF & FS & FF & 0.10465 \\ \hline
        6 & FF & FF & FS & FS & 0.08914 & 70 & FS & FF & FS & FS & 0.01826 & 134 & SF & FF & FS & FS & 0.10465 & 198 & SS & FF & FS & FS & 0.02143 \\ \hline
        7 & FF & FF & FS & SF & 0.51092 & 71 & FS & FF & FS & SF & 0.10465 & 135 & SF & FF & FS & SF & 0.59977 & 199 & SS & FF & FS & SF & 0.12284 \\
\hline
    \end{tabular}
    \caption{Scenario Table}
    \label{tab:scenarioTable}
\end{table}
    
\newpage

\section{Detailed Results of the Sensitivity Analysis}

\subsection{Relaxed Budget}

\begin{figure}[H]
\centering
%
\caption{Relaxed Budget}
\label{fig:Budget300-case}
\end{figure}

Figure~\ref{fig:Budget300-case} represents the case where the annual budget is set at 300 million USD.
The optimization process took 5,888 seconds, resulting in an objective function value of 13,586,698,671 USD. After rounding the variables, the objective function value increased to 13,615,884,770 USD.

\begin{figure}[H]
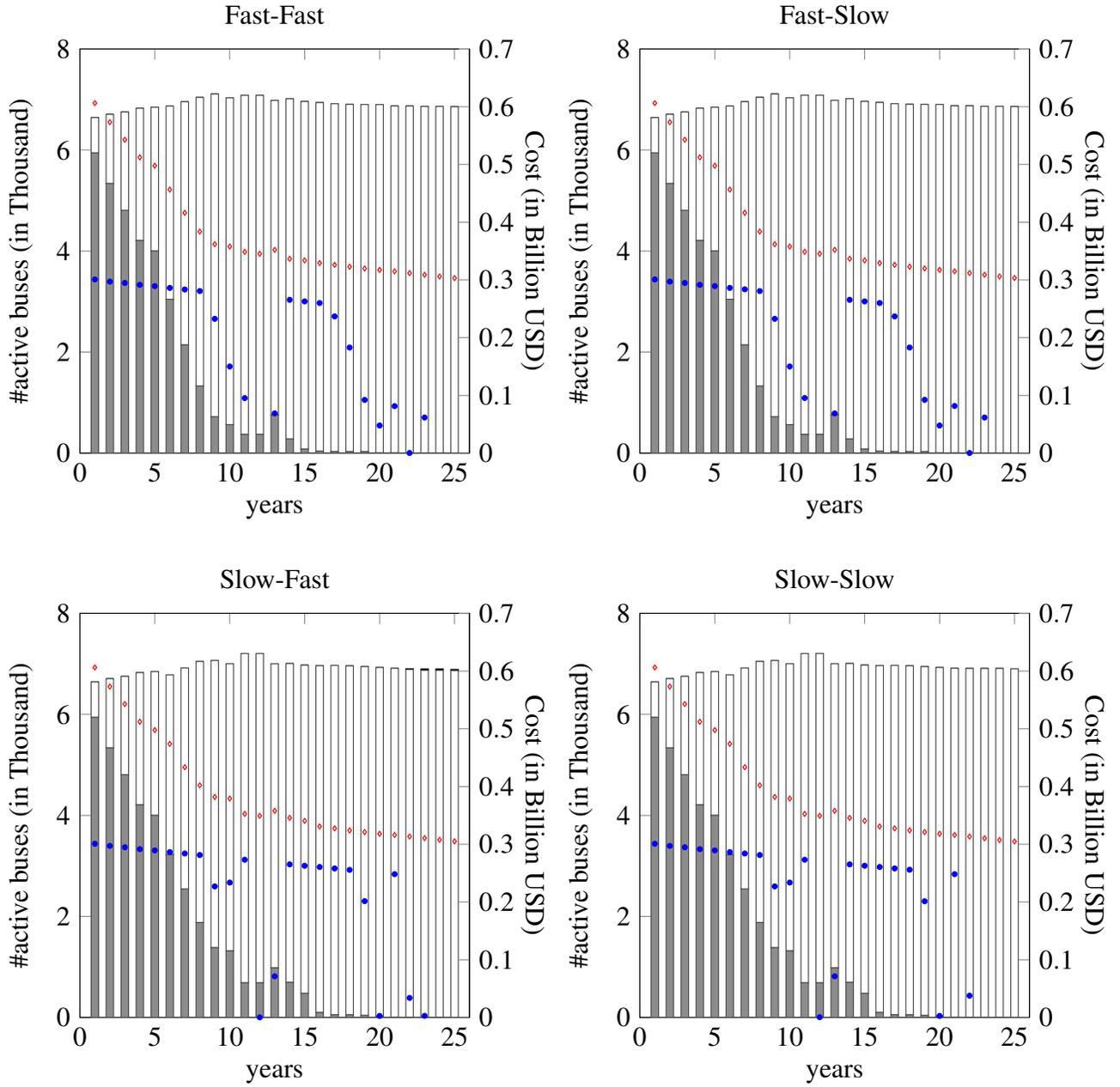

\centering

\begin{subfigure}{0.475\textwidth}
\centering
%
\label{fig:Budget300-Slow-Slow-scenario}
\end{subfigure}

\caption{Relaxed Budget -four specific scenarios.}

\label{fig:Budget300-4scenarios-adj}

\end{figure}

\subsection{Strict Emission Constraints}

\begin{figure}[H]
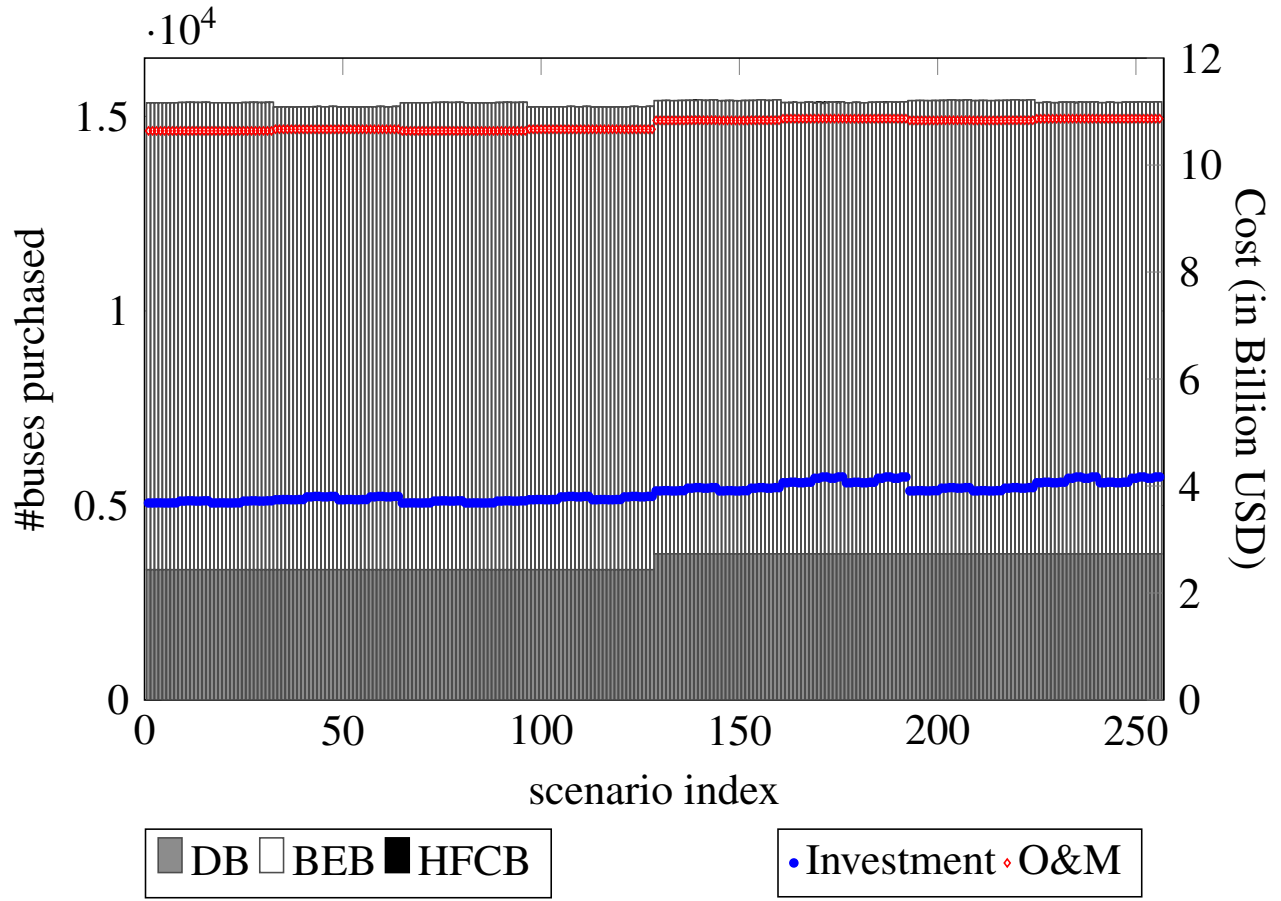

\centering
%
\caption{Strict Emission Constraints}
\label{fig:EmissionStrict-case}
\end{figure}

Figure~\ref{fig:EmissionStrict-case} represents the case with only zero emission limits in the final period, and without any intermediate targets.
The optimization process took 4,877 seconds, resulting in an objective function value of 14,133,211,491 USD. After rounding the variables, the objective function value increased to 14,165,320,512 USD.

\begin{figure}[H]
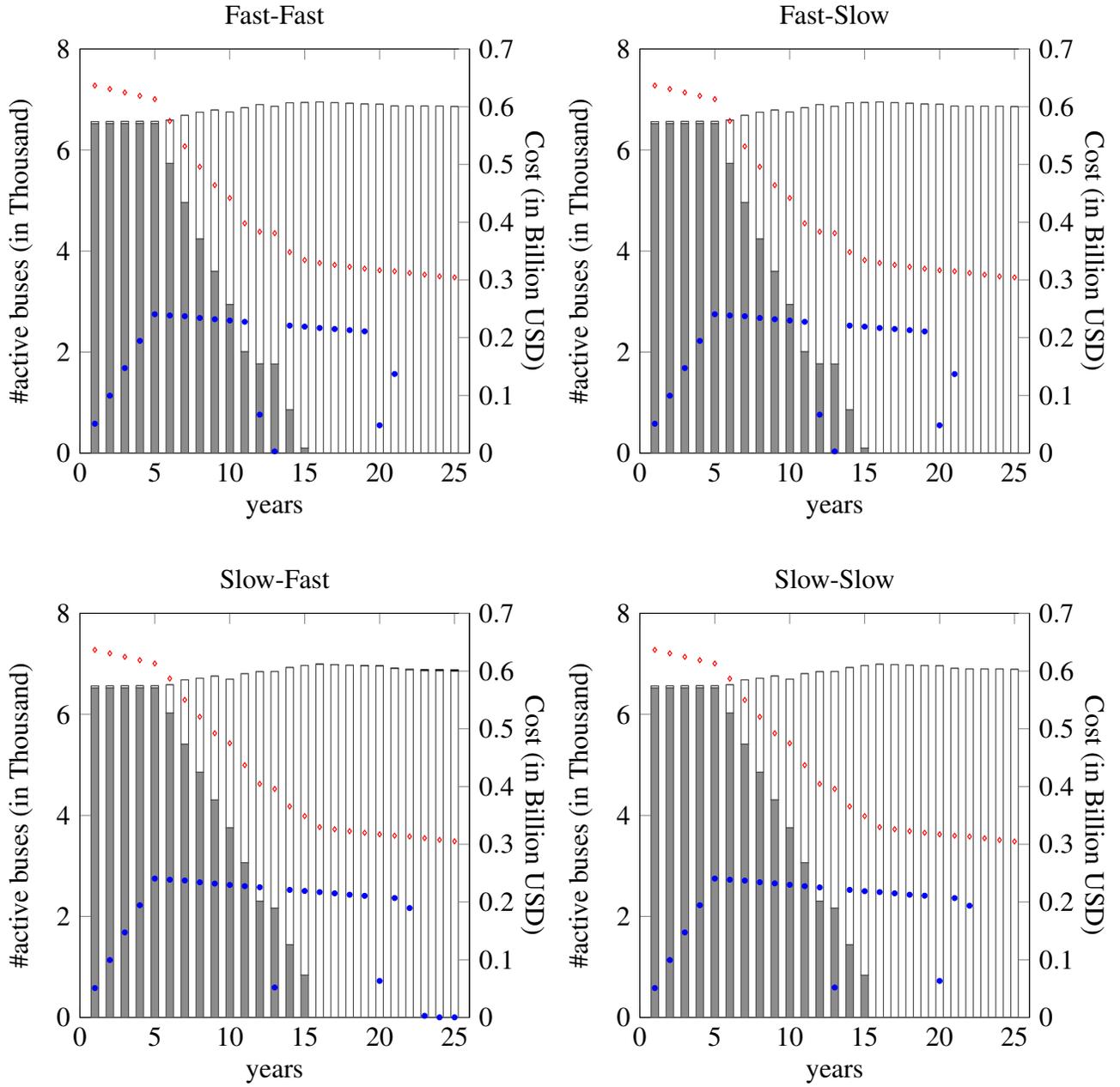

\centering

\begin{subfigure}{0.475\textwidth}
\centering
%
\label{fig:EmissionStrict-Slow-Slow-scenario}
\end{subfigure}

\caption{Strict Emission Constraints}

\label{fig:EmissionStrict-4scenarios-adj}

\end{figure}

\subsection{Reduced Hydrogen Price}

\begin{figure}[H]
\centering
%
\caption{Reduced Hydrogen  Price}
\label{fig:Hydrogen2-case}
\end{figure}

Figure~\ref{fig:Hydrogen2-case} represents the case where hydrogen price is reduced to 2 USD per \unit{\kg}.
The optimization process took 6,133 seconds, resulting in an objective function value of 14,099,073,373 USD. After rounding the variables, the objective function value increased to 14,129,697,437 USD.

\begin{figure}[H]
\centering

\begin{subfigure}{0.475\textwidth}
\centering
%
\label{fig:Hydrogen2-Slow-Slow-scenario}
\end{subfigure}

\caption{Reduced Hydrogen Price}

\label{fig:Hydrogen2-4scenarios-adj}

\end{figure}

\section{Detailed Results of Model Comparison}

\subsection{Base Case - Simplified Assignments}
\label{sec:base-simplified}

\begin{figure}[H]
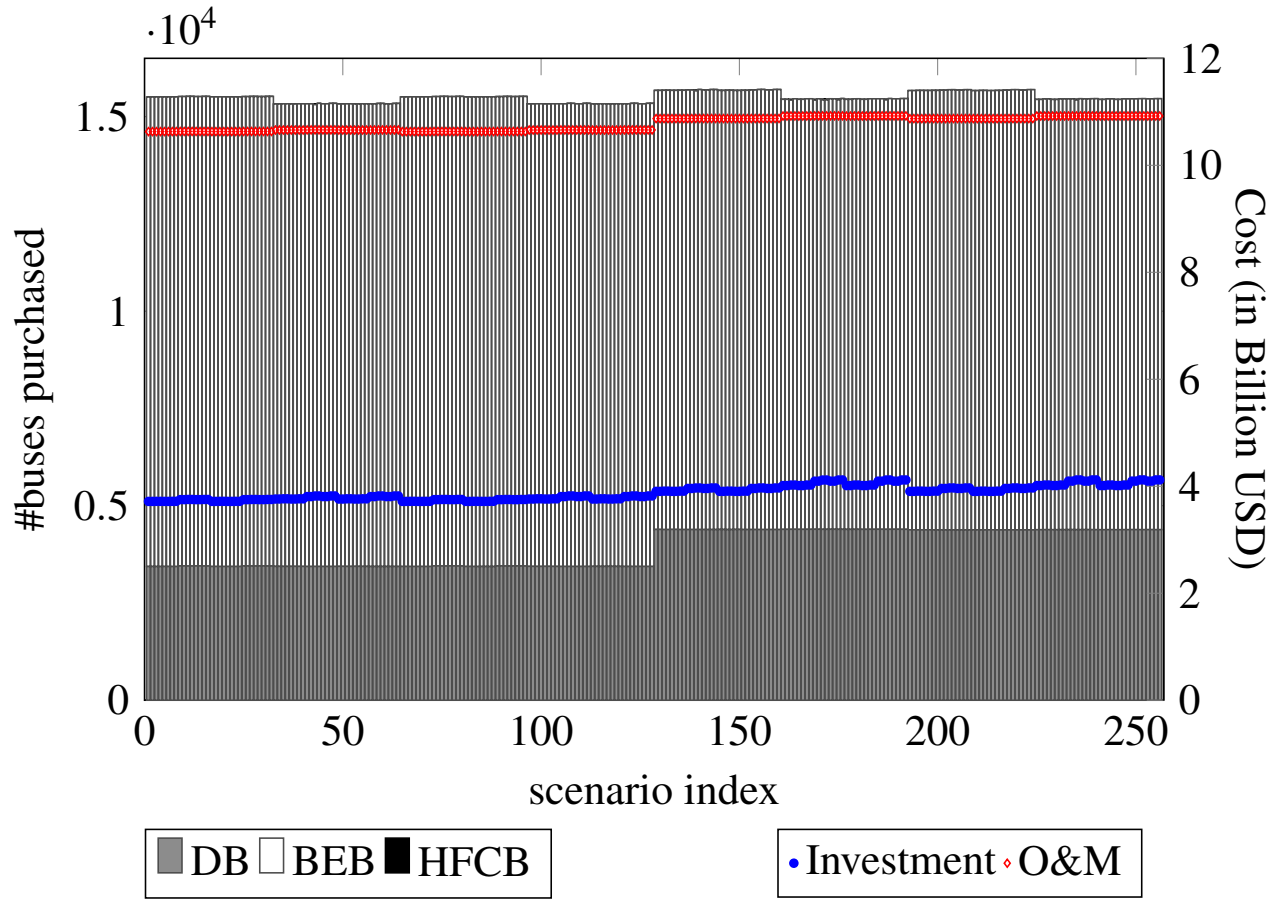

\centering
%
\caption{Base Case with Simplified Assignments}
\label{fig:BaseSimplified-case}
\end{figure}

Figure~\ref{fig:BaseSimplified-case} represents the simplified assignments version of the Base Case. In this case, we assume that for route-bus length pairs with a higher demand in the winter, the buses assigned to a cluster will remain in that cluster on all other seasons. However, we adjust the O\&M costs based on the seasonal difference in demand.
The optimization process took 1,770 seconds, resulting in an objective function value of 14,111,317,685 USD. After rounding the variables, the objective function value increased to 14,136,930,170 USD.

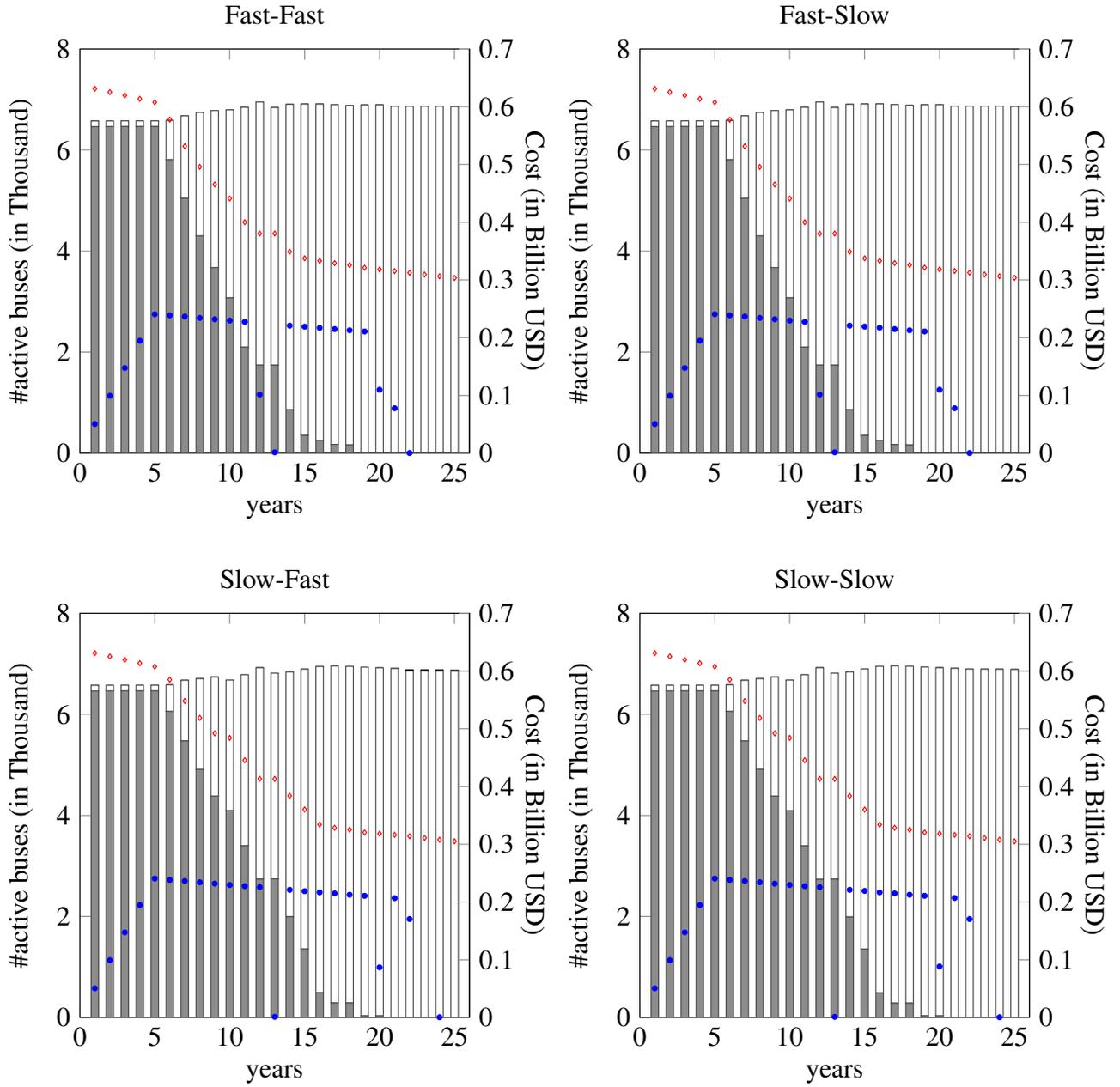
\begin{figure}[H]
\centering

\begin{subfigure}{0.475\textwidth}
\centering
\begin{tikzpicture}[scale=1.1]
      \pgfplotsset{every axis y label/.append style={rotate=0,yshift=-0.4cm}}
\begin{axis}[
    title={Fast-Fast},
    x=6pt,
    xmin=0,
    xmax=26,
    ymin = 0,
    ymax = 8,     ytick pos = left,
    bar width=3pt, 
    ybar stacked,
    ylabel={\#active buses (in Thousand)},    
    xlabel={years},
]

 \addplot+ [ybar, fill=gray!90,draw=black!70] coordinates {
(1,6.462)
(2,6.463)
(3,6.464)
(4,6.464)
(5,6.464)
(6,5.809)
(7,5.045)
(8,4.3)
(9,3.673)
(10,3.073)
(11,2.101)
(12,1.745)
(13,1.743)
(14,0.858)
(15,0.351)
(16,0.254)
(17,0.172)
(18,0.163)
(19,0)
(20,0)
(21,0)
(22,0.001)
(23,0.001)
(24,0.001)
(25,0)
};

 \addplot+ [ybar, fill=gray!0,draw=black!70] coordinates {
(1,0.112)
(2,0.112)
(3,0.112)
(4,0.112)
(5,0.112)
(6,0.776)
(7,1.63)
(8,2.442)
(9,3.108)
(10,3.721)
(11,4.743)
(12,5.204)
(13,5.097)
(14,6.046)
(15,6.558)
(16,6.655)
(17,6.724)
(18,6.719)
(19,6.891)
(20,6.893)
(21,6.864)
(22,6.863)
(23,6.863)
(24,6.863)
(25,6.858)
};

 \addplot+ [ybar, black] coordinates {
(1,0)
(2,0)
(3,0)
(4,0)
(5,0)
(6,0)
(7,0)
(8,0)
(9,0)
(10,0)
(11,0)
(12,0)
(13,0)
(14,0)
(15,0)
(16,0)
(17,0)
(18,0)
(19,0)
(20,0)
(21,0)
(22,0)
(23,0)
(24,0)
(25,0)
};

\end{axis}
      \pgfplotsset{every axis y label/.append style={rotate=180,yshift=7.2cm}}
    \begin{axis}[
    x=6pt,
    xmin=0,
    xmax=26,
    ymin=0,
    ymax=0.7,     ytick={0, 0.1, 0.2, 0.3, 0.4, 0.5, 0.6,0.7},
    axis y line*=right,
    axis x line=none,
    ylabel=Cost (in Billion USD),
    ]

 \addplot[mark=*,blue, mark size=1pt, only marks] coordinates {
(1,0.05048266)
(2,0.09924571)
(3,0.14745078)
(4,0.19472864)
(5,0.24061138)
(6,0.23873053)
(7,0.23665899)
(8,0.23409796)
(9,0.23189643)
(10,0.22968143)
(11,0.22727773)
(12,0.10149470)
(13,0.00143449)
(14,0.22084643)
(15,0.21900437)
(16,0.21693159)
(17,0.21491124)
(18,0.21281451)
(19,0.21083503)
(20,0.10976564)
(21,0.07772931)
(22,0.00024538)
};

 \addplot[mark=diamond,red, mark size=1pt, only marks] coordinates {
(1,0.63108967)
(2,0.62523277)
(3,0.61952414)
(4,0.61357882)
(5,0.60769478)
(6,0.57771534)
(7,0.53153809)
(8,0.49585085)
(9,0.46533538)
(10,0.44099925)
(11,0.40005792)
(12,0.38037550)
(13,0.38042311)
(14,0.34882175)
(15,0.33741242)
(16,0.33298548)
(17,0.32902287)
(18,0.32571347)
(19,0.32123345)
(20,0.31815234)
(21,0.31532749)
(22,0.31234899)
(23,0.30931977)
(24,0.30638207)
(25,0.30345342)
};

    \end{axis}
\end{tikzpicture}
\label{fig:BaseSimplified-Fast-Fast-scenario}
\end{subfigure}
\hfill
\begin{subfigure}{0.475\textwidth}
\centering
\begin{tikzpicture}[scale=1.1]
      \pgfplotsset{every axis y label/.append style={rotate=0,yshift=-.4cm}}
\begin{axis}[
    title={Fast-Slow},
    x=6pt,
    xmin=0,
    xmax=26,
    ymin=0,
    ymax = 8,     ytick pos = left,
    bar width=3pt, 
    ybar stacked,
    ylabel={\#active buses (in Thousand)},    
    xlabel={years},
]
 \addplot+ [ybar, fill=gray!90,draw=black!70] coordinates {
(1,6.462)
(2,6.463)
(3,6.464)
(4,6.464)
(5,6.464)
(6,5.809)
(7,5.045)
(8,4.3)
(9,3.673)
(10,3.073)
(11,2.101)
(12,1.745)
(13,1.743)
(14,0.858)
(15,0.351)
(16,0.254)
(17,0.172)
(18,0.163)
(19,0)
(20,0)
(21,0)
(22,0.001)
(23,0.001)
(24,0.001)
(25,0)
};

 \addplot+ [ybar, fill=gray!0,draw=black!70] coordinates {
(1,0.112)
(2,0.112)
(3,0.112)
(4,0.112)
(5,0.112)
(6,0.776)
(7,1.63)
(8,2.442)
(9,3.108)
(10,3.721)
(11,4.743)
(12,5.204)
(13,5.097)
(14,6.046)
(15,6.559)
(16,6.658)
(17,6.727)
(18,6.722)
(19,6.894)
(20,6.896)
(21,6.867)
(22,6.866)
(23,6.866)
(24,6.866)
(25,6.861)
};

 \addplot+ [ybar, black] coordinates {
(1,0)
(2,0)
(3,0)
(4,0)
(5,0)
(6,0)
(7,0)
(8,0)
(9,0)
(10,0)
(11,0)
(12,0)
(13,0)
(14,0)
(15,0)
(16,0)
(17,0)
(18,0)
(19,0)
(20,0)
(21,0)
(22,0)
(23,0)
(24,0)
(25,0)
};

\end{axis}
      \pgfplotsset{every axis y label/.append style={rotate=180,yshift=7.2cm}}
    \begin{axis}[
    x=6pt,
    xmin=0,
    xmax=26,
    ymin=0,
    ymax=0.7,     ytick={0, 0.1, 0.2, 0.3, 0.4, 0.5, 0.6,0.7},
    axis y line*=right,
    axis x line=none,
    ylabel=Cost (in Billion USD),
    ]

 \addplot[mark=*,blue, mark size=1pt, only marks] coordinates {
(1,0.05048266)
(2,0.09924571)
(3,0.14745078)
(4,0.19472864)
(5,0.24061138)
(6,0.23873053)
(7,0.23665899)
(8,0.23409796)
(9,0.23189643)
(10,0.22968143)
(11,0.22727773)
(12,0.10149470)
(13,0.00143449)
(14,0.22084643)
(15,0.21920180)
(16,0.21728104)
(17,0.21491124)
(18,0.21281451)
(19,0.21083503)
(20,0.10976564)
(21,0.07772931)
(22,0.00024538)
};

 \addplot[mark=diamond,red, mark size=1pt, only marks] coordinates {
(1,0.63108967)
(2,0.62523277)
(3,0.61952414)
(4,0.61357882)
(5,0.60769478)
(6,0.57749158)
(7,0.53146653)
(8,0.49571611)
(9,0.46533888)
(10,0.44094912)
(11,0.39998881)
(12,0.38015902)
(13,0.38039059)
(14,0.34886652)
(15,0.33744322)
(16,0.33309899)
(17,0.32922374)
(18,0.32582885)
(19,0.32138083)
(20,0.31830649)
(21,0.31554275)
(22,0.31248819)
(23,0.30948642)
(24,0.30653893)
(25,0.30367049)
};

    \end{axis}
\end{tikzpicture}
\label{fig:BaseSimplified-Fast-Slow-scenario}
\end{subfigure}

\begin{subfigure}{0.475\textwidth}
\centering
\begin{tikzpicture}[scale=1.1]
      \pgfplotsset{every axis y label/.append style={rotate=0,yshift=-.4cm}}
\begin{axis}[
    title={Slow-Fast},
    x=6pt,
    xmin=0,
    xmax=26,
    ymin = 0,
    ymax= 8,
    ytick pos = left,
    bar width=3pt,
    ybar stacked,
    ylabel={\#active buses (in Thousand)},    
    xlabel={years},
]

 \addplot+ [ybar, fill=gray!90,draw=black!70] coordinates {
(1,6.462)
(2,6.463)
(3,6.464)
(4,6.464)
(5,6.464)
(6,6.059)
(7,5.474)
(8,4.911)
(9,4.383)
(10,4.097)
(11,3.402)
(12,2.739)
(13,2.743)
(14,1.994)
(15,1.359)
(16,0.49)
(17,0.289)
(18,0.289)
(19,0.033)
(20,0.033)
(21,0.003)
(22,0.003)
(23,0.003)
(24,0.005)
(25,0)
};

 \addplot+ [ybar, fill=gray!0,draw=black!70] coordinates {
(1,0.112)
(2,0.112)
(3,0.112)
(4,0.112)
(5,0.112)
(6,0.525)
(7,1.201)
(8,1.795)
(9,2.356)
(10,2.581)
(11,3.382)
(12,4.184)
(13,4.072)
(14,4.846)
(15,5.536)
(16,6.457)
(17,6.669)
(18,6.66)
(19,6.901)
(20,6.888)
(21,6.905)
(22,6.857)
(23,6.856)
(24,6.853)
(25,6.853)

};

 \addplot+ [ybar, black] coordinates {
(1,0)
(2,0)
(3,0)
(4,0)
(5,0)
(6,0)
(7,0)
(8,0)
(9,0)
(10,0)
(11,0)
(12,0)
(13,0)
(14,0)
(15,0)
(16,0)
(17,0)
(18,0)
(19,0)
(20,0)
(21,0)
(22,0.021)
(23,0.021)
(24,0.021)
(25,0.021)
};

\end{axis}
      \pgfplotsset{every axis y label/.append style={rotate=180,yshift=7.2cm}}
    \begin{axis}[
    x=6pt,
    xmin=0,
    xmax=26,
    ymin=0,
    ymax=0.7,     ytick={0, 0.1, 0.2, 0.3, 0.4, 0.5, 0.6,0.7},
    axis y line*=right,
    axis x line=none,
    ylabel=Cost (in Billion USD),
    ]

 \addplot[mark=*,blue, mark size=1pt, only marks] coordinates {
(1,0.05048266)
(2,0.09924571)
(3,0.14745078)
(4,0.19472864)
(5,0.24061138)
(6,0.23836232)
(7,0.23635007)
(8,0.23408439)
(9,0.23190931)
(10,0.22951133)
(11,0.22748023)
(12,0.22567303)
(13,0.00106982)
(14,0.22113980)
(15,0.21880523)
(16,0.21672042)
(17,0.21488785)
(18,0.21262321)
(19,0.21077073)
(20,0.08684569)
(21,0.20673112)
(22,0.17049527)
(24,0.00040122)
};

 \addplot[mark=diamond,red, mark size=1pt, only marks] coordinates {
(1,0.63108967)
(2,0.62523277)
(3,0.61952414)
(4,0.61357882)
(5,0.60769478)
(6,0.58512714)
(7,0.54807991)
(8,0.51894948)
(9,0.49212034)
(10,0.48423410)
(11,0.44558209)
(12,0.41349340)
(13,0.41314439)
(14,0.38397929)
(15,0.36036121)
(16,0.33411062)
(17,0.32851961)
(18,0.32531837)
(19,0.32056255)
(20,0.31841968)
(21,0.31632226)
(22,0.31400269)
(23,0.31097660)
(24,0.30790580)
(25,0.30522029)
};

    \end{axis}
\end{tikzpicture}
\label{fig:BaseSimplified-Slow-Fast-scenario}
\end{subfigure}
\hfill
\begin{subfigure}{0.475\textwidth}
\centering
\begin{tikzpicture}[scale=1.1]
      \pgfplotsset{every axis y label/.append style={rotate=0,yshift=-.4cm}}
\begin{axis}[
    title={Slow-Slow},
    x=6pt,
    xmin=0,
    xmax=26,
    ymin = 0,
    ymax=  8,
    ytick pos = left,
    bar width=3pt, 
    ybar stacked,
    ylabel={\#active buses (in Thousand)},    
    xlabel={years},
]

 \addplot+ [ybar, fill=gray!90,draw=black!70] coordinates {
(1,6.462)
(2,6.463)
(3,6.464)
(4,6.464)
(5,6.464)
(6,6.059)
(7,5.473)
(8,4.912)
(9,4.383)
(10,4.092)
(11,3.397)
(12,2.734)
(13,2.738)
(14,1.989)
(15,1.354)
(16,0.485)
(17,0.284)
(18,0.284)
(19,0.033)
(20,0.033)
(21,0.003)
(22,0.003)
(23,0.003)
(24,0.005)
(25,0)
};

 \addplot+ [ybar, fill=gray!0,draw=black!70] coordinates {
(1,0.112)
(2,0.112)
(3,0.112)
(4,0.112)
(5,0.112)
(6,0.525)
(7,1.201)
(8,1.794)
(9,2.356)
(10,2.586)
(11,3.387)
(12,4.189)
(13,4.077)
(14,4.851)
(15,5.544)
(16,6.465)
(17,6.676)
(18,6.668)
(19,6.903)
(20,6.89)
(21,6.908)
(22,6.893)
(23,6.892)
(24,6.889)
(25,6.889)
};

 \addplot+ [ybar, black] coordinates {
(1,0)
(2,0)
(3,0)
(4,0)
(5,0)
(6,0)
(7,0)
(8,0)
(9,0)
(10,0)
(11,0)
(12,0)
(13,0)
(14,0)
(15,0)
(16,0)
(17,0)
(18,0)
(19,0)
(20,0)
(21,0)
(22,0)
(23,0)
(24,0)
(25,0)
};

\end{axis}
      \pgfplotsset{every axis y label/.append style={rotate=180,yshift=7.2cm}}
    \begin{axis}[
    x=6pt,
    xmin=0,
    xmax=26,
    ymin=0,
    ymax=0.7,     ytick={0, 0.1, 0.2, 0.3, 0.4, 0.5, 0.6,0.7},
      axis y line*=right,
      axis x line=none,
      ylabel=Cost (in Billion USD),
    ]

 \addplot[mark=*,blue, mark size=1pt, only marks] coordinates {
(1,0.05048266)
(2,0.09924571)
(3,0.14745078)
(4,0.19472864)
(5,0.24061138)
(6,0.23836232)
(7,0.23635007)
(8,0.23420926)
(9,0.23178562)
(10,0.22957156)
(11,0.22748023)
(12,0.22567303)
(13,0.00106982)
(14,0.22113980)
(15,0.21927134)
(16,0.21672042)
(17,0.21484530)
(18,0.21262321)
(19,0.21079714)
(20,0.08841331)
(21,0.20702991)
(22,0.17043186)
(24,0.00040122)
};

 \addplot[mark=diamond,red, mark size=1pt, only marks] coordinates {
(1,0.63108967)
(2,0.62523277)
(3,0.61952414)
(4,0.61357882)
(5,0.60769478)
(6,0.58510067)
(7,0.54806438)
(8,0.51887732)
(9,0.49206994)
(10,0.48422051)
(11,0.44543200)
(12,0.41327049)
(13,0.41308166)
(14,0.38374725)
(15,0.36020727)
(16,0.33390864)
(17,0.32848411)
(18,0.32518791)
(19,0.32057500)
(20,0.31846494)
(21,0.31621611)
(22,0.31399722)
(23,0.31096977)
(24,0.30790683)
(25,0.30511319)
};

    \end{axis}
\end{tikzpicture}
\label{fig:BaseSimplified-Slow-Slow-scenario}
\end{subfigure}

\caption{Base Case - Simplified Assignments}

\label{fig:BaseSimplified-4scenarios-adj}

\end{figure}

\newpage

\subsection{Extended Scenario Tree - 3 By 2}

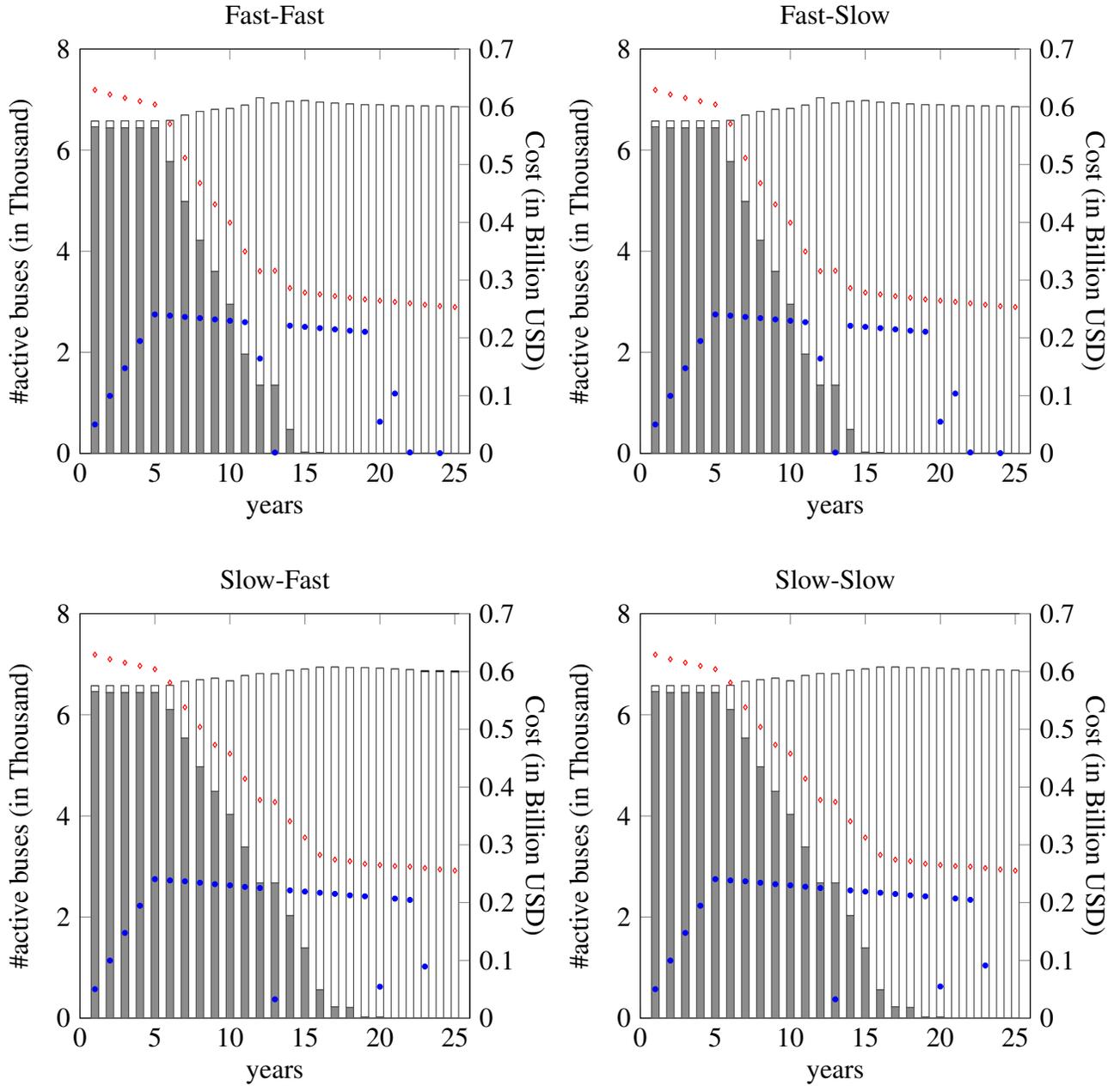
\begin{figure}[H]
\centering

\begin{subfigure}{0.475\textwidth}
\centering
\begin{tikzpicture}[scale=1.1]
      \pgfplotsset{every axis y label/.append style={rotate=0,yshift=-0.4cm}}
\begin{axis}[
    title={Fast-Fast},
    x=6pt,
    xmin=0,
    xmax=26,
    ymin = 0,
    ymax = 8,     ytick pos = left,
    bar width=3pt, 
    ybar stacked,
    ylabel={\#active buses (in Thousand)},    
    xlabel={years},
]

 \addplot+ [ybar, fill=gray!90,draw=black!70] coordinates {
(1,6.461)
(2,6.438)
(3,6.439)
(4,6.44)
(5,6.44)
(6,5.772)
(7,4.982)
(8,4.216)
(9,3.601)
(10,2.949)
(11,1.966)
(12,1.348)
(13,1.348)
(14,0.475)
(15,0.025)
(16,0.021)
(17,0)
(18,0)
(19,0)
(20,0)
(21,0)
(22,0.007)
(23,0.007)
(24,0.009)
(25,0)
};

 \addplot+ [ybar, fill=gray!0,draw=black!70] coordinates {
(1,0.111)
(2,0.138)
(3,0.138)
(4,0.138)
(5,0.138)
(6,0.812)
(7,1.709)
(8,2.549)
(9,3.203)
(10,3.872)
(11,4.922)
(12,5.687)
(13,5.581)
(14,6.49)
(15,6.954)
(16,6.927)
(17,6.929)
(18,6.912)
(19,6.897)
(20,6.896)
(21,6.873)
(22,6.865)
(23,6.865)
(24,6.862)
(25,6.857)
};

 \addplot+ [ybar, black] coordinates {
(1,0)
(2,0)
(3,0)
(4,0)
(5,0)
(6,0)
(7,0)
(8,0)
(9,0)
(10,0)
(11,0)
(12,0)
(13,0)
(14,0)
(15,0)
(16,0)
(17,0)
(18,0)
(19,0)
(20,0)
(21,0)
(22,0)
(23,0)
(24,0)
(25,0)

};

\end{axis}
      \pgfplotsset{every axis y label/.append style={rotate=180,yshift=7.2cm}}
    \begin{axis}[
    x=6pt,
    xmin=0,
    xmax=26,
    ymin=0,
    ymax=0.7,     ytick={0, 0.1, 0.2, 0.3, 0.4, 0.5, 0.6,0.7},
    axis y line*=right,
    axis x line=none,
    ylabel=Cost (in Billion USD),
    ]

 \addplot[mark=*,blue, mark size=1pt, only marks] coordinates {
(1,0.05013665)
(2,0.09964426)
(3,0.14754888)
(4,0.19463147)
(5,0.24061138)
(6,0.23855519)
(7,0.23629269)
(8,0.23426018)
(9,0.23199491)
(10,0.22968344)
(11,0.22728089)
(12,0.16419533)
(13,0.00151247)
(14,0.22111847)
(15,0.21910629)
(16,0.21688347)
(17,0.21479647)
(18,0.21250556)
(19,0.21064064)
(20,0.05484681)
(21,0.10363122)
(22,0.00171769)
(24,0.00048146)
};

 \addplot[mark=diamond,red, mark size=1pt, only marks] coordinates {
(1,0.62911821)
(2,0.62122838)
(3,0.61525206)
(4,0.60965821)
(5,0.60400238)
(6,0.57030364)
(7,0.51166972)
(8,0.46796108)
(9,0.43113054)
(10,0.39941315)
(11,0.34955435)
(12,0.31556520)
(13,0.31622482)
(14,0.28617926)
(15,0.27835896)
(16,0.27544759)
(17,0.27227746)
(18,0.26941108)
(19,0.26669484)
(20,0.26438581)
(21,0.26228992)
(22,0.25984741)
(23,0.25735334)
(24,0.25500694)
(25,0.25336264)
};

    \end{axis}
\end{tikzpicture}
\label{fig:ExtendedTree-Fast-Fast-scenario}
\end{subfigure}
\hfill
\begin{subfigure}{0.475\textwidth}
\centering
\begin{tikzpicture}[scale=1.1]
      \pgfplotsset{every axis y label/.append style={rotate=0,yshift=-.4cm}}
\begin{axis}[
    title={Fast-Slow},
    x=6pt,
    xmin=0,
    xmax=26,
    ymin=0,
    ymax = 8,     ytick pos = left,
    bar width=3pt, 
    ybar stacked,
    ylabel={\#active buses (in Thousand)},    
    xlabel={years},
]
 \addplot+ [ybar, fill=gray!90,draw=black!70] coordinates {
(1,6.461)
(2,6.438)
(3,6.439)
(4,6.44)
(5,6.44)
(6,5.772)
(7,4.982)
(8,4.216)
(9,3.601)
(10,2.949)
(11,1.966)
(12,1.348)
(13,1.348)
(14,0.475)
(15,0.025)
(16,0.021)
(17,0)
(18,0)
(19,0)
(20,0)
(21,0)
(22,0.007)
(23,0.007)
(24,0.009)
(25,0)
};

 \addplot+ [ybar, fill=gray!0,draw=black!70] coordinates {
(1,0.111)
(2,0.138)
(3,0.138)
(4,0.138)
(5,0.138)
(6,0.812)
(7,1.709)
(8,2.549)
(9,3.203)
(10,3.872)
(11,4.922)
(12,5.687)
(13,5.581)
(14,6.489)
(15,6.953)
(16,6.926)
(17,6.928)
(18,6.911)
(19,6.897)
(20,6.896)
(21,6.873)
(22,6.865)
(23,6.865)
(24,6.862)
(25,6.857)
};

 \addplot+ [ybar, black] coordinates {
(1,0)
(2,0)
(3,0)
(4,0)
(5,0)
(6,0)
(7,0)
(8,0)
(9,0)
(10,0)
(11,0)
(12,0)
(13,0)
(14,0)
(15,0)
(16,0)
(17,0)
(18,0)
(19,0)
(20,0)
(21,0)
(22,0)
(23,0)
(24,0)
(25,0)
};

\end{axis}
      \pgfplotsset{every axis y label/.append style={rotate=180,yshift=7.2cm}}
    \begin{axis}[
    x=6pt,
    xmin=0,
    xmax=26,
    ymin=0,
    ymax=0.7,     ytick={0, 0.1, 0.2, 0.3, 0.4, 0.5, 0.6,0.7},
    axis y line*=right,
    axis x line=none,
    ylabel=Cost (in Billion USD),
    ]

 \addplot[mark=*,blue, mark size=1pt, only marks] coordinates {
(1,0.05013665)
(2,0.09964426)
(3,0.14754888)
(4,0.19463147)
(5,0.24061138)
(6,0.23855519)
(7,0.23629269)
(8,0.23426018)
(9,0.23199491)
(10,0.22968344)
(11,0.22728089)
(12,0.16419533)
(13,0.00151247)
(14,0.22092489)
(15,0.21910629)
(16,0.21688347)
(17,0.21479647)
(18,0.21250556)
(19,0.21077296)
(20,0.05484681)
(21,0.10363123)
(22,0.00171769)
(24,0.00048146)
};

 \addplot[mark=diamond,red, mark size=1pt, only marks] coordinates {
(1,0.62911821)
(2,0.62122838)
(3,0.61525206)
(4,0.60965821)
(5,0.60400238)
(6,0.57036642)
(7,0.51146760)
(8,0.46788875)
(9,0.43104959)
(10,0.39937994)
(11,0.34943790)
(12,0.31549463)
(13,0.31620994)
(14,0.28621766)
(15,0.27838217)
(16,0.27533496)
(17,0.27224662)
(18,0.26944356)
(19,0.26668096)
(20,0.26443173)
(21,0.26240743)
(22,0.25984736)
(23,0.25738707)
(24,0.25492241)
(25,0.25339791)

};

    \end{axis}
\end{tikzpicture}
\label{fig:ExtendedTree-Fast-Slow-scenario}
\end{subfigure}

\begin{subfigure}{0.475\textwidth}
\centering
\begin{tikzpicture}[scale=1.1]
      \pgfplotsset{every axis y label/.append style={rotate=0,yshift=-.4cm}}
\begin{axis}[
    title={Slow-Fast},
    x=6pt,
    xmin=0,
    xmax=26,
    ymin = 0,
    ymax= 8,
    ytick pos = left,
    bar width=3pt,
    ybar stacked,
    ylabel={\#active buses (in Thousand)},    
    xlabel={years},
]

 \addplot+ [ybar, fill=gray!90,draw=black!70] coordinates {
(1,6.461)
(2,6.438)
(3,6.439)
(4,6.44)
(5,6.44)
(6,6.105)
(7,5.538)
(8,4.972)
(9,4.486)
(10,4.03)
(11,3.386)
(12,2.676)
(13,2.676)
(14,2.027)
(15,1.388)
(16,0.561)
(17,0.224)
(18,0.211)
(19,0.024)
(20,0.024)
(21,0)
(22,0)
(23,0.001)
(24,0.001)
(25,0)
};

 \addplot+ [ybar, fill=gray!0,draw=black!70] coordinates {
(1,0.111)
(2,0.138)
(3,0.138)
(4,0.138)
(5,0.138)
(6,0.477)
(7,1.126)
(8,1.72)
(9,2.236)
(10,2.64)
(11,3.389)
(12,4.139)
(13,4.136)
(14,4.852)
(15,5.516)
(16,6.379)
(17,6.719)
(18,6.724)
(19,6.909)
(20,6.897)
(21,6.906)
(22,6.892)
(23,6.85)
(24,6.85)
(25,6.844)
};

 \addplot+ [ybar, black] coordinates {
(1,0)
(2,0)
(3,0)
(4,0)
(5,0)
(6,0)
(7,0)
(8,0)
(9,0)
(10,0)
(11,0)
(12,0)
(13,0)
(14,0)
(15,0)
(16,0)
(17,0)
(18,0)
(19,0)
(20,0)
(21,0)
(22,0)
(23,0.021)
(24,0.021)
(25,0.021)
};

\end{axis}
      \pgfplotsset{every axis y label/.append style={rotate=180,yshift=7.2cm}}
    \begin{axis}[
    x=6pt,
    xmin=0,
    xmax=26,
    ymin=0,
    ymax=0.7,     ytick={0, 0.1, 0.2, 0.3, 0.4, 0.5, 0.6,0.7},
    axis y line*=right,
    axis x line=none,
    ylabel=Cost (in Billion USD),
    ]

 \addplot[mark=*,blue, mark size=1pt, only marks] coordinates {
(1,0.05013665)
(2,0.09964426)
(3,0.14754888)
(4,0.19463147)
(5,0.24061138)
(6,0.23840301)
(7,0.23687918)
(8,0.23444062)
(9,0.23194368)
(10,0.22989906)
(11,0.22720068)
(12,0.22526105)
(13,0.03257936)
(14,0.22116026)
(15,0.21909470)
(16,0.21708354)
(17,0.21521989)
(18,0.21256318)
(19,0.21094005)
(20,0.05441506)
(21,0.20689595)
(22,0.20463024)
(23,0.08935349)
};

 \addplot[mark=diamond,red, mark size=1pt, only marks] coordinates {
(1,0.62911821)
(2,0.62122838)
(3,0.61525206)
(4,0.60965821)
(5,0.60400238)
(6,0.58072850)
(7,0.53797044)
(8,0.50406129)
(9,0.47303087)
(10,0.45778238)
(11,0.41427878)
(12,0.37763589)
(13,0.37404760)
(14,0.34085113)
(15,0.31258488)
(16,0.28256743)
(17,0.27440660)
(18,0.27157192)
(19,0.26707110)
(20,0.26503389)
(21,0.26317988)
(22,0.26219611)
(23,0.25985399)
(24,0.25744213)
(25,0.25541627)
};

    \end{axis}
\end{tikzpicture}
\label{fig:ExtendedTree-Slow-Fast-scenario}
\end{subfigure}
\hfill
\begin{subfigure}{0.475\textwidth}
\centering
\begin{tikzpicture}[scale=1.1]
      \pgfplotsset{every axis y label/.append style={rotate=0,yshift=-.4cm}}
\begin{axis}[
    title={Slow-Slow},
    x=6pt,
    xmin=0,
    xmax=26,
    ymin = 0,
    ymax=  8,
    ytick pos = left,
    bar width=3pt, 
    ybar stacked,
    ylabel={\#active buses (in Thousand)},    
    xlabel={years},
]

 \addplot+ [ybar, fill=gray!90,draw=black!70] coordinates {
(1,6.461)
(2,6.438)
(3,6.439)
(4,6.44)
(5,6.44)
(6,6.105)
(7,5.538)
(8,4.972)
(9,4.486)
(10,4.03)
(11,3.386)
(12,2.676)
(13,2.676)
(14,2.027)
(15,1.388)
(16,0.561)
(17,0.224)
(18,0.211)
(19,0.024)
(20,0.024)
(21,0)
(22,0)
(23,0.001)
(24,0.001)
(25,0)
};

 \addplot+ [ybar, fill=gray!0,draw=black!70] coordinates {
(1,0.111)
(2,0.138)
(3,0.138)
(4,0.138)
(5,0.138)
(6,0.477)
(7,1.126)
(8,1.721)
(9,2.237)
(10,2.641)
(11,3.391)
(12,4.141)
(13,4.138)
(14,4.854)
(15,5.518)
(16,6.381)
(17,6.719)
(18,6.724)
(19,6.908)
(20,6.897)
(21,6.907)
(22,6.894)
(23,6.886)
(24,6.886)
(25,6.88)
};

 \addplot+ [ybar, black] coordinates {
(1,0)
(2,0)
(3,0)
(4,0)
(5,0)
(6,0)
(7,0)
(8,0)
(9,0)
(10,0)
(11,0)
(12,0)
(13,0)
(14,0)
(15,0)
(16,0)
(17,0)
(18,0)
(19,0)
(20,0)
(21,0)
(22,0)
(23,0)
(24,0)
(25,0)
};

\end{axis}
      \pgfplotsset{every axis y label/.append style={rotate=180,yshift=7.2cm}}
    \begin{axis}[
    x=6pt,
    xmin=0,
    xmax=26,
    ymin=0,
    ymax=0.7,     ytick={0, 0.1, 0.2, 0.3, 0.4, 0.5, 0.6,0.7},
      axis y line*=right,
      axis x line=none,
      ylabel=Cost (in Billion USD),
    ]

 \addplot[mark=*,blue, mark size=1pt, only marks] coordinates {
(1,0.05013665)
(2,0.09964426)
(3,0.14754888)
(4,0.19463147)
(5,0.24061138)
(6,0.23840301)
(7,0.23687918)
(8,0.23444062)
(9,0.23194368)
(10,0.22989906)
(11,0.22748195)
(12,0.22526105)
(13,0.03257936)
(14,0.22116026)
(15,0.21909470)
(16,0.21708354)
(17,0.21496725)
(18,0.21256318)
(19,0.21094005)
(20,0.05474485)
(21,0.20711706)
(22,0.20484924)
(23,0.09116438)
};

 \addplot[mark=diamond,red, mark size=1pt, only marks] coordinates {
(1,0.62911821)
(2,0.62122838)
(3,0.61525206)
(4,0.60965821)
(5,0.60400238)
(6,0.58070628)
(7,0.53801069)
(8,0.50402959)
(9,0.47307563)
(10,0.45782195)
(11,0.41443841)
(12,0.37768503)
(13,0.37426302)
(14,0.34075971)
(15,0.31257009)
(16,0.28260095)
(17,0.27433539)
(18,0.27161279)
(19,0.26724300)
(20,0.26505656)
(21,0.26321077)
(22,0.26225816)
(23,0.25975665)
(24,0.25727174)
(25,0.25532443)
};

    \end{axis}
\end{tikzpicture}
\label{fig:ExtendedTree-Slow-Slow-scenario}
\end{subfigure}

\caption{Extended Scenario Tree - 3 By 2}

\label{fig:ExtendedTree-4scenarios-adj}

\end{figure}

We present the results of the Extended Scenario Tree - 3 By 2 case in Figure~\ref{fig:ExtendedTree-4scenarios-adj}. In this case, three branches for BEB technological improvements and 2 branches for HFCB technological improvements are considered in each stage. Similar to \ref{sec:base-simplified}, we simplify the assignment decisions for route-bus length combinations having higher demand in the winter.
CPU time is 17,673 seconds, the total expected costs is 13,265,506,356 USD, which increases to 13,288,586,164 USD after rounding the variables.

\subsection{Extended Scenario Tree - 2 By 3}

\begin{figure}[H]
\centering

\begin{subfigure}{0.475\textwidth}
\centering
\begin{tikzpicture}[scale=1.1]
      \pgfplotsset{every axis y label/.append style={rotate=0,yshift=-0.4cm}}
\begin{axis}[
    title={Fast-Fast},
    x=6pt,
    xmin=0,
    xmax=26,
    ymin = 0,
    ymax = 8,     ytick pos = left,
    bar width=3pt, 
    ybar stacked,
    ylabel={\#active buses (in Thousand)},    
    xlabel={years},
]

 \addplot+ [ybar, fill=gray!90,draw=black!70] coordinates {
(1,6.462)
(2,6.463)
(3,6.464)
(4,6.464)
(5,6.464)
(6,5.809)
(7,5.045)
(8,4.3)
(9,3.673)
(10,3.073)
(11,2.101)
(12,1.745)
(13,1.743)
(14,0.858)
(15,0.351)
(16,0.254)
(17,0.172)
(18,0.163)
(19,0)
(20,0)
(21,0)
(22,0.001)
(23,0.001)
(24,0.001)
(25,0)
};

 \addplot+ [ybar, fill=gray!0,draw=black!70] coordinates {
(1,0.112)
(2,0.112)
(3,0.112)
(4,0.112)
(5,0.112)
(6,0.776)
(7,1.63)
(8,2.442)
(9,3.108)
(10,3.721)
(11,4.743)
(12,5.204)
(13,5.097)
(14,6.047)
(15,6.56)
(16,6.657)
(17,6.726)
(18,6.721)
(19,6.894)
(20,6.896)
(21,6.866)
(22,6.865)
(23,6.865)
(24,6.865)
(25,6.86)
};

 \addplot+ [ybar, black] coordinates {
(1,0)
(2,0)
(3,0)
(4,0)
(5,0)
(6,0)
(7,0)
(8,0)
(9,0)
(10,0)
(11,0)
(12,0)
(13,0)
(14,0)
(15,0)
(16,0)
(17,0)
(18,0)
(19,0)
(20,0)
(21,0)
(22,0)
(23,0)
(24,0)
(25,0)

};

\end{axis}
      \pgfplotsset{every axis y label/.append style={rotate=180,yshift=7.2cm}}
    \begin{axis}[
    x=6pt,
    xmin=0,
    xmax=26,
    ymin=0,
    ymax=0.7,     ytick={0, 0.1, 0.2, 0.3, 0.4, 0.5, 0.6,0.7},
    axis y line*=right,
    axis x line=none,
    ylabel=Cost (in Billion USD),
    ]

 \addplot[mark=*,blue, mark size=1pt, only marks] coordinates {
(1,0.05048266)
(2,0.09924571)
(3,0.14745078)
(4,0.19472864)
(5,0.24061138)
(6,0.23873053)
(7,0.23668358)
(8,0.23409796)
(9,0.23189643)
(10,0.22968143)
(11,0.22727773)
(12,0.10149470)
(13,0.00143449)
(14,0.22104585)
(15,0.21900437)
(16,0.21693362)
(17,0.21491124)
(18,0.21281451)
(19,0.21112760)
(20,0.10976564)
(21,0.07773957)
(22,0.00024538)
};

 \addplot[mark=diamond,red, mark size=1pt, only marks] coordinates {
(1,0.63108967)
(2,0.62532034)
(3,0.61943815)
(4,0.61348375)
(5,0.60786673)
(6,0.57736147)
(7,0.53146254)
(8,0.49589785)
(9,0.46545857)
(10,0.44087617)
(11,0.40003000)
(12,0.38024984)
(13,0.38049440)
(14,0.34882146)
(15,0.33751948)
(16,0.33311863)
(17,0.32916703)
(18,0.32582483)
(19,0.32138326)
(20,0.31817974)
(21,0.31548206)
(22,0.31246056)
(23,0.30952065)
(24,0.30655277)
(25,0.30360349)
};

    \end{axis}
\end{tikzpicture}
\label{fig:ExtendedTree2-Fast-Fast-scenario}
\end{subfigure}
\hfill
\begin{subfigure}{0.475\textwidth}
\centering
\begin{tikzpicture}[scale=1.1]
      \pgfplotsset{every axis y label/.append style={rotate=0,yshift=-.4cm}}
\begin{axis}[
    title={Fast-Slow},
    x=6pt,
    xmin=0,
    xmax=26,
    ymin=0,
    ymax = 8,     ytick pos = left,
    bar width=3pt, 
    ybar stacked,
    ylabel={\#active buses (in Thousand)},    
    xlabel={years},
]
 \addplot+ [ybar, fill=gray!90,draw=black!70] coordinates {
(1,6.462)
(2,6.463)
(3,6.464)
(4,6.464)
(5,6.464)
(6,5.809)
(7,5.045)
(8,4.3)
(9,3.673)
(10,3.073)
(11,2.101)
(12,1.745)
(13,1.743)
(14,0.858)
(15,0.351)
(16,0.254)
(17,0.172)
(18,0.163)
(19,0)
(20,0)
(21,0)
(22,0.001)
(23,0.001)
(24,0.001)
(25,0)
};

 \addplot+ [ybar, fill=gray!0,draw=black!70] coordinates {
(1,0.112)
(2,0.112)
(3,0.112)
(4,0.112)
(5,0.112)
(6,0.776)
(7,1.63)
(8,2.442)
(9,3.108)
(10,3.721)
(11,4.743)
(12,5.204)
(13,5.097)
(14,6.046)
(15,6.559)
(16,6.656)
(17,6.725)
(18,6.721)
(19,6.894)
(20,6.896)
(21,6.866)
(22,6.865)
(23,6.865)
(24,6.865)
(25,6.86)
};

 \addplot+ [ybar, black] coordinates {
(1,0)
(2,0)
(3,0)
(4,0)
(5,0)
(6,0)
(7,0)
(8,0)
(9,0)
(10,0)
(11,0)
(12,0)
(13,0)
(14,0)
(15,0)
(16,0)
(17,0)
(18,0)
(19,0)
(20,0)
(21,0)
(22,0)
(23,0)
(24,0)
(25,0)
};

\end{axis}
      \pgfplotsset{every axis y label/.append style={rotate=180,yshift=7.2cm}}
    \begin{axis}[
    x=6pt,
    xmin=0,
    xmax=26,
    ymin=0,
    ymax=0.7,     ytick={0, 0.1, 0.2, 0.3, 0.4, 0.5, 0.6,0.7},
    axis y line*=right,
    axis x line=none,
    ylabel=Cost (in Billion USD),
    ]

 \addplot[mark=*,blue, mark size=1pt, only marks] coordinates {
(1,0.05048266)
(2,0.09924571)
(3,0.14745078)
(4,0.19472864)
(5,0.24061138)
(6,0.23873053)
(7,0.23668358)
(8,0.23409796)
(9,0.23189643)
(10,0.22968143)
(11,0.22727773)
(12,0.10149470)
(13,0.00143449)
(14,0.22084652)
(15,0.21900437)
(16,0.21693362)
(17,0.21491124)
(18,0.21295214)
(19,0.21112760)
(20,0.10976564)
(21,0.07773957)
(22,0.00024538)
};

 \addplot[mark=diamond,red, mark size=1pt, only marks] coordinates {
(1,0.63108967)
(2,0.62532034)
(3,0.61943815)
(4,0.61348375)
(5,0.60786673)
(6,0.57747113)
(7,0.53135996)
(8,0.49583556)
(9,0.46528435)
(10,0.44095843)
(11,0.39988835)
(12,0.38025352)
(13,0.38031450)
(14,0.34888739)
(15,0.33750348)
(16,0.33299751)
(17,0.32898354)
(18,0.32581648)
(19,0.32134678)
(20,0.31817262)
(21,0.31543778)
(22,0.31248911)
(23,0.30929356)
(24,0.30634614)
(25,0.30360374)
};

    \end{axis}
\end{tikzpicture}
\label{fig:ExtendedTree2-Fast-Slow-scenario}
\end{subfigure}

\begin{subfigure}{0.475\textwidth}
\centering
\begin{tikzpicture}[scale=1.1]
      \pgfplotsset{every axis y label/.append style={rotate=0,yshift=-.4cm}}
\begin{axis}[
    title={Slow-Fast},
    x=6pt,
    xmin=0,
    xmax=26,
    ymin = 0,
    ymax= 8,
    ytick pos = left,
    bar width=3pt,
    ybar stacked,
    ylabel={\#active buses (in Thousand)},    
    xlabel={years},
]

 \addplot+ [ybar, fill=gray!90,draw=black!70] coordinates {
(1,6.462)
(2,6.463)
(3,6.464)
(4,6.464)
(5,6.464)
(6,6.059)
(7,5.474)
(8,4.911)
(9,4.383)
(10,4.097)
(11,3.402)
(12,2.739)
(13,2.743)
(14,1.995)
(15,1.359)
(16,0.49)
(17,0.289)
(18,0.289)
(19,0.033)
(20,0.033)
(21,0.003)
(22,0.003)
(23,0.003)
(24,0.005)
(25,0)
};

 \addplot+ [ybar, fill=gray!0,draw=black!70] coordinates {
(1,0.112)
(2,0.112)
(3,0.112)
(4,0.112)
(5,0.112)
(6,0.525)
(7,1.201)
(8,1.795)
(9,2.356)
(10,2.581)
(11,3.382)
(12,4.184)
(13,4.072)
(14,4.846)
(15,5.536)
(16,6.457)
(17,6.669)
(18,6.66)
(19,6.901)
(20,6.888)
(21,6.871)
(22,6.856)
(23,6.855)
(24,6.852)
(25,6.852)
};

 \addplot+ [ybar, black] coordinates {
(1,0)
(2,0)
(3,0)
(4,0)
(5,0)
(6,0)
(7,0)
(8,0)
(9,0)
(10,0)
(11,0)
(12,0)
(13,0)
(14,0)
(15,0)
(16,0)
(17,0)
(18,0)
(19,0)
(20,0)
(21,0.021)
(22,0.021)
(23,0.021)
(24,0.021)
(25,0.021)
};

\end{axis}
      \pgfplotsset{every axis y label/.append style={rotate=180,yshift=7.2cm}}
    \begin{axis}[
    x=6pt,
    xmin=0,
    xmax=26,
    ymin=0,
    ymax=0.7,     ytick={0, 0.1, 0.2, 0.3, 0.4, 0.5, 0.6,0.7},
    axis y line*=right,
    axis x line=none,
    ylabel=Cost (in Billion USD),
    ]

 \addplot[mark=*,blue, mark size=1pt, only marks] coordinates {
(1,0.05048266)
(2,0.09924571)
(3,0.14745078)
(4,0.19472864)
(5,0.24061138)
(6,0.23836232)
(7,0.23636992)
(8,0.23408439)
(9,0.23190931)
(10,0.22951133)
(11,0.22748023)
(12,0.22567303)
(13,0.00106982)
(14,0.22113980)
(15,0.21880523)
(16,0.21672260)
(17,0.21465794)
(18,0.21262321)
(19,0.21078263)
(20,0.08684569)
(21,0.20676697)
(22,0.16981038)
(24,0.00040122)
};

 \addplot[mark=diamond,red, mark size=1pt, only marks] coordinates {
(1,0.63108967)
(2,0.62532034)
(3,0.61943815)
(4,0.61348375)
(5,0.60786673)
(6,0.58513020)
(7,0.54810991)
(8,0.51894849)
(9,0.49216061)
(10,0.48432154)
(11,0.44565711)
(12,0.41349967)
(13,0.41324508)
(14,0.38404298)
(15,0.36023786)
(16,0.33408360)
(17,0.32855334)
(18,0.32530540)
(19,0.32052892)
(20,0.31849899)
(21,0.31638498)
(22,0.31398508)
(23,0.31098906)
(24,0.30791881)
(25,0.30525809)
};

    \end{axis}
\end{tikzpicture}
\label{fig:ExtendedTree2-Slow-Fast-scenario}
\end{subfigure}
\hfill
\begin{subfigure}{0.475\textwidth}
\centering
\begin{tikzpicture}[scale=1.1]
      \pgfplotsset{every axis y label/.append style={rotate=0,yshift=-.4cm}}
\begin{axis}[
    title={Slow-Slow},
    x=6pt,
    xmin=0,
    xmax=26,
    ymin = 0,
    ymax=  8,
    ytick pos = left,
    bar width=3pt, 
    ybar stacked,
    ylabel={\#active buses (in Thousand)},    
    xlabel={years},
]

 \addplot+ [ybar, fill=gray!90,draw=black!70] coordinates {
(1,6.462)
(2,6.463)
(3,6.464)
(4,6.464)
(5,6.464)
(6,6.059)
(7,5.473)
(8,4.912)
(9,4.383)
(10,4.092)
(11,3.397)
(12,2.734)
(13,2.738)
(14,1.99)
(15,1.354)
(16,0.485)
(17,0.284)
(18,0.284)
(19,0.033)
(20,0.033)
(21,0.003)
(22,0.003)
(23,0.003)
(24,0.005)
(25,0)
};

 \addplot+ [ybar, fill=gray!0,draw=black!70] coordinates {
(1,0.112)
(2,0.112)
(3,0.112)
(4,0.112)
(5,0.112)
(6,0.525)
(7,1.201)
(8,1.794)
(9,2.356)
(10,2.587)
(11,3.388)
(12,4.19)
(13,4.078)
(14,4.852)
(15,5.545)
(16,6.466)
(17,6.677)
(18,6.668)
(19,6.903)
(20,6.89)
(21,6.907)
(22,6.892)
(23,6.891)
(24,6.888)
(25,6.888)
};

 \addplot+ [ybar, black] coordinates {
(1,0)
(2,0)
(3,0)
(4,0)
(5,0)
(6,0)
(7,0)
(8,0)
(9,0)
(10,0)
(11,0)
(12,0)
(13,0)
(14,0)
(15,0)
(16,0)
(17,0)
(18,0)
(19,0)
(20,0)
(21,0)
(22,0)
(23,0)
(24,0)
(25,0)
};

\end{axis}
      \pgfplotsset{every axis y label/.append style={rotate=180,yshift=7.2cm}}
    \begin{axis}[
    x=6pt,
    xmin=0,
    xmax=26,
    ymin=0,
    ymax=0.7,     ytick={0, 0.1, 0.2, 0.3, 0.4, 0.5, 0.6,0.7},
      axis y line*=right,
      axis x line=none,
      ylabel=Cost (in Billion USD),
    ]

 \addplot[mark=*,blue, mark size=1pt, only marks] coordinates {
(1,0.05048266)
(2,0.09924571)
(3,0.14745078)
(4,0.19472864)
(5,0.24061138)
(6,0.23836232)
(7,0.23636992)
(8,0.23420926)
(9,0.23178562)
(10,0.22957156)
(11,0.22748023)
(12,0.22567303)
(13,0.00106982)
(14,0.22113980)
(15,0.21927134)
(16,0.21672260)
(17,0.21485064)
(18,0.21262321)
(19,0.21080974)
(20,0.08841331)
(21,0.20672867)
(22,0.17043186)
(24,0.00040122)
};

 \addplot[mark=diamond,red, mark size=1pt, only marks] coordinates {
(1,0.63108967)
(2,0.62532034)
(3,0.61943815)
(4,0.61348375)
(5,0.60786673)
(6,0.58511771)
(7,0.54810991)
(8,0.51896437)
(9,0.49221031)
(10,0.48404599)
(11,0.44548858)
(12,0.41326697)
(13,0.41312438)
(14,0.38372603)
(15,0.36023064)
(16,0.33399242)
(17,0.32841749)
(18,0.32519154)
(19,0.32051702)
(20,0.31848325)
(21,0.31603412)
(22,0.31386026)
(23,0.31082324)
(24,0.30781231)
(25,0.30518164)
};

    \end{axis}
\end{tikzpicture}
\label{fig:ExtendedTree2-Slow-Slow-scenario}
\end{subfigure}

\caption{Extended Scenario Tree - 2 By 3}

\label{fig:ExtendedTree2-4scenarios-adj}

\end{figure}
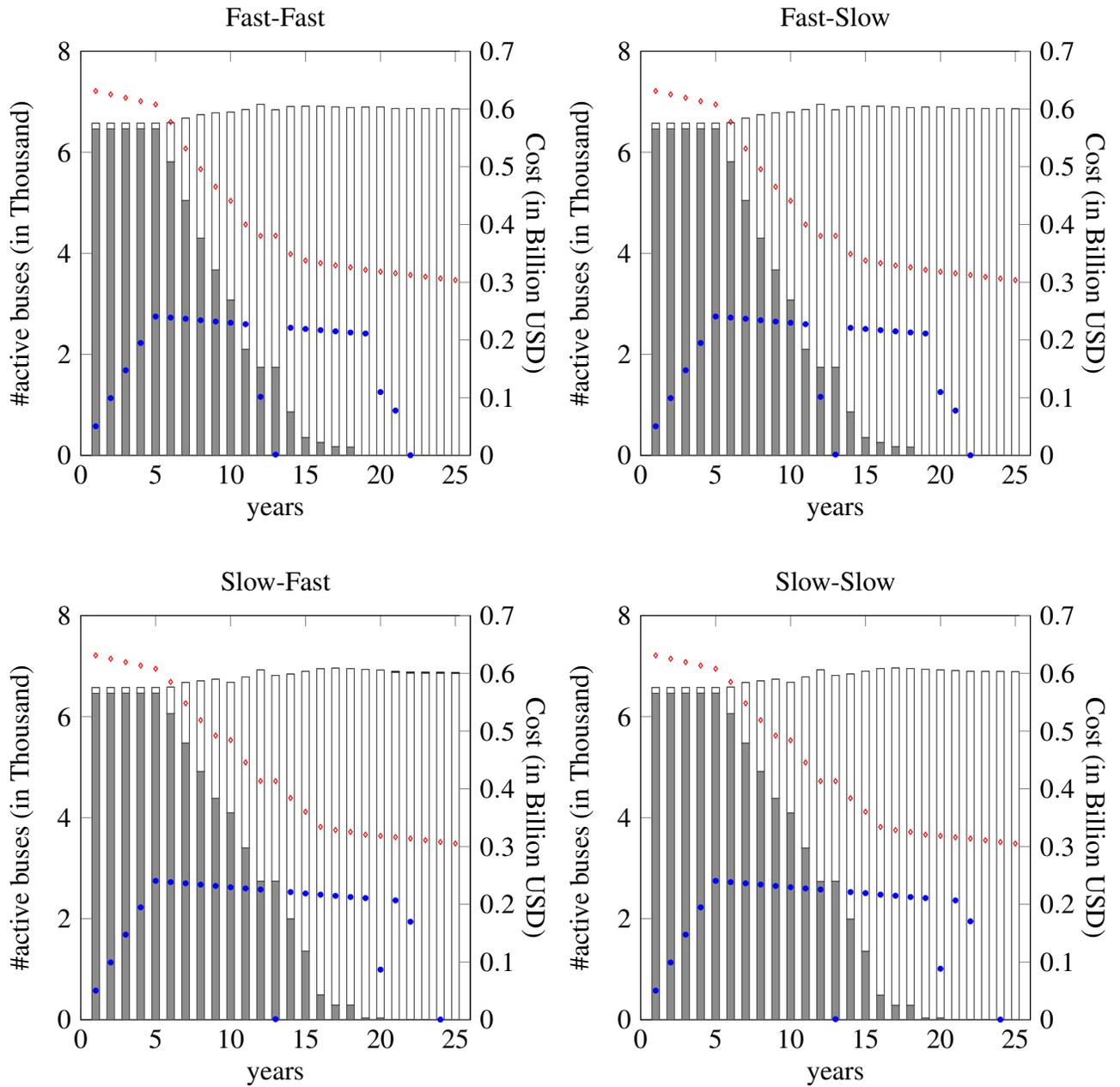

We present the results of the Extended Scenario Tree - 2 By 3 case in Figure~\ref{fig:ExtendedTree2-4scenarios-adj}. 
In this case, two branches for BEB technological improvements and three branches for HFCB technological improvements are considered in each stage. We used simplified assignment decisions, similar to the simplified base case (Section~\ref{sec:base-simplified}).

CPU time is 26,039 seconds and the total expected cost is 14,111,360,587 USD, which increases to 14,137,237,734 USD after rounding the variables.

\end{document}